\documentclass[letterpaper]{article}
\usepackage[english]{babel}

\usepackage{authblk}

\usepackage{amsmath,amsthm}
\usepackage{amssymb}
\usepackage{float}
\usepackage{enumitem}
\usepackage{empheq}

\makeatletter

\newskip\@bigflushglue \@bigflushglue = -100pt plus 1fil

\def\bigcentering{\let\\\@centercr\rightskip\@bigflushglue%
\leftskip\@bigflushglue
\parindent\z@\parfillskip\z@skip}

\makeatother

\usepackage{graphics,graphicx,xcolor}
\usepackage{standalone}
\usepackage{tikz}
\usetikzlibrary{calc,shapes.arrows,decorations.pathmorphing,decorations.pathreplacing,decorations.markings,arrows,positioning,patterns,shapes.multipart}

\definecolor{rouge}{RGB}{255,77,77}
\definecolor{vert}{RGB}{0,178,102}
\definecolor{jaune}{RGB}{255,255,0}
\definecolor{violet}{RGB}{208,32,144}
\definecolor{orange}{RGB}{255,140,0}
\definecolor{bleu}{RGB}{0,0,205}

\theoremstyle{plain}

\newtheorem*{theorem*}{Theorem}
\newtheorem{theorem}{Theorem}[section]
\newtheorem{lemma}[theorem]{Lemma}
\newtheorem{proposition}[theorem]{Proposition}

\newtheorem*{questions*}{Questions}
\newtheorem{example}[theorem]{Example}
\newtheorem{remark}[theorem]{Remark}

\def\A{\mathcal{A}}

\def\F{\mathcal{F}}
\def\H{\mathbb{H}}
\def\N{\mathbb{N}}

\def\R{\mathbb{R}}
\def\Z{\mathbb{Z}}

\newcommand{\define}[1]{\emph{#1}}

\newcommand{\lettertone}{\vbox to 12pt{\hbox to 3pt{
\begin{tikzpicture}[scale=0.3]
\draw (0,0) rectangle (2,1);
\draw[color=bleu,pattern color=bleu,pattern=north west lines] (0,0) rectangle (0.5,0.5);
\draw[color=bleu,pattern color=bleu,pattern=north west lines] (1.5,0) rectangle (2,0.5);
\end{tikzpicture}
}}}

\newcommand{\letterttwo}{\vbox to 12pt{\hbox to 3pt{
\begin{tikzpicture}[scale=0.3]
\draw (0,0) rectangle (2,1);
\draw[color=bleu,pattern color=bleu,pattern=north west lines] (0.5,0) rectangle (2,0.5); 
\end{tikzpicture}
}}}

\newcommand{\lettertthree}{\vbox to 12pt{\hbox to 3pt{
\begin{tikzpicture}[scale=0.3]
\draw (0,0) rectangle (2,1);
\draw[color=bleu,pattern color=bleu,pattern=north west lines] (0,0) rectangle (2,0.5); 
\end{tikzpicture}
}}}

\newcommand{\lettertfour}{\vbox to 12pt{\hbox to 3pt{
\begin{tikzpicture}[scale=0.3]
\draw (0,0) rectangle (2,1);
\draw[color=bleu,pattern color=bleu,pattern=north west lines] (0,0) rectangle (1.5,0.5); 
\end{tikzpicture}
}}}

\newcommand{\letterbone}{\vbox to 12pt{\hbox to 3pt{
\begin{tikzpicture}[scale=0.3]
\draw (0,0) rectangle (2,1);
\draw[color=bleu,pattern color=bleu,pattern=north west lines] (1,0.5) rectangle (2,1); 
\end{tikzpicture}
}}}

\newcommand{\letterbtwo}{\vbox to 12pt{\hbox to 3pt{
\begin{tikzpicture}[scale=0.3]
\draw (0,0) rectangle (2,1);
\draw[color=bleu,pattern color=bleu,pattern=north west lines] (0,0.5) rectangle (2,1);
\end{tikzpicture}
}}}

\newcommand{\letterbthree}{\vbox to 12pt{\hbox to 3pt{
\begin{tikzpicture}[scale=0.3]
\draw (0,0) rectangle (2,1);
\draw[color=bleu,pattern color=bleu,pattern=north west lines] (0,0.5) rectangle (1,1); 
\end{tikzpicture}
}}}

\newcommand{\lettermone}{\vbox to 12pt{\hbox to 3pt{
\begin{tikzpicture}[scale=0.3]
\draw (0,0) rectangle (2,1);
\draw[color=bleu,pattern color=bleu,pattern=north west lines] (0.5,0) -- (2,0) -- (2,1) -- (1,1) -- cycle ;
\end{tikzpicture}
}}}

\newcommand{\lettermtwo}{\vbox to 12pt{\hbox to 3pt{
\begin{tikzpicture}[scale=0.3]
\draw (0,0) rectangle (2,1);
\draw[color=bleu,pattern color=bleu,pattern=north west lines] (1.5,0) -- (2,0) -- (2,1) -- (1,1) -- cycle ; 
\draw[color=bleu,pattern color=bleu,pattern=north west lines] (0,0) -- (0.5,0) -- (0,0.5) -- cycle ; 
\end{tikzpicture}
}}}

\newcommand{\lettermthree}{\vbox to 12pt{\hbox to 3pt{
\begin{tikzpicture}[scale=0.3]
\draw (0,0) rectangle (2,1);
\draw[color=bleu,pattern color=bleu,pattern=north west lines] (0,0) rectangle (2,1); 
\end{tikzpicture}
}}}

\newcommand{\lettermfour}{\vbox to 12pt{\hbox to 3pt{
\begin{tikzpicture}[scale=0.3]
\draw (0,0) rectangle (2,1);
\draw[color=bleu,pattern color=bleu,pattern=north west lines] (0,0) -- (1.5,0) -- (1,1) -- (0,1) -- cycle ; 
\end{tikzpicture}
}}}

\newcommand{\lettermfive}{\vbox to 12pt{\hbox to 3pt{
\begin{tikzpicture}[scale=0.3]
\draw (0,0) rectangle (2,1);
\draw[color=bleu,pattern color=bleu,pattern=north west lines] (0,0) -- (2,0) -- (2,0.5) -- (1,1) -- (0,1) -- cycle ;
\end{tikzpicture}
}}}

\newcommand{\tile}[4]{
\draw[thick,fill=#4] (#1,#2) -- (#1,#2+3) -- (#1+#3*5,#2+3) -- (#1+#3*5,#2) -- (#1+#3*4,#2) -- (#1+#3*3,#2) -- (#1+#3*2,#2) -- (#1+#3*1,#2) -- (#1,#2) -- cycle ;
}

\title{Tilings of the hyperbolic plane of substitutive origin as subshifts of finite type on Baumslag-Solitar groups $BS(1,n)$.}
\date{}
\author[1]{Nathalie Aubrun}
\author[2]{Michael Schraudner}
\affil[1]{LRI, Univ. Paris-Sud, CNRS, INRIA, Universit\'e Paris-Saclay, F-91405 Orsay, France}
\affil[2]{Centro de Modelamiento Matem\'atico, Universidad de Chile, Av. Blanco Encalada 2120, Piso 7, Santiago de Chile}
\begin{document}
 
\maketitle 
 
\begin{abstract}
We present a technique to lift some tilings of the discrete hyperbolic plane --tilings defined by a 1D substitution-- into a zero entropy subshift of finite type (SFT) on non-abelian amenable Baumslag-Solitar groups $BS(1,n)$ for $n\geq2$. For well chosen hyperbolic tilings, this SFT is also aperiodic and minimal. As an application we construct a strongly aperiodic SFT on $BS(1,n)$ with a hierarchical structure, which is an analogue of Robinson's construction on~$\Z^2$ or Goodman-Strauss's on~$\H_2$.
\end{abstract}

\section*{Introduction}
\label{section.introduction}

Given a finitely generated group $G$ and a finite alphabet $\A$, the set $\A^G$ can be endowed with the pro discrete topology. The group $G$ acts --from the left-- by translation on $\A^G$, and a classical result states that with this topology, the set $\A^G$ is a compact set. Elements of $\A^G$ are called configurations, and can be thought of as colorings of the group $G$ with colors chosen among the finite set~$\A$. Symbolic dynamics on $G$ studies subshifts, i.e. sets of configurations of $G$ that are both invariant under the action of $G$, and closed for the product topology. Equivalently, subshifts can be defined as sets of configurations that avoid some set of finite patterns --this set of forbidden patterns not necessarily being finite. If the set of forbidden patterns can be chosen finite, the subshift is called a subshift of finite type (SFT for short). SFTs are of particular interest for two reasons. First, because they can be used to model real-life phenomena defined through local interactions. Second, because their finite description allows to ask complexity and computability questions

\medskip

Among questions related to SFTs, an open problem is to characterize groups $G$ which admit a strongly aperiodic SFT, i.e. an SFT such that all its configurations have a trivial stabilizer. The same question for general subshifts --non SFT-- has been answered positively~\cite{GaoJacksonSeward2009,AubrunBarbieriThomasse2018}. Free groups cannot possess strongly aperiodic SFTs~\cite{MullerSchupp1985}, and a finitely generated and recursively presented group with an aperiodic SFT necessarily has a decidable word problem~\cite{Jeandel2015aperiodic}. Groups that are known to admit a strongly aperiodic SFT are $\Z^2$~\cite{Robinson1971} and $\Z^d$ for $d>2$~\cite{CulikKari1996}, fundamental groups of oriented surfaces~\cite{CohenGoodmanStrauss2017}, hyperbolic groups~\cite{CohenGSRieck2017}, discrete Heisenberg group~\cite{SahinSchraudnerUgarcovici2020} and more generally groups that can be written as a semi-direct product $G=\Z^2\rtimes_\phi H$, provided $G$ has decidable word problem~\cite{BarbieriSablik2019}, and amenable Baumslag-Solitar groups~\cite{EsnayMoutot2020}.

\medskip

There are several techniques to produce aperiodicity inside tilings, but in the particular setting of SFTs, two of them stand out. The first goes back to Robinson~\cite{Robinson1971}: the aperiodicity of Robinson's SFT follows from the hierarchical structure shared by all configurations. Indeed, in every configuration one sees different levels of squares, where all squares of a given level are arranged according to a lattice, and are gathered by four to produce a square of higher level. Since there are arbitrarily large squares, no non-trivial translation can keep the configuration unchanged. The second technique goes back to Kari~\cite{Kari1996}: a very simple aperiodic dynamical system is encoded into an SFT, so that configurations correspond to orbits of the dynamical system. Aperiodicity of the SFT comes from both the aperiodicity of the dynamical system and the clever encoding.

\medskip

In this article we focus on the particular case where $G$ is a Baumslag-Solitar group. Baumslag-Solitar groups are examples of HNN extensions, a fundamental construction in combinatorial group theory that plays for instance a key role in the proof of Higman embedding theorem~\cite{LyndonSchupp2001}. There are described by the two generators and one relation presentation
\[
 BS(m,n):=\langle a,t\mid t^{-1}a^mt=a^n\rangle,
\]
where $m,n$ are two positive integers. The whole class of $BS(m,n)$ is clearly separated into two subclasses with radically different behaviors: the groups $BS(1,n)$ are solvable 
hence amenable, while the $BS(m,n)$ groups with $m,n>1$ contain free subgroups 
and are consequently non amenable. Each of these two subclasses is also well-understood from a geometrical point of view, since it is known under which condition two $BS(m,n)$ are quasi-isometric:
\begin{itemize}
 \item groups $BS(1,n)$ and $BS(1,n')$ are quasi-isometric if and only if $n$ and $n'$ have a common power~\cite{FarbMosher1998} --and in this case, the two groups are even commensurable;
 \item groups $BS(m,n)$ and $BS(m',n')$ are quasi-isometric as soon as $2\leq m < n$ and $2\leq m'<n'$~\cite{Whyte2001}.
\end{itemize}
In~\cite{AubrunKari2013} the example of aperiodic SFT that is given is a construction adapted from~\cite{Kari1996}. This SFT is proved to be weakly aperiodic --all configurations have infinite orbit-- in~\cite{AubrunKari2013} for all Baumslag-Solitar groups, and not strongly aperiodic for the non amenable ones. In~\cite{EsnayMoutot2020} the authors show that this SFT is actually strongly aperiodic for amenable Baumslag-Solitar groups.

\medskip

In this article, we will only consider non abelian amenable Baumslag-Solitar groups $BS(1,n)$ with $n\geq2$. Our main result consists in a construction that allows to lift some particular tilings of the discrete hyperbolic plane --tilings defined from a 1D substitution-- into an SFT on $BS(1,n)$ (Theorem~\ref{theorem.sigma_tilings_and_X_sigma}). If the substitution is well chosen, then the SFT is moreover minimal, strongly aperiodic and has zero entropy (Propositions~\ref{proposition.X_sigma_minimal},~\ref{proposition.X_sigma_aperiodic} and~\ref{proposition.X_zero_entropy}). As an application, we construct a strongly aperiodic SFT on $BS(1,n)$ with a hierarchical structure, that goes back to Robinson's construction on $\Z^2$~\cite{Robinson1971} and Goodman-Strauss's construction on the discrete hyperbolic plane~\cite{GS2010}. We strongly believe that this hierarchical SFT will also serve as a basis for more sophisticated constructions.

\section{Subshifts on amenable Baumslag-Solitar groups}
\label{section.subshifts_on_BS}

\subsection{Subshifts on finitely generated groups}
\label{subsection.1.1}

In this article $\A$ denotes a finite set, and $G$ is a finitely generated group. The identity of $G$ is denoted $1_G$. Endowed with the prodiscrete topology, the set $\A^G$ is a compact and metrizable space. Elements of $\A^G$ are called configurations, and can be though of as colorings of the group $G$ by the finite alphabet $\A$. The shift map is the natural left action of $G$ on $\A^G$ by translation 
$$\mathfrak{S}:
 \left(
 \begin{array}{lll}
  G\times \A^G & \to & \A^G\\
  (g,x) & \mapsto & \mathfrak{S}_g(x)
 \end{array}
 \right)
 $$
where $\mathfrak{S}_g(x)$ is the configuration such that $\left(\mathfrak{S}_g(x)\right)_h = x_{g^{-1}\cdot h}$ for every $h\in G$. The dynamical system $\left(\A^G,\mathfrak{S}\right)$ is called a full-shift. Subshifts are subsystems of the full-shift, i.e. subsets of $\A^G$ that are both closed for the prodiscrete topology and invariant under the shift action. Interestingly this dynamical definition of subshifts coincides with a combinatorial one. A pattern if a finite configuration $p\in\A^S$, where $S$ is a finite subset of $G$ called the support of $p$. A pattern $p\in\A^S$ appears in a configuration $x\in \A^G$ if there exists a group element $g\in G$ such that $p_h=\left(\mathfrak{S}_g(x)\right)_h$ for every $h\in S$. Otherwise we say that $x$ avoids $p$. Given a set of patterns $\F$, define $X_\F$ as the set of configurations that avoid all patterns in $\F$. A set of configurations $Y$ is a subshift if and only if there exists some $\F$ such that $Y=X_\F$. A subshift of finite type --SFT for short-- is a subshift for which the set of forbidden patterns $\F$ can be chosen finite. A sofic subshift is a subshift which is the image of an SFT under a continuous and shift commuting map --such maps are called morphisms.

\medskip

Let $x\in \A^G$ be a configuration. The orbit of $x$ is the set of configurations $orb_\mathfrak{S}(x)=\left\{\mathfrak{S}_g(x)\mid g\in G \right\}$, and the stabilizer of $x$ is the subgroup $stab_\mathfrak{S}(x)=\left\{\ g\in G\mid \mathfrak{S}_g(x)=x\right\}$. 
A subshift $X\subseteq \A^G$ is weakly aperiodic if for every configuration $x\in X$, $|orb_\mathfrak{S}(x)|=\infty$. A subshift $X\subseteq \A^G$ is strongly aperiodic if for every configuration $x\in X$, $stab_\mathfrak{S}(x)=\left\{ 1_G\right\}$. The subshift $X$ is minimal if it does not contain a proper subshift or, equivalently, if the orbit of each configuration $x\in X$ is dense.

\subsection{Baumslag-Solitar groups}
\label{subsection.BS}

Given two positive integers $m$ and $n$, define the finitely presented group
$$BS(m,n):=\langle a,t\mid t^{-1}a^mt=a^n\rangle,$$
known as the Baumslag-Solitar group with parameters $m,n$. In this article we are only interested in the case $m=1$, which makes these groups solvable hence amenable. For the class of amenable Baumslag-Solitar groups, there exists a short normal form, as stated below in Proposition~\ref{proposition.normal_form}. 


\begin{proposition}\label{proposition.normal_form}
Every element of $BS(1,n)$ can be uniquely written as one of the following words 
\begin{itemize}
 \item $a^jt^{-k}$ with $k\geq0$ and $j\in\Z$,
 \item or $t^ia^j$ with $i>0$ and $j\in\Z$,
 \item or $t^ia^jt^{-k}$ with $i>0$, $k>0$ and $n\nmid j$.
\end{itemize}
\end{proposition}

\begin{proof}
Using the relation $at=ta^n$, one deduces that $t^{-1}a=a^nt^{-1}$. Combining these two rewriting rules, any word on $\{a,a^{-1},t,t^{-1}\}$ can be rewritten as some $t^ia^jt^{-k}$ with $i,k\geq0$ and $j\in\Z$. If $i=0$ or $k=0$ then we are done. Otherwise if $i,k>0$ and $j=nj'$, write $t^ia^jt^{-k}=t^{i-1}a^{j'}t^{-k+1}$ and iterate the process as soon as $i,k>0$ and $n|j$. 

It remains to prove uniqueness. Note that for a given $g$, all words $w$ in $a,a^{-1},t,t^{-1}$ that represent $g$ should share the same quantity $|w|_{t}-|w|_{t^{-1}}$. So if  $g$ can be written as 
$w=t^ia^jt^{-k}$ and $w'=t^{i'}a^{j'}t^{-k'}$ with $i,i',k,k'\geq0$, and $j,j'\in\Z$, then either $i=i'=0$ and $k=k'$, or $k=k'=0$ and $i=i'$, or $i,k,i',k'>0$ and $n\nmid j$, $n\nmid j'$ and $i-k=i'-k'$. In the two first cases, we have directly $i=i'$, $j=j'$ and $k=k'$. In the third case, assume that $i'=i+\ell$ and $k'=k+\ell$ with $\ell\geq0$. Then necessarily $t^\ell a^{j'}t^{-\ell}=a^j$ with $n\nmid j$ and $n\nmid j'$. This is not possible unless $\ell=0$, so we finally get $i=i'$, $j=j'$ and $k=k'$, and the normal form is unique.
\end{proof}

Remark that this normal form does not have minimal length: for instance the element $t^2a^{n^2-1}$ could have been written using less generators as $a^{-1}t^2a$. 

\medskip

Consider the semi-direct product $\Z[\frac{1}{n}]\rtimes_f\Z$ given by
$$
f:
 \left(
 \begin{array}{lll}
  \Z & \to & Aut\left(\Z[\frac{1}{n}]\right)\\
  k & \mapsto & x\mapsto n^k\cdot x
 \end{array}
 \right).
$$

Using the normal form of Proposition~\ref{proposition.normal_form}, we prove that $BS(1,n)$ and $\Z[\frac{1}{n}]\rtimes_f\Z$ are isomorphic, through

$$
\Phi:
 \left(
 \begin{array}{lll}
  BS(1,n) & \to & \Z[\frac{1}{n}]\rtimes_f\Z\\
  t^{i}a^{j}t^{-k} & \mapsto & \left(j\cdot n^{-i},k-i\right)
 \end{array}
 \right).
$$

\begin{remark}\label{remark.phi_well_defined}
The application $\Phi$ is well-defined. Suppose an element $g\in BS(1,n)$ can be written as $w=t^{i}a^{n^{\alpha}j}t^{-k}$ for some $\alpha\geq1$ and $n\nmid j$. Using that $a^{n^{\alpha}j}=t^{-\alpha}a^j t^{\alpha}$, we deduce the word $w'=t^{i-\alpha}a^{j}t^{-k+\alpha}$ also represents the element $g$. And since
\begin{align*}
\Phi(t^{i}a^{n^{\alpha}j}t^{-k}) &=\left(n^{\alpha}j\cdot n^{-i},k-i\right) \\
 &=\left(j\cdot n^{-(i-\alpha)},(k-\alpha)-(i-\alpha)\right) \\
 &= \Phi(t^{i-\alpha}a^{j}t^{-(k-\alpha)}) \\
 &= \Phi(g),
\end{align*}
so that $\Phi$ is well-defined.
\end{remark}

\begin{proposition}\label{proposition.phi_isomorphism}
 The map $\Phi$ is an isomorphism.
\end{proposition}

\begin{proof}
We first check that $\Phi$ is a group morphism. Take $g,h$ two group elements in $G$, and assume they can be written as $g=t^{i}a^jt^{-k}$ and $h=t^{i'}a^{j'}t^{-k'}$. Remind that the law for the semi-direct product is $(h_1,k_1)\cdot(h_2,k_2)=\left(h_1+f(k_1)(h2),k_1+k_2\right)$ for $(h_i,k_i)\in\Z[\frac{1}{n}]\rtimes_f\Z$. Then on the one hand
\begin{align*}
\Phi(gh) &=\Phi\left(t^{i}a^jt^{-k}\cdot t^{i'}a^{j'}t^{-k'}\right) \\
 &=\Phi\left(t^{i}a^jt^{i'}\cdot t^{-k}a^{j'}t^{-k'}\right) \\ 
 &=\Phi\left(t^{i+i'}a^{jn^{i'}}\cdot a^{j'n^{k}}t^{-k-k'}\right) \\
 &=\left((jn^{i'}+j'n^{k})n^{-(i+i')},k+k'-i-i'\right)
\end{align*}
where the last equality comes from Remark~\ref{remark.phi_well_defined}, and on the other hand
\begin{align*}
\Phi(g)\Phi(h) &=\left(j\cdot n^{-i},k-i\right)\cdot\left(j'\cdot n^{-i'},k'-i'\right) \\
 &=\left(j\cdot n^{-i}+n^{k-i}\cdot j' n^{-i'},k+k'-i-i'\right) \\
 &=\Phi(gh),
\end{align*}
which proves that $\Phi$ is a group morphism.

Assume $g=t^{i}a^{nj+r}t^{-k}$ and $h=t^{i'}a^{nj'+r'}t^{-k'}$ and suppose $\Phi(g)=\Phi(h)$ with $0\leq r,r'<n$ --$r\neq0$ if $i,j>0$ and $r'\neq0$ if $i',j'>0$. Then 
$$\left( (nj+r)\cdot n^{-i},k-i \right)=\left( (nj'+r')\cdot n^{-i'},k'-i' \right),$$
which directly leads to $i=i'$, $j=j'$ and $k=k'$, hence $g=h$.
\end{proof}

\subsection{Structure of $BS(1,n)$}
\label{subsection.structure_BS}

The Cayley graph of $BS(1,n)$ with generating set $\{a,t,a^{-1},t^{-1}\}$ is made of several sheets that merge $n$ by $n$, so that the global structure these sheets are arranged looks like an $(n+1)$-regular tree. Each of these sheet is quasi-isometric to the hyperbolic plane $\mathbb{H}^2$. From now on, we only consider the case $n=2$, but the attentive reader will be convinced that the rest of the article easily transpose to the general case $n>2$
.
\medskip

We call \define{rectangle} any finite subset of $BS(1,2)$ of the form
$$R_{k,\ell}:=\left\{ t^{\ell}a^it^{-j} \mid i\in [0;(k+1)\cdot2^{\ell-1}-1]\text{ and }j\in[0;\ell]\right\}.$$

  \begin{figure}[!ht]
  \begin{bigcenter}
  \begin{tikzpicture}[scale=0.25]

\begin{scope}[rotate=20]
 
\begin{scope}[color=bleu!50]
\foreach \i [count = \xi] in {0,...,46}{
	\node[scale=0.75] at (\i,0) {\textbullet};
	\draw (\i,0) -- (\i+1,0);
	}
\node[scale=0.75] at (47,0) {\textbullet};	
\foreach \i [count = \xi] in {0,...,22}{
	\node[scale=0.75] at (2*\i,5) {\textbullet};
	\draw (2*\i,5) -- (2*\i+2,5);
	\draw (2*\i,5) -- (2*\i,0);
	}	
\node[scale=0.75] at (2*23,5) {\textbullet};
\draw (2*23,5) -- (2*23,0);
\foreach \i [count = \xi] in {0,...,10}{
	\node[scale=0.75] at (4*\i,10) {\textbullet};
	\draw (4*\i,10) -- (4*\i+4,10);
	\draw (4*\i,10) -- (4*\i,5);
	}
\node[scale=0.75] at (4*11,10) {\textbullet};
\draw (4*11,10) -- (4*11,5);
\foreach \i [count = \xi] in {0,...,4}{
	\node[scale=0.75] at (8*\i,15) {\textbullet};
	\draw (8*\i,15) -- (8*\i+8,15);
	\draw (8*\i,15) -- (8*\i,10);
	}
\node[scale=0.75] at (8*5,15) {\textbullet};
\draw (8*5,15) -- (8*5,10);
\foreach \i [count = \xi] in {0,...,1}{
	\node[scale=0.75] at (16*\i,20) {\textbullet};
	\draw[very thick] (16*\i,20) -- (16*\i+16,20);
	\draw (16*\i,20) -- (16*\i,15);
	}

\node[scale=0.75] at (16*2,20) {\textbullet};
\draw (16*2,20) -- (16*2,15);

\begin{scope}[decoration={
    markings,
    mark=at position 0.5 with {\arrow{latex}}}
    ] 
\draw[postaction={decorate}] (0,20) -- (0,15);
\draw node at (-1,17.5) {$t$};
\draw[postaction={decorate}] (16,20) -- (16,15);
\draw node at (15,17.5) {$t$};
\draw[postaction={decorate}] (0,15) -- (8,15);
\draw node at (4,16) {$a$};
\draw[postaction={decorate}] (8,15) -- (16,15);
\draw node at (12,16) {$a$};
\draw[postaction={decorate}] (0,20) -- (16,20);
\draw node at (8,21) {$a$};
\end{scope}

\end{scope}



\begin{scope}[color=rouge]
\foreach \i [count = \xi] in {0,...,22}{
	\node[scale=0.75] at (2*\i+1+2,4) {\textbullet};
	\draw (2*\i+1+2,4) -- (2*\i+2+1+2,4);
	\draw (2*\i+1+2,4) -- (2*\i+1,0);
	}
\node[scale=0.75] at (2*23+1+2,4) {\textbullet};
\draw (2*23+1+2,4) -- (2*23+1,0);
\foreach \i [count = \xi] in {0,...,10}{
	\node[scale=0.75] at (4*\i+2+1+2+2,8) {\textbullet};
	\draw (4*\i+2+1+2+2,8) -- (4*\i+4+2+1+2+2,8);
	\draw (4*\i+2+1+2+2,8) -- (4*\i+2+1+2,4);
	}
\node[scale=0.75] at (4*11+2+1+2+2,8) {\textbullet};
\draw (4*11+2+1+2+2,8) -- (4*11+2+1+2,4);
\foreach \i [count = \xi] in {0,...,4}{
	\node[scale=0.75] at (8*\i+4+2+1+2+4,12) {\textbullet};
	\draw (8*\i+4+2+1+2+4,12) -- (8*\i+8+4+2+1+2+4,12);
	\draw (8*\i+4+2+1+2+4,12) -- (8*\i+4+2+1+4,8);
	}
\node[scale=0.75] at (8*5+4+2+1+2+4,12) {\textbullet};
\draw (8*5+4+2+1+2+4,12) -- (8*5+4+2+1+4,8);
\foreach \i [count = \xi] in {0,...,1}{
	\node[scale=0.75] at (16*\i+8+4+2+1+8,16) {\textbullet};
	\draw[very thick] (16*\i+8+4+2+1+8,16) -- (16*\i+16+8+4+2+1+8,16);
	\draw (16*\i+8+4+2+1+8,16) -- (16*\i+8+4+2+1+2+4,12);
	}
\node[scale=0.75] at (16*2+8+4+2+1+8,16) {\textbullet};
\draw (16*2+8+4+2+1+8,16) -- (16*2+8+4+2+1+2+4,12);
\end{scope}

\end{scope}

\end{tikzpicture}
  \end{bigcenter}
  \caption{The rectangle $R_{2,4}$ in $BS(1,2)$. The leftmost sheet is pictured in blue, the rightmost in red. The $2^{4}-2=14$ sheets in between are not pictured.}
  \label{figure.rectangle}
  \end{figure}
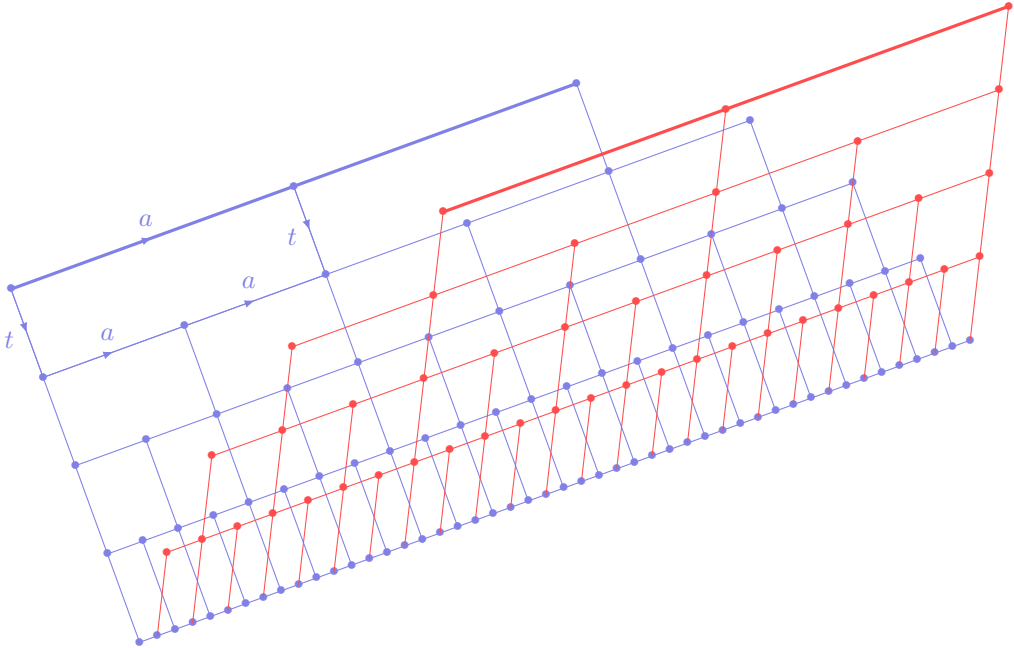

%
%

%
%
%
%

These rectangles are a classical example of a F\o lner sequence for $BS(1,2)$. As so, they will be used in the next section to define the topological entropy for subshifts on $BS(1,2)$.

\begin{proposition}\label{proposition.rectangles_Folner}
 The sequence $\left(R_{k,\ell}\right)_{k,\ell\in\N}$ is a F\o lner sequence.
\end{proposition}
\subsection{Entropy}
\label{subsection.entropy_BS}

The \define{language of size $(k,\ell)$} of a subshift $X\subseteq \A^{BS(1,2)}$ is the set of patterns with support $R_{k,\ell}$ that appear in some configuration of $X$
$$\mathcal{L}_{k,\ell}(X):=\left\{p\in \A^{R_{k,\ell}}\mid \exists x\in X,g\in BS(1,2)\text{ s.t. }\forall h \in R_{k,\ell} ,x_{gh} = p_h \right\}.$$
The \define{language} of a subshift $X\subseteq \A^{BS(1,2)}$ is given by
$$\mathcal{L}(X):=\bigcup_{k,\ell\in\N}\mathcal{L}_{k,\ell}(X),$$
and is the set of globally admissible patterns; i.e. finite patterns that can be extended into a valid configuration.

The \define{entropy} of a subshift $X\subseteq \A^{BS(1,2)}$ is given by
$$h(X):=\lim_{k,\ell\to\infty} \frac{\log\left(\sharp \mathcal{L}_{k,\ell}(X)\right)}{|R_{k,\ell}|}.$$

Since $BS(1,2)$ is amenable, this limit exists and does not depend on the choice of the F\o lner
sequence~\cite{OrnsteinWeiss1987}. Moreover, we have that 
$$h(X)=\limsup_{k,\ell\to\infty} \frac{\log\left(\sharp \mathcal{L}_{k,\ell}(X)\right)}{|R_{k,\ell}|}.$$

%
%
%
%
%
%
%

\section{Substitutions and their tilings}
\label{section.substitutions}

A --non-deterministic-- \define{substitution} is $\sigma=(\A,R)$, with $\A=\left\{a_1,\dots,a_n\right\}$ a finite alphabet of size $n$ and $R\subset\A\times\A^*$ the set of rules of $\sigma$. A substitution $\sigma$ is \define{deterministic} if for every letter $a\in\A$, there exists only one rule $(a,\sigma(a))\in R$. In the deterministic setting, we associate to a substitution $\sigma$ its \define{incidence matrix} $M_\sigma\in \mathcal{M}_{n}(\N)$ given by $M_\sigma(i,j)=|\sigma(a_i)|_{a_j}$ the number of occurrences of $a_j$ inside the word $\sigma(a_i)$. A substitution $\sigma$ is \define{primitive} if there exists a power of $M_\sigma$ for which all the entries are positive. Perron Frobenius theorem implies that the incidence matrix of a primitive substitution admits an \define{expanding eigenvalue}, i.e. an eigenvalue $\lambda$ which is real, greater than one and greater than the modulus of all the other eigenvalues. Similarly, and following~\cite{ABM2019}, we say that a non-deterministic substitution $\sigma=(\A,R)$ has an expanding eigenvalue if there exists a real $\lambda>1$ and a vector $v:\A\to\R^{+*}$ such that, for every $(a,w)\in R$ :
    \[
    \lambda\cdot v(a)=\sum_{i=1}^{|w|} v(w_i).
    \]
In the sequel we treat only the case of deterministic substitutions in order to simplify notations, but all results presented still hold for non-deterministic substitutions.

\begin{example} The deterministic substitution given by $\sigma(a)=aab$ and $\sigma(b)=ba$ has an expanding eigenvalue $\lambda =\frac{3+\sqrt{5}}{2}$ and associated eigenvector $v=\left( \frac{1+\sqrt{5}}{2},1\right)$.\end{example}

Let $\sigma:\A\to\A$ be a primitive substitution with an expanding eigenvalue $\lambda>1$. The \define{substitutive subshift} $Z_\sigma\subset \A^\Z$ is a two-sided one-dimensional subshift defined by
$$
Z_\sigma:=\left\{z\in\A^\Z\mid \forall w\sqsubset z,\;\exists a\in\A\text{ and }n\in\N\text{ s.t. }w\sqsubset \sigma^n(w)\right\}.
$$

We now describe orbits of substitutions as tilings of $\R^2$, as previously done in~\cite{ABM2019}. A \define{tile} is a compact subset of $\R^2$ with non-empty interior. If $T$ is a set of tiles, not necessarily finite, a \define{$T$-tiling} or \define{tiling by $T$} is a collection of translated copies of tiles which have pairwise disjoint interiors and whose union is the entire $\R^2$.

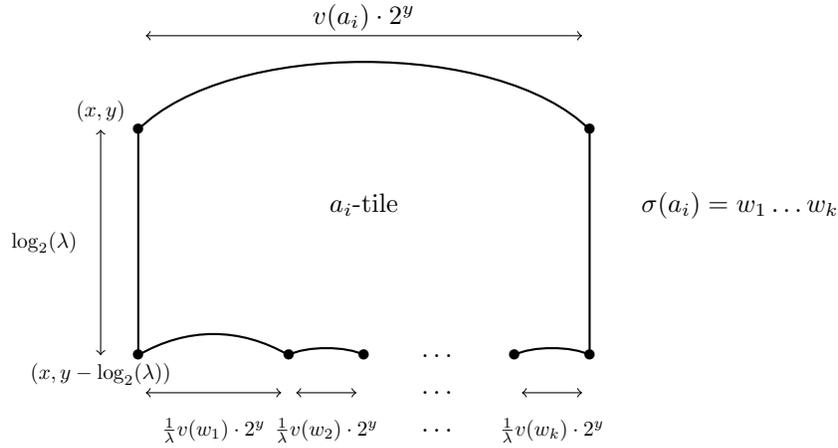
\begin{figure}[!ht]
\begin{bigcenter}
\begin{tikzpicture}[scale=1]
\draw[thick] (0,0) -- (0,3);
\draw[thick] (0,3) to [controls=+(45:1.7) and +(135:1.7)] (6,3);
\draw[thick] (6,3) -- (6,0);
\draw[thick] (0,0) to [controls=+(30:0.75) and +(150:0.75)] (2,0);
\draw[thick] (2,0) to [controls=+(30:0.25) and +(150:0.25)] (3,0);
\draw[thick] (5,0) to [controls=+(30:0.25) and +(150:0.25)] (6,0);
\foreach \Point in {(0,0), (0,3), (6,3), (6,0), (2,0), (3,0), (5,0)}{
    \node[scale=1] at \Point {\textbullet};
}
\draw node at (3,2) {$a_i$-tile};
\draw node at (-0.5,3.25) {\scalebox{0.8}{$(x,y)$}};
\draw node at (-0.5,-0.25) {\scalebox{0.8}{$(x,y-\log_2(\lambda))$}};

\draw[<->] (0.1,4.25) -- (5.9,4.25);
\draw node at (3,4.5) {$v(a_i)\cdot 2^y$};

\draw[<->] (-0.5,0) -- (-0.5,3);
\draw node at (-1.25,1.5) {\scalebox{0.8}{$\log_2(\lambda)$}};

\draw[<->] (0.1,-0.5) -- (1.9,-0.5);
\draw node at (1,-1) {\scalebox{0.8}{$\frac{1}{\lambda}v(w_1)\cdot 2^y$}};

\draw[<->] (2.1,-0.5) -- (2.9,-0.5);
\draw node at (2.5,-1) {\scalebox{0.8}{$\frac{1}{\lambda}v(w_2)\cdot 2^y$}};

\draw node at (4,-0.5) {$\dots$};
\draw node at (4,0) {$\dots$};
\draw node at (4,-1) {$\dots$};

\draw[<->] (5.1,-0.5) -- (5.9,-0.5);
\draw node at (5.5,-1) {\scalebox{0.8}{$\frac{1}{\lambda}v(w_k)\cdot 2^y$}};

\draw node at (8,2) {$\sigma(a_i)=w_1\dots w_k$};
\end{tikzpicture}
\end{bigcenter}
\caption{An $a_i$-tile for some letter $a_i\in \A$ with $\sigma(a_i)=w_1\dots w_k$.}
\label{figure.ab_tile}
\end{figure}

Let $\sigma$ be a primitive substitution with an expanding eigenvalue $\lambda>1$ and let $v\in(\R^{+})^{|\A|}$ be an associated eigenvector. For every letter $a\in \A$, define the \define{$a$-tile in position $(x,y)\in\R^2$}\label{def.sigma_tiles} as the square polygon with $|w|+3$ edges pictured in \ref{figure.ab_tile}, where $\sigma(a)=w_{1}\dots w_{k}$ (horizontal edges are curved to be more visible, but are in fact just straight lines).

\begin{remark}\label{remark.tile_rectangular}
The length of the top edge and the sum of lengths of bottom edges of this tile are the same. Since $M_\sigma\cdot v=\lambda\cdot v$, one has
$$\sum_{j=1}^{k}\frac{1}{\lambda} v(w_i)\cdot 2^y =\frac{2^y}{\lambda}\cdot\lambda\cdot v(a)=v(a)\cdot 2^y,$$
so that the bottom right vertex $(x+\frac{1}{\lambda}\left(v(w_1)+\dots+v(w_k)\right) 2^y,y-\log_2(\lambda))$ is indeed $(x+v(a)\cdot 2^y,y-\log_2(\lambda))$.
\end{remark}

\begin{figure}[!ht]
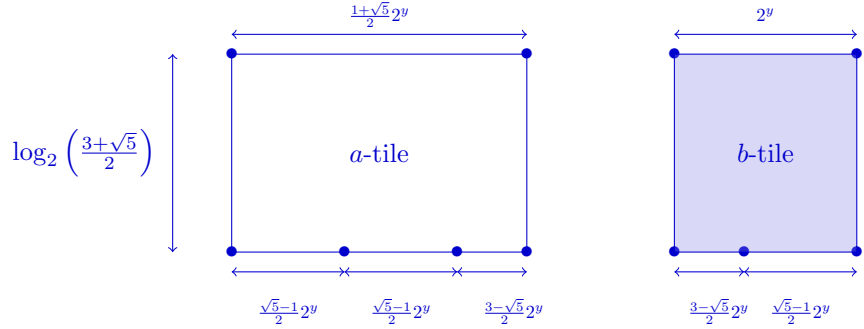

\begin{bigcenter}

\end{bigcenter}
\caption{The $\sigma$-tiles for $\sigma(a) = aab, \sigma(b) = ba$. For this substitution, $\lambda=\frac{3+\sqrt{5}}{2}$ and $v=\left( \frac{1+\sqrt{5}}{2},1\right)$.}
\label{figure.example_sigma_tiles}
\end{figure}

\begin{proposition}[\cite{ABM2019}, Proposition 5]\label{proposition.orbits_implies_tiling}
	If a substitution $\sigma$ has an expanding eigenvalue, then there exists a $\sigma$-tiling.
\end{proposition}

\begin{figure}[!ht]
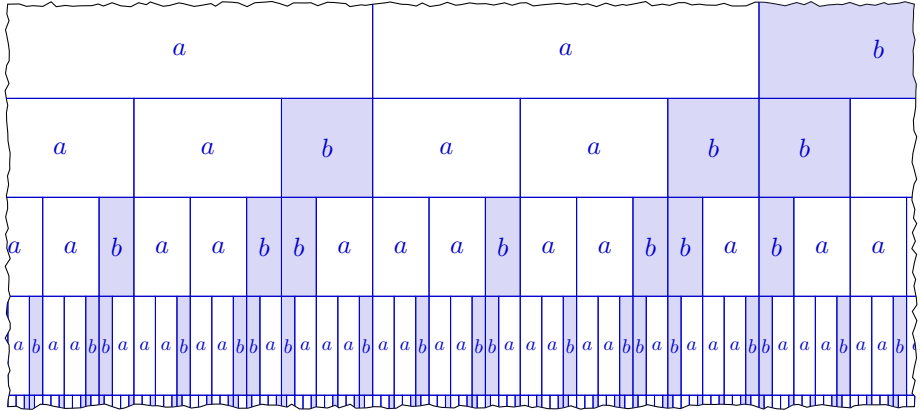

\begin{bigcenter}
%
\end{bigcenter}
\caption{A finite piece of a $\sigma$-tiling for $\sigma(a) = aab, \sigma(b) = ba$.}
\label{figure.example_sigma_tiling}
\end{figure}

If $\tau$ is a $\sigma$-tiling of $\R^2$, it can be seen as a mapping $\tau:\R^2\to\A$ with additional properties. Since $\A^{\R^2}$ is a compact space, from any sequence of $\sigma$-tilings one can extract a sequence that converges to some mapping $f:\R^2\to\A$.

\begin{proposition}\label{proposition.limit_sigma_tilings}
 If $(\tau_n)_{n\in\N}$ is a sequence of $\sigma$-tilings of $\R^2$ converging to some $f:\R^2\to\A$, then $\tau$ is a $\sigma$-tiling of $\R^2$.
\end{proposition}

\begin{proof}
Let $(\tau_n)_{n\in\N}$ be a sequence of $\sigma$-tilings of $\R^2$ converging to some $f:\R^2\to\A$. By definition, for every $k\in\N$, there exists a $N_k\in\N$ which satisfies that
$$\text{ for every }n\geq N_k,\; (\tau_n)_{|[-k;k]^2}=f_{|[-k;k]^2}.$$
Define a collection $\tau$ of $\sigma$-tiles as follows: an $a_i$-tile in position $(x,y)$ is in $\tau$ iff $\tau_{N_k}$ contains an $a_i$-tile in position $(x,y)$, where 
$$k:=\lceil\max\left(|x|,|y|,|y-\log_2(\lambda)|,|x+v_i\cdot2^y|\right)\rceil.$$
One can easily check that this collection $\tau$ is a tiling of $\R^2$, and that $\tau$ seen as a mapping $:\R^2\to\A$ coincide with $f$.
\end{proof}

%

\section{$\sigma$-tilings as a subshift on $BS(1,2)$}
\label{section.tilings_as_subshifts}

For all this section we fix $\A$ a finite alphabet with at least two letters $|\A|\geq2$ and $\sigma:\A\to\A^*$ a primitive substitution with an expanding eigenvalue $\lambda>1$. We denote $M_\sigma$ the matrix associated with $\sigma$ and $v\in(\R^{+})^{|\A|}$ an eigenvector associated with $\lambda$, so that $M_\sigma\cdot v=\lambda\cdot v$. In the sequel we will need a technical additional assumption, that we now describe. 

\medskip

First note that for every $k\in\N^*$, the substitution $\sigma^k$ has expanding eigenvalue $\lambda^k>1$ and the same eigenvector $v$. Thus choosing $k$ big enough and since $\lambda>1$, we can have $\log_2(\lambda^k)\geq3$. We then modify $\sigma$ to its $k$th power $\sigma^k$, and to shorten the notation we also denote it $\sigma$. Second, since $Cv$ also satisfy $M(Cv)=\lambda(Cv)$ for every $C>0$, we may assume without loss of generality that $v(a)\geq 3$ for every letter $a\in\A$. If this holds and that moreover $\lambda\geq 8$, we say that $\sigma$ and $v$ \define{have the unique size property}.

\medskip

We then introduce the notion of $\Phi(g)$-box for a group element $g\in BS(1,2)$, and list some useful properties on these boxes. If $g\in BS(1,2)$ and $g$ has normal form $g=t^ia^jt^{-k}$, the \define{$\Phi(g)$-box} is the rectangle $\Phi(g)+[0,2^{-i}[\times ]-1,0]$ in $\mathbb{R}^2$. A $\Phi$-box is simply a $\Phi(g)$-box for some group element $g\in BS(1,2)$.

\begin{remark}
    \label{remak.Phi_box_partition}
    The set of the $\Phi(g)$-boxes for all elements $g$ from a same sheet of $BS(1,n)$ is a partition of $\R^2$.
\end{remark}

    This remark is an important ingredient of the construction of Section~\ref{section.SFT_sigma_tilings}. Indeed, it is enough to know the content of all $\Phi$-boxes from the same sheet to know the content of the entire plane.

\begin{remark}\label{remark.content_box}
For a given $\sigma$-tiling, every $\Phi(g)$-box contains exactly one among
\begin{itemize}
 \item a vertical line of the $\sigma$-tiling;
 \item an horizontal line of the $\sigma$-tiling;
 \item a cross (intersection of a vertical and an horizontal lines) of the $\sigma$-tiling;
 \item a $T$ (a vertical line that start immediately below an horizontal line) of the $\sigma$-tiling;
 \item nothing.
\end{itemize}
\end{remark} 

\begin{proposition}\label{proposition.unique_size_property}
 If $\sigma$ and $v$ have the unique size property, then
 \begin{enumerate}
  \item a $\sigma$-tile vertically intersects $h$ or $h+1$ $\Phi$-boxes , where $h:=\lfloor\log_2(\lambda)\rfloor\geq 3$;
  \item for every $a\in\A$, the top of an $a$-tile intersects $w$ $\Phi$-boxes horizontally, where $3\leq v(a)\leq w< 2v(a)$;
  \item the number of $\Phi$-boxes horizontally intersected by  the bottom of a $\sigma$-tile uniquely determines 
  the number of $\Phi$-boxes vertically intersected.
 \end{enumerate}
\end{proposition}

\begin{proof}
 The first two items are direct consequences of the definitions of $\Phi$-boxes and $\sigma$-tiles. The third one is a little bit more subtile. Assume that an $a$-tile intersects $N$ $\Phi$-boxes on its bottom. Denote $\alpha:=\lceil\frac{N}{2^h}\rceil$ and $\beta:=\lceil\frac{N}{2^{h+1}}\rceil$ so that $\alpha$ and $\beta$ are the unique non-negative integers such that
 \begin{align*}
 \alpha &< \frac{N}{2^h} \leq \alpha\\
 \beta &< \frac{N}{2^{h+1}} \leq \beta.
 \end{align*}
Note that multiplying by $2$ the line immediately above leads to $\alpha=2\beta$. So either the $a$-tile vertically intersects $h$ $\Phi$-boxes, and then it horizontally intersects $\alpha$ $\Phi$-boxes on its top, or the $a$-tile vertically intersects $h+1$ $\Phi$-boxes, and then it horizontally intersects $\beta$ $\Phi$-boxes on its top. So the second item in Proposition~\ref{proposition.unique_size_property} thus implies that 
 \begin{align*}
 v(a)\leq \alpha< 2v(a)\\
 v(a)\leq  \beta< 2v(a).
 \end{align*}
But since $\alpha=2\beta$, one cannot have both at the same time and only one holds, which gives the number of $\Phi$-boxes vertically intersected by the $a$-tile.
\end{proof}


From now on we assume $\sigma$ and $v$ have the unique size property. Define the subshift $Y_\sigma\subset \A^{BS(1,2)}$ as
 \[Y_\sigma:=\left\{ y\in \A^{BS(1,2)}\mid \exists \tau\text{ a }\sigma\text{-tiling s.t. }y_g=\tau(\Phi(gat))\;\forall g\in BS(1,2)\right\}.\]
 
One may have expected to read $\tau(\Phi(g)) = \pi(x_g)$ instead of $\tau(\Phi(gat)) = \pi(x_g)$ in the definition of $Y_\sigma$. This is due to the definition of $\Phi$-box : if a $\Phi(g)$-box can see for instance the top left corner of a $\sigma$-tile (tiles $t_1$ or $t_2$), then $gat$ --i.e. the bottom right corner of the $\Phi(g)$-box-- is the only element that is certainly covered by this $\sigma$-tile in the $\sigma$-tiling.
 
\begin{proposition}\label{proposition.Y_sigma_subshift}
  The set $Y_\sigma$ is a subshift.
 \end{proposition}
 
 \begin{proof}
  
 We first prove that the set of configurations $Y_\sigma$ is $\mathfrak{S}$-invariant. Let $y\in Y_\sigma$, then there exists $\tau$ a $\sigma$-tiling of $\R^2$ such that $y_g=\tau(\Phi(gat))$ for every $g\in BS(1,2)$. Let $h\in BS(1,2)$ and denote $z:=\mathfrak{S}_h(y)$. Then  
 \begin{align*}
 z_g &= \left(\mathfrak{S}_h(y)\right)_g \\
  &= y_{h^{-1}\cdot g}  \\
  &= \tau\left( \Phi(h^{-1}\cdot gat)\right) \\
  &= \tau\left( (\alpha + 2^\beta x,\beta+y) \right) \text{ where }\Phi(h^{-1})=(\alpha,\beta)\text{ and }\Phi(gat)=(x,y)\\
  &= \left[ S_{ \Phi(h^{-1})}(\tau)\right](g),
\end{align*}
where $S$ is the $\R^2$-action given by 
$$S:
 \left(
 \begin{array}{lll}
  \R^2\times \R^2 & \to & \R^2\\
  \left((\alpha,\beta),(x,y)\right) & \mapsto & (\alpha + 2^\beta x,\beta+y)
 \end{array}
 \right).
 $$
 Since $\tau':=S_{ \Phi(h^{-1})}(\tau)$ is also a $\sigma$-tiling if $\tau$ is, we have that $z_g=\tau'(\Phi(gat))$ for every $g\in BS(1,2)$ and thus $z=\mathfrak{S}_h(y)$ is also in $Y_\sigma$.
 
 \medskip 
  
 We now check that $Y$ is closed. Take $(y_n)_{n\in\N}$ a sequence of configurations of $Y_\sigma$ that converges toward some $y\in\A^{BS(1,2)}$. Assume that for every $n\in\N$, $\tau_n$ is a $\sigma$-tiling such that $\tau_n(\Phi(g))=(y_n)_g$ for every $g\in BS(1,2)$. By compactness of $\A^{\R^2}$, there exists $\phi$ an extraction such that the subsequence $\left(\tau_{\phi(n)}\right)$ is converging toward some $\tau:\R^2\to\A$. Moreover Proposition~\ref{proposition.limit_sigma_tilings} implies that $\tau$ is a $\sigma$-tiling of $\R^2$. Define now the configuration $z\in Y$ by $\tau(\Phi(gat))=z_g$ for every $g\in BS(1,2)$. One can check by a standard compactness argument that $y=z$. Finally the set $Y_\sigma$ is $\mathfrak{S}$-invariant and closed, hence it is a subshift.
 \end{proof}    
 
 The subshift $Y_\sigma$ is for sure effectively closed, since one can enumerate, by increasing size, the complement of the language $\mathcal{L}(Y_\sigma)$. Note also that for some substitutions $\sigma$, the subshift $Y_\sigma$ is not an SFT. This is the case for Example~\ref{example.sigma_Y_not_SFT}.
 
 \begin{example}\label{example.sigma_Y_not_SFT}
 Let $\A=\{0,1\}$ and $\sigma$ the length 3 substitution given by $\sigma(0)=010$ and $\sigma(1)=101$. This substitution has two fixpoints $\omega=(01)^\Z$ and $\omega'=(10)^\Z$, and $\sigma$ is not recognizable: $\omega=\dots010101.0101010\dots$ can be parsed either as $\dots(010)(101).(010)(101)(0\dots$ or as $\dots0)(101)(01.0)(101)(010)\dots$ and thus has two different pre-images under $\sigma$. We use this property to show that $Y_\sigma$ is not an SFT. By contradiction assume that $Y_\sigma$ is an SFT, and $\F\subset \A^{R_{m,n}}$ is a finite set of forbidden patterns that defines it: $Y_\sigma=X_{\F}$. Then for $j> n$, the identity and the element $t^{j}at^{-j}$ cannot be both inside a same rectangle $g\cdot R_{m,n}$.
 
 \medskip
 
 We construct two configurations $x$ and $y$ in $Y_\sigma$, that arise from the two $\sigma$-tilings $\tau$ and $\tau'$ pictured on Figure~\ref{figure.proof_Y_sigma_non_SFT} (left): the $\sigma$-tilings $\tau$ and $\tau'$ coincide on the bottom half plane $R\times\R^{-}$, but the choice to de-substitute the fixpoint $\omega$ immediately above that is different, so that $\tau$ and $\tau'$ never match again on the top half plane $R\times\R^{+}$. These two configurations $x$ and $y$ thus coincide on every group element $h=t^{i}a^jt^{-k}$ with $i>k$ and only there, see Figure~\ref{figure.proof_Y_sigma_non_SFT} (middle).

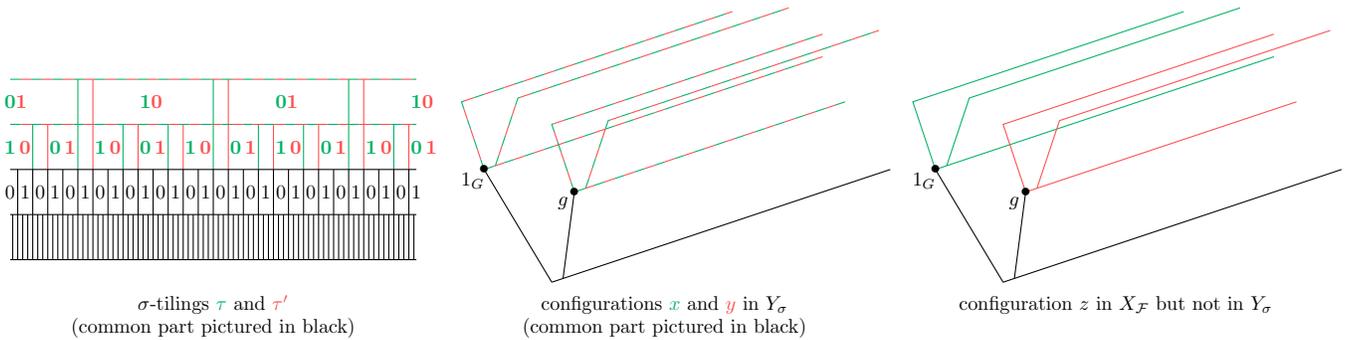
\begin{figure}[!ht]
\begin{bigcenter}
\begin{tikzpicture}[scale=0.6]

\begin{scope}[shift={(0,0)}]
 
\draw[vert,dash pattern= on 4pt off 4pt] (-4.5,5) -- (4.5,5);
\draw[rouge,dash pattern= on 4pt off 4pt,dash phase=4pt] (-4.5,5) -- (4.5,5);
\foreach \i in {-3,0,3}{
	\draw[vert] (\i,5) -- (\i,4);
	} 
\foreach \i in {-3,0,3}{
	\draw[rouge] (\i+1/3,5) -- (\i+1/3,4);
	} 	
\draw[vert,dash pattern= on 4pt off 4pt] (-4.5,4) -- (4.5,4);
\draw[rouge,dash pattern= on 4pt off 4pt,dash phase=4pt] (-4.5,4) -- (4.5,4);
\foreach \i in {-4,...,4}{
	\draw[vert] (\i,4) -- (\i,3);
	}	
\foreach \i in {-4,...,4}{
	\draw[rouge] (\i+1/3,4) -- (\i+1/3,3);
	}	
\draw (-4.5,3) -- (4.5,3);
\foreach \i in {-13,...,13}{
	\draw[] (\i/3,3) -- (\i/3,2);
	} 
\draw (-4.5,2) -- (4.5,2);
\foreach \i in {-40,...,40}{
	\draw[] (\i/9,2) -- (\i/9,1);
	} 
\draw (-4.5,1) -- (4.5,1);

\node[scale=0.75] at (-4.5,4.5) {\textcolor{vert}{\boldmath $0$}};
\node[scale=0.75] at (-1.5,4.5) {\textcolor{vert}{\boldmath $1$}};
\node[scale=0.75] at (1.5,4.5) {\textcolor{vert}{\boldmath $0$}};
\node[scale=0.75] at (4.5,4.5) {\textcolor{vert}{\boldmath $1$}};
\node[scale=0.75] at (-4.5,3.5) {\textcolor{vert}{\boldmath $1$}};
\node[scale=0.75] at (-3.5,3.5) {\textcolor{vert}{\boldmath $0$}};
\node[scale=0.75] at (-2.5,3.5) {\textcolor{vert}{\boldmath $1$}};
\node[scale=0.75] at (-1.5,3.5) {\textcolor{vert}{\boldmath $0$}};
\node[scale=0.75] at (-0.5,3.5) {\textcolor{vert}{\boldmath $1$}};
\node[scale=0.75] at (0.5,3.5) {\textcolor{vert}{\boldmath $0$}};
\node[scale=0.75] at (1.5,3.5) {\textcolor{vert}{\boldmath $1$}};
\node[scale=0.75] at (2.5,3.5) {\textcolor{vert}{\boldmath $0$}};
\node[scale=0.75] at (3.5,3.5) {\textcolor{vert}{\boldmath $1$}};
\node[scale=0.75] at (4.5,3.5) {\textcolor{vert}{\boldmath $0$}};
\node[scale=0.75] at (-4.5,2.5) {$0$};
\node[scale=0.75] at (-4.166,2.5) {$1$};
\node[scale=0.75] at (-3.833,2.5) {$0$};
\node[scale=0.75] at (-3.5,2.5) {$1$};
\node[scale=0.75] at (-3.166,2.5) {$0$};
\node[scale=0.75] at (-2.833,2.5) {$1$};
\node[scale=0.75] at (-2.5,2.5) {$0$};
\node[scale=0.75] at (-2.166,2.5) {$1$};
\node[scale=0.75] at (-1.833,2.5) {$0$};
\node[scale=0.75] at (-1.5,2.5) {$1$};
\node[scale=0.75] at (-1.166,2.5) {$0$};
\node[scale=0.75] at (-0.833,2.5) {$1$};
\node[scale=0.75] at (-0.5,2.5) {$0$};
\node[scale=0.75] at (-0.166,2.5) {$1$};
\node[scale=0.75] at (0.166,2.5) {$0$};
\node[scale=0.75] at (0.5,2.5) {$1$};
\node[scale=0.75] at (0.833,2.5) {$0$};
\node[scale=0.75] at (1.166,2.5) {$1$};
\node[scale=0.75] at (1.5,2.5) {$0$};
\node[scale=0.75] at (1.833,2.5) {$1$};
\node[scale=0.75] at (2.166,2.5) {$0$};
\node[scale=0.75] at (2.5,2.5) {$1$};
\node[scale=0.75] at (2.833,2.5) {$0$};
\node[scale=0.75] at (3.166,2.5) {$1$};
\node[scale=0.75] at (3.5,2.5) {$0$};
\node[scale=0.75] at (3.833,2.5) {$1$};
\node[scale=0.75] at (4.166,2.5) {$0$};
\node[scale=0.75] at (4.5,2.5) {$1$};

\begin{scope}[shift={(0.25,0)}]
\node[scale=0.75] at (-4.5,4.5) {\textcolor{rouge}{\boldmath $1$}};
\node[scale=0.75] at (-1.5,4.5) {\textcolor{rouge}{\boldmath $0$}};
\node[scale=0.75] at (1.5,4.5) {\textcolor{rouge}{\boldmath $1$}};
\node[scale=0.75] at (4.5,4.5) {\textcolor{rouge}{\boldmath $0$}};
\end{scope}
\begin{scope}[shift={(0.33,0)}]
\node[scale=0.75] at (-4.5,3.5) {\textcolor{rouge}{\boldmath $0$}};
\node[scale=0.75] at (-3.5,3.5) {\textcolor{rouge}{\boldmath $1$}};
\node[scale=0.75] at (-2.5,3.5) {\textcolor{rouge}{\boldmath $0$}};
\node[scale=0.75] at (-1.5,3.5) {\textcolor{rouge}{\boldmath $1$}};
\node[scale=0.75] at (-0.5,3.5) {\textcolor{rouge}{\boldmath $0$}};
\node[scale=0.75] at (0.5,3.5) {\textcolor{rouge}{\boldmath $1$}};
\node[scale=0.75] at (1.5,3.5) {\textcolor{rouge}{\boldmath $0$}};
\node[scale=0.75] at (2.5,3.5) {\textcolor{rouge}{\boldmath $1$}};
\node[scale=0.75] at (3.5,3.5) {\textcolor{rouge}{\boldmath $0$}};
\node[scale=0.75] at (4.5,3.5) {\textcolor{rouge}{\boldmath $1$}};
\end{scope}

\node[scale=0.75] at (0,0) {$\sigma$-tilings \textcolor{vert}{$\tau$} and \textcolor{rouge}{$\tau'$}};
\node[scale=0.75] at (0,-0.5) {(common part pictured in black)};
\end{scope}

\begin{scope}[shift={(10,0)}]

\draw[vert,dash pattern= on 4pt off 4pt] (-4,3) -- (3.5,5.5);
\draw[rouge,dash pattern= on 4pt off 4pt,dash phase=4pt] (-4,3) -- (3.5,5.5);

\draw[vert,dash pattern= on 4pt off 4pt] (-2,2.5) -- (4,4.5);
\draw[rouge,dash pattern= on 4pt off 4pt,dash phase=4pt] (-2,2.5) -- (4,4.5);

\draw[vert,dash pattern= on 4pt off 4pt] (-4,3) -- (-4.5,4.5) -- (1.5,6.5);
\draw[vert,dash pattern= on 4pt off 4pt] (-4+0.25,3+0.25/3) -- (-3.5+0.25,4.5+0.25/3) -- (-3.5+0.25+6,4.5+0.25/3+2);
\begin{scope}[shift={(2,-0.5)}]
\draw[vert,dash pattern= on 4pt off 4pt] (-4,3) -- (-4.5,4.5) -- (1.5,6.5);
\draw[vert,dash pattern= on 4pt off 4pt] (-4+0.25,3+0.25/3) -- (-3.5+0.25,4.5+0.25/3) -- (-3.5+0.25+6,4.5+0.25/3+2);
\end{scope}

\draw[rouge,dash pattern= on 4pt off 4pt,dash phase=4pt] (-4,3) -- (-4.5,4.5) -- (1.5,6.5);
\draw[rouge,dash pattern= on 4pt off 4pt,dash phase=4pt] (-4+0.25,3+0.25/3) -- (-3.5+0.25,4.5+0.25/3) -- (-3.5+0.25+6,4.5+0.25/3+2);
\begin{scope}[shift={(2,-0.5)}]
\draw[rouge,dash pattern= on 4pt off 4pt,dash phase=4pt] (-4,3) -- (-4.5,4.5) -- (1.5,6.5);
\draw[rouge,dash pattern= on 4pt off 4pt,dash phase=4pt] (-4+0.25,3+0.25/3) -- (-3.5+0.25,4.5+0.25/3) -- (-3.5+0.25+6,4.5+0.25/3+2);
\end{scope}

\draw[] (-2.5,0.5) -- (-4,3);
\draw[] (-2.5+0.25,0.5+0.25/3) -- (-2,2.5);
\draw[] (-2.5,0.5) -- (5,3);

\node[scale=0.75] at (0,0) {configurations \textcolor{vert}{$x$} and \textcolor{rouge}{$y$} in $Y_\sigma$};
\node[scale=0.75] at (0,-0.5) {(common part pictured in black)};

\node[scale=0.75] at (-4,3) {$\bullet$};
\node[scale=0.75] at (-4.25,2.75) {$1_G$};

\node[scale=0.75] at (-2,2.5) {$\bullet$};
\node[scale=0.75] at (-2.25,2.25) {$g$};
\end{scope}

\begin{scope}[shift={(20,0)}]

\draw[vert] (-4,3) -- (3.5,5.5);

\draw[rouge] (-2,2.5) -- (4,4.5);

\draw[vert] (-4,3) -- (-4.5,4.5) -- (1.5,6.5);
\draw[vert] (-4+0.25,3+0.25/3) -- (-3.5+0.25,4.5+0.25/3) -- (-3.5+0.25+6,4.5+0.25/3+2);

\begin{scope}[shift={(2,-0.5)}]
\draw[rouge] (-4,3) -- (-4.5,4.5) -- (1.5,6.5);
\draw[rouge] (-4+0.25,3+0.25/3) -- (-3.5+0.25,4.5+0.25/3) -- (-3.5+0.25+6,4.5+0.25/3+2);
\end{scope}

\draw[] (-2.5,0.5) -- (-4,3);
\draw[] (-2.5+0.25,0.5+0.25/3) -- (-2,2.5);
\draw[] (-2.5,0.5) -- (5,3);

\node[scale=0.75] at (0,0) {configuration $z$ in $X_\mathcal{F}$ but not in $Y_\sigma$}; 

\node[scale=0.75] at (-4,3) {$\bullet$};
\node[scale=0.75] at (-4.25,2.75) {$1_G$};

\node[scale=0.75] at (-2,2.5) {$\bullet$};
\node[scale=0.75] at (-2.25,2.25) {$g$};
\end{scope}

\end{tikzpicture}
\end{bigcenter}
\caption{The substitution given by $\sigma(0)=010$ and $\sigma(1)=101$ is not recognizable, hence the subshift $Y_\sigma$ is not SFT.}
\label{figure.proof_Y_sigma_non_SFT}
\end{figure}
 
 We now fix some $j> n$, define $g:=t^{j}at^{-j}$ and a new configuration $z\in\{0,1\}^{BS(1,2)}$ by

$$
\left\{\begin{array}{rl}
z_h & :=y_h \text{ if }h=g\cdot a^k \cdot t^{-\ell}\text{ for some }k\in\Z\text{ and }\ell\in\N;\\
z_h & :=x_h \text{ otherwise.}
\end{array}\right.
$$

 Then this configuration $z$ belongs to the subshift $X_\F$, since every pattern $z_{|h\cdot R_{m,n}}$ that appears in $z$ also appears in either $y$, if $h=g\cdot a^k \cdot t^{-\ell}$ for some $k\in\Z$ and $\ell\in\N$, or $x$ for every other group element $h$. Nevertheless, the configuration $z$ does not belong to $Y_\sigma$, since a $\sigma$-tiling that would correspond to $z$ should be at the same time $\tau$ and $\tau'$, which is impossible. Hence the subshift $Y_\sigma$ is not SFT.
 \end{example}
 
 In the rest of the article, we wonder what more can be said about $Y_\sigma$. In Section~\ref{section.SFT_sigma_tilings}, we define an SFT cover $X_\sigma$ for $Y_\sigma$. Section~\ref{section.X_sigma_SFT_cover_Y_sigma} is devoted to prove that the construction of $X_\sigma$ is correct. In Section~\ref{section.dynamical_properties_X_sigma} we then discuss dynamical properties of $Y_\sigma$.

\section{An SFT on $BS(1,2)$ to encode $\sigma$-tilings}
\label{section.SFT_sigma_tilings} 
 
We now describe an SFT $X_\sigma$ on $BS(1,2)$ that encodes $\sigma$-tilings, i.e. it is an SFT cover for the subshift $Y_\sigma$. The SFT is made so that a $\sigma$-tiling corresponds to a configuration in $X_\sigma$ and vice-versa (see Theorem~\ref{theorem.sigma_tilings_and_X_sigma}). We first define integers $h$ and $w_a$ for every letter $a\in\A$ by
\begin{itemize}
 \item $h:=\lfloor\log_2(\lambda)\rfloor-1$;
 \item for every $a\in\A$, $w_a:=\lfloor v(a)\rfloor-1$.
\end{itemize}
Since we have chosen $\lambda\geq8$, we have $h\geq2$. Similarly since the eigenvector $v$ was chosen so that $v(a)\geq3$ for every $a\in\A$, we get that $w_a\geq2$ for every letter $a\in\A$. 

\paragraph{Description of the alphabet}
 
Denote $M:=\max_{a\in\A} |\sigma(a)|$. We define a finite alphabet $\mathcal{R}$ that contains letters like $t_1\left(b,j,a,i\right)$, $m_3\left(a,i\right)$ or $b_2\left(a,i\right)$. More precisely, we consider the alphabet $\mathcal{R}=\mathcal{T}\cup\mathcal{M}\cup\mathcal{B}$, where $\mathcal{T}$ stands for \emph{top}, $\mathcal{B}$ for \emph{bottom} and $\mathcal{M}$ for \emph{middle}, where
\begin{align*}
\mathcal{T}&:=t_1\cup t_2 \cup t_3 \cup t_4\\
\mathcal{M}&:=m_1\cup m_2 \cup m_3 \cup m_4 \cup m_5\\
\mathcal{B}&:=b_1\cup b_2 \cup b_3
\end{align*}
with
\begin{align*}
\alpha:=&\left\{\alpha(b,j,a,i):a,b\in\A,i,j=1\dots M\right\}\\
\beta:=&\left\{\beta(a,i):a\in\A,i=1\dots M\right\}
\end{align*}
for $\alpha\in\{t_1,m_2\}$ and every $\beta\in\{t_2,t_3,t_4,m_1,m_3,m_4,m_5,b_1,b_2,b_3\}$. The letters themselves are the pentagonal Wang tiles described below.

\begin{bigcenter}
\begin{tikzpicture}[xscale=1,yscale=1]
 
\node[scale=1.25] at (-1.5,0.5) {$\mathcal{T}$};

\begin{scope}[shift={(1.5,0)}]
\begin{scope}[shift={(0,0)}]
\draw (0,0) rectangle (2,1);
\draw[color=bleu,pattern color=bleu,pattern=north west lines] (0,0) rectangle (0.5,0.5);
\draw[color=bleu,pattern color=bleu,pattern=north west lines] (1.5,0) rectangle (2,0.5);
\foreach \i in {(0,0),(1,0),(2,0),(0,1),(2,1)}{
	\node[scale=0.75] at \i {\textbullet};
	}
\node[scale=0.75] at (1,1.25) {$t_1(b,j,a,i)$};
\node[draw,rectangle, draw=bleu, fill=bleu!10, text width=2em,align=center, rounded corners=2, minimum height=1em,scale=0.4] at (0.25,0.25) {\textcolor{bleu}{$(b,j)$}};
\node[draw,rectangle, draw=bleu, fill=bleu!10, text width=2em,align=center, rounded corners=2, minimum height=1em,scale=0.4] at (1.75,0.25) {\textcolor{bleu}{$(a,i)$}};
\end{scope}


\begin{scope}[shift={(3,0)}]
\draw (0,0) rectangle (2,1);
\draw[color=bleu,pattern color=bleu,pattern=north west lines] (0.5,0) rectangle (2,0.5); 
\foreach \i in {(0,0),(1,0),(2,0),(0,1),(2,1)}{
	\node[scale=0.75] at \i {\textbullet};
	}
\node[scale=0.75] at (1,1.25) {$t_2(a,i)$};
\node[draw,rectangle, draw=bleu, fill=bleu!10, text width=2em,align=center, rounded corners=2, minimum height=1em,scale=0.6] at (1.25,0.25) {\textcolor{bleu}{$(a,i)$}};
\end{scope}


\begin{scope}[shift={(6,0)}]
\draw (0,0) rectangle (2,1);
\draw[color=bleu,pattern color=bleu,pattern=north west lines] (0,0) rectangle (2,0.5); 
\foreach \i in {(0,0),(1,0),(2,0),(0,1),(2,1)}{
	\node[scale=0.75] at \i {\textbullet};
	}
\node[scale=0.75] at (1,1.25) {$t_3(a,i)$};
\node[draw,rectangle, draw=bleu, fill=bleu!10, text width=2em,align=center, rounded corners=2, minimum height=1em,scale=0.6] at (1,0.25) {\textcolor{bleu}{$(a,i)$}};
\end{scope}


\begin{scope}[shift={(9,0)}]
\draw (0,0) rectangle (2,1);

\draw[color=bleu,pattern color=bleu,pattern=north west lines] (0,0) rectangle (1.5,0.5); 
\foreach \i in {(0,0),(1,0),(2,0),(0,1),(2,1)}{
	\node[scale=0.75] at \i {\textbullet};
	}
\node[scale=0.75] at (1,1.25) {$t_4(a,i)$};
\node[draw,rectangle, draw=bleu, fill=bleu!10, text width=2em,align=center, rounded corners=2, minimum height=1em,scale=0.6] at (0.75,0.25) {\textcolor{bleu}{$(a,i)$}};
\end{scope}

%
\end{scope}

\node[scale=1.25] at (-1.5,-1.5) {$\mathcal{M}$};

\begin{scope}[shift={(0,-2)}]
\draw (0,0) rectangle (2,1);
\draw[color=bleu,pattern color=bleu,pattern=north west lines] (0.5,0) -- (2,0) -- (2,1) -- (1,1) -- cycle ;
\foreach \i in {(0,0),(1,0),(2,0),(0,1),(2,1)}{
	\node[scale=0.75] at \i {\textbullet};
	}
\node[scale=0.75] at (1,1.25) {$m_1(a,i)$};
\node[draw,rectangle, draw=bleu, fill=bleu!10, text width=2em,align=center, rounded corners=2, minimum height=1em,scale=0.6] at (1.5,0.5) {\textcolor{bleu}{$(a,i)$}};
\end{scope}


\begin{scope}[shift={(3,-2)}]
\draw (0,0) rectangle (2,1);
\draw[color=bleu,pattern color=bleu,pattern=north west lines] (1.5,0) -- (2,0) -- (2,1) -- (1,1) -- cycle ; 
\draw[color=bleu,pattern color=bleu,pattern=north west lines] (0,0) -- (0.5,0) -- (0,0.5) -- cycle ; 
\foreach \i in {(0,0),(1,0),(2,0),(0,1),(2,1)}{
	\node[scale=0.75] at \i {\textbullet};
	}
\node[scale=0.75] at (1,1.25) {$m_2(b,j,a,i)$};
\node[draw,rectangle, draw=bleu, fill=bleu!10, text width=2em,align=center, rounded corners=2, minimum height=1em,scale=0.6] at (1.6,0.6) {\textcolor{bleu}{$(a,i)$}};
\node[draw,rectangle, draw=bleu, fill=bleu!10, text width=2em,align=center, rounded corners=2, minimum height=1em,scale=0.4] at (0.25,0.25) {\textcolor{bleu}{$(b,j)$}};
\end{scope}


\begin{scope}[shift={(6,-2)}]
\draw (0,0) rectangle (2,1);
\draw[color=bleu,pattern color=bleu,pattern=north west lines] (0,0) rectangle (2,1); 
\foreach \i in {(0,0),(1,0),(2,0),(0,1),(2,1)}{
	\node[scale=0.75] at \i {\textbullet};
	}
\node[scale=0.75] at (1,1.25) {$m_3(a,i)$};
\node[draw,rectangle, draw=bleu, fill=bleu!10, text width=2em,align=center, rounded corners=2, minimum height=1em,scale=0.6] at (1,0.5) {\textcolor{bleu}{$(a,i)$}};
\end{scope}



\begin{scope}[shift={(9,-2)}]
\draw (0,0) rectangle (2,1);
\draw[color=bleu,pattern color=bleu,pattern=north west lines] (0,0) -- (1.5,0) -- (1,1) -- (0,1) -- cycle ;   
\foreach \i in {(0,0),(1,0),(2,0),(0,1),(2,1)}{
	\node[scale=0.75] at \i {\textbullet};
	}
\node[scale=0.75] at (1,1.25) {$m_4(a,i)$};
\node[draw,rectangle, draw=bleu, fill=bleu!10, text width=2em,align=center, rounded corners=2, minimum height=1em,scale=0.6] at (0.5,0.5) {\textcolor{bleu}{$(a,i)$}};
\end{scope}



\begin{scope}[shift={(12,-2)}]
\draw (0,0) rectangle (2,1);
\draw[color=bleu,pattern color=bleu,pattern=north west lines] (0,0) -- (2,0) -- (2,0.5) -- (1,1) -- (0,1) -- cycle ;
\foreach \i in {(0,0),(1,0),(2,0),(0,1),(2,1)}{
	\node[scale=0.75] at \i {\textbullet};
	}
\node[scale=0.75] at (1,1.25) {$m_5(a,i)$};
\node[draw,rectangle, draw=bleu, fill=bleu!10, text width=2em,align=center, rounded corners=2, minimum height=1em,scale=0.6] at (1,0.5) {\textcolor{bleu}{$(a,i)$}};
\end{scope}

\node[scale=1.25] at (-1.5,-3.5) {$\mathcal{B}$};

\begin{scope}[shift={(3,-4)}]
\draw (0,0) rectangle (2,1);
\draw[color=bleu,pattern color=bleu,pattern=north west lines] (1,0.5) rectangle (2,1); 
\foreach \i in {(0,0),(1,0),(2,0),(0,1),(2,1)}{
	\node[scale=0.75] at \i {\textbullet};
	}
\node[scale=0.75] at (1,1.25) {$b_1(a,i)$};
\node[draw,rectangle, draw=bleu, fill=bleu!10, text width=2em,align=center, rounded corners=2, minimum height=1em,scale=0.6] at (1.5,0.75) {\textcolor{bleu}{$(a,i)$}};
\end{scope}


\begin{scope}[shift={(6,-4)}]
\draw (0,0) rectangle (2,1);
\draw[color=bleu,pattern color=bleu,pattern=north west lines] (0,0.5) rectangle (2,1); 
\foreach \i in {(0,0),(1,0),(2,0),(0,1),(2,1)}{
	\node[scale=0.75] at \i {\textbullet};
	}
\node[scale=0.75] at (1,1.25) {$b_2(a,i)$};
\node[draw,rectangle, draw=bleu, fill=bleu!10, text width=2em,align=center, rounded corners=2, minimum height=1em,scale=0.6] at (1,0.75) {\textcolor{bleu}{$(a,i)$}};
\end{scope}


\begin{scope}[shift={(9,-4)}]
\draw (0,0) rectangle (2,1);
\draw[color=bleu,pattern color=bleu,pattern=north west lines] (0,0.5) rectangle (1,1); 
\foreach \i in {(0,0),(1,0),(2,0),(0,1),(2,1)}{
	\node[scale=0.75] at \i {\textbullet};
	}
\node[scale=0.75] at (1,1.25) {$b_3(a,i)$};
\node[draw,rectangle, draw=bleu, fill=bleu!10, text width=2em,align=center, rounded corners=2, minimum height=1em,scale=0.6] at (0.5,0.75) {\textcolor{bleu}{$(a,i)$}};
\end{scope}

\end{tikzpicture}

\end{bigcenter}

In the case of a letter $t_1\left(b,j,a,i\right)$ or $m_2\left(b,j,a,i\right)$, we adopt the convention that $(a,i)$ is associated with the rightmost shaded area, and $(b,j)$ with the leftmost one.

\paragraph{Intuition on the alphabet}

Informally, every letter in $r\in\mathcal{R}$ gives some local  information about the $\sigma$-tiling encoded. Assume $x\in X_\sigma$ and $g\in BS(1,2)$ with normal form $g=t^ia^jt^{-k}$, so that $\Phi(g)=(j\cdot 2^{-i},k-i)$. Then the letter $r=x_g$ describes what the $\Phi(g)$-box
geometrically captures from the $\sigma$-tiling (which $\sigma$-tile is seen, and where it is with horizontal error $2^{k-i}$) and also tells in which position the $\sigma$-tile in from the other $\sigma$-tile immediately above it. 

 \begin{figure}[!ht]
\begin{bigcenter}
 \begin{tikzpicture}[scale=0.7]

\draw[color=vert,pattern color=vert,pattern=north west lines] (1,0) rectangle (2,1);
\draw[color=vert,very thick] (1.3,-0.5) -- (1.3,0.8) -- (2.5,0.8);
\draw (0,0) rectangle (2,1);
\node[scale=0.75] at (0.75,2.25) {$x_g=t_1(.,.,a,.)$};
\node[scale=0.75] at (0.65,1.75) {$\lettertone$};
\node[scale=0.75,rotate=90] at (1,1.325) {$\Leftrightarrow$};
\node[scale=0.5] at (0,1) {\textbullet};
\node[anchor=east,scale=0.75] at (0,1) {$\Phi(g)$};
\node[draw,rectangle, draw=vert, fill=vert!10,align=center, rounded corners=2, minimum height=1em,scale=0.75] at (2.25,0) {\textcolor{vert}{\boldmath $a$-tile}};

\begin{scope}[shift={(5,0)}]
\draw[color=vert,pattern color=vert,pattern=north west lines] (0,0) rectangle (1,1);
\draw[color=vert,very thick] (0.3,-0.5) -- (0.3,0.8) -- (2.5,0.8);
\draw (0,0) rectangle (2,1);
\node[scale=0.75] at (0.75,2.25) {$x_g=t_2(a,.)$};
\node[scale=0.75] at (0.65,1.75) {$\letterttwo$};
\node[scale=0.75,rotate=90] at (1,1.325) {$\Leftrightarrow$};
\node[scale=0.5] at (0,1) {\textbullet};
\node[anchor=east,scale=0.75] at (0,1) {$\Phi(g)$};
\node[draw,rectangle, draw=vert, fill=vert!10,align=center, rounded corners=2, minimum height=1em,scale=0.75] at (1.25,0) {\textcolor{vert}{\boldmath $a$-tile}};
\end{scope}

\begin{scope}[shift={(10,0)}]
\draw[color=vert,pattern color=vert,pattern=north west lines] (0,0) rectangle (2,1);
\draw[color=vert,very thick] (-0.5,0.3) -- (2.5,0.3);
\draw (0,0) rectangle (2,1);
\node[scale=0.75] at (0.75,2.25) {$x_g=t_3(a,.)$};
\node[scale=0.75] at (0.65,1.75) {$\lettertthree$};
\node[scale=0.75,rotate=90] at (1,1.325) {$\Leftrightarrow$};
\node[scale=0.5] at (0,1) {\textbullet};
\node[anchor=east,scale=0.75] at (0,1) {$\Phi(g)$};
\node[draw,rectangle, draw=vert, fill=vert!10,align=center, rounded corners=2, minimum height=1em,scale=0.75] at (1,-0.15) {\textcolor{vert}{\boldmath $a$-tile}};
\end{scope}

\begin{scope}[shift={(15,0)}]
\draw[color=vert,pattern color=vert,pattern=north west lines] (0,0) rectangle (2,1);
\draw[color=vert,very thick] (-0.5,0.3) -- (2.5,0.3);
\draw (0,0) rectangle (2,1);
\node[scale=0.75] at (0.75,2.25) {$x_g=t_4(a,.)$};
\node[scale=0.75] at (0.65,1.75) {$\lettertfour$};
\node[scale=0.75,rotate=90] at (1,1.325) {$\Leftrightarrow$};
\node[scale=0.5] at (0,1) {\textbullet};
\node[anchor=east,scale=0.75] at (0,1) {$\Phi(g)$};
\node[draw,rectangle, draw=vert, fill=vert!10,align=center, rounded corners=2, minimum height=1em,scale=0.75] at (1,-0.15) {\textcolor{vert}{\boldmath $a$-tile}};
\end{scope}


\begin{scope}[shift={(3,-4)}]
\draw[color=vert,pattern color=vert,pattern=north west lines] (0,0) rectangle (2,1);
\draw[color=vert,very thick] (0.7,-0.5) -- (0.7,1.5);
\draw (0,0) rectangle (2,1);
\node[scale=0.75] at (0.75,2.25) {$x_g=b_1(a,.)$};
\node[scale=0.75] at (0.65,1.75) {$\letterbone$};
\node[scale=0.75,rotate=90] at (1,1.325) {$\Leftrightarrow$};
\node[scale=0.5] at (0,1) {\textbullet};
\node[anchor=east,scale=0.75] at (0,1) {$\Phi(g)$};
\node[draw,rectangle, draw=vert, fill=vert!10,align=center, rounded corners=2, minimum height=1em,scale=0.75] at (1.5,0.5) {\textcolor{vert}{\boldmath $a$-tile}};
\end{scope}

\begin{scope}[shift={(8,-4)}]
\draw (0,0) rectangle (2,1);
\node[scale=0.75] at (0.75,2.25) {$x_g=b_2(a,.)$};
\node[scale=0.75] at (0.65,1.75) {$\letterbtwo$};
\node[scale=0.75,rotate=90] at (1,1.325) {$\Leftrightarrow$};
\node[scale=0.5] at (0,1) {\textbullet};
\node[anchor=east,scale=0.75] at (0,1) {$\Phi(g)$};
\node[draw,rectangle, draw=vert, fill=vert!10,align=center, rounded corners=2, minimum height=1em,scale=0.75] at (1,0.5) {\textcolor{vert}{\boldmath $a$-tile}};
\end{scope}

\begin{scope}[shift={(13,-4)}]
\draw (0,0) rectangle (2,1);
\node[scale=0.75] at (0.75,2.25) {$x_g=b_3(a,.)$};
\node[scale=0.75] at (0.65,1.75) {$\letterbthree$};
\node[scale=0.75,rotate=90] at (1,1.325) {$\Leftrightarrow$};
\node[scale=0.5] at (0,1) {\textbullet};
\node[anchor=east,scale=0.75] at (0,1) {$\Phi(g)$};
\node[draw,rectangle, draw=vert, fill=vert!10,align=center, rounded corners=2, minimum height=1em,scale=0.75] at (1,0.5) {\textcolor{vert}{\boldmath $a$-tile}};
\end{scope}


\begin{scope}[shift={(-3,-8)}]
\draw[color=vert,pattern color=vert,pattern=north west lines] (0,0) rectangle (1,1);
\draw[color=vert,very thick] (0.7,0) -- (0.7,1);
\draw (0,0) rectangle (2,1);
\draw[color=vert,very thick] (0.7,-0.5) -- (0.7,1.5);
\draw (0,0) rectangle (2,1);
\node[scale=0.75] at (0.75,2.25) {$x_g=m_1(a,.)$};
\node[scale=0.75] at (0.65,1.75) {$\lettermone$};
\node[scale=0.75,rotate=90] at (1,1.325) {$\Leftrightarrow$};
\node[scale=0.5] at (0,1) {\textbullet};
\node[anchor=east,scale=0.75] at (0,1) {$\Phi(g)$};
\node[draw,rectangle, draw=vert, fill=vert!10,align=center, rounded corners=2, minimum height=1em,scale=0.75] at (1.5,0.5) {\textcolor{vert}{\boldmath $a$-tile}};
\end{scope}

\begin{scope}[shift={(2,-8)}]
\draw[color=vert,pattern color=vert,pattern=north west lines] (1,0) rectangle (2,1);
\draw[color=vert,very thick] (1.2,-0.5) -- (1.2,1.5);
\draw (0,0) rectangle (2,1);
\node[scale=0.75] at (0.75,2.25) {$x_g=m_2(a,.)$};
\node[scale=0.75] at (0.65,1.75) {$\lettermtwo$};
\node[scale=0.75,rotate=90] at (1,1.325) {$\Leftrightarrow$};
\node[scale=0.5] at (0,1) {\textbullet};
\node[anchor=east,scale=0.75] at (0,1) {$\Phi(g)$};
\node[draw,rectangle, draw=vert, fill=vert!10,align=center, rounded corners=2, minimum height=1em,scale=0.75] at (2,0.5) {\textcolor{vert}{\boldmath $a$-tile}};
\end{scope}

\begin{scope}[shift={(7,-8)}]
\draw (0,0) rectangle (2,1);
\node[scale=0.75] at (0.75,2.25) {$x_g=m_3(a,.)$};
\node[scale=0.75] at (0.65,1.75) {$\lettermthree$};
\node[scale=0.75,rotate=90] at (1,1.325) {$\Leftrightarrow$};
\node[scale=0.5] at (0,1) {\textbullet};
\node[anchor=east,scale=0.75] at (0,1) {$\Phi(g)$};
\node[draw,rectangle, draw=vert, fill=vert!10,align=center, rounded corners=2, minimum height=1em,scale=0.75] at (1,0.5) {\textcolor{vert}{\boldmath $a$-tile}};
\end{scope}

\begin{scope}[shift={(12,-8)}]
\draw (0,0) rectangle (2,1);
\node[scale=0.75] at (0.75,2.25) {$x_g=m_4(a,.)$};
\node[scale=0.75] at (0.65,1.75) {$\lettermfour$};
\node[scale=0.75,rotate=90] at (1,1.325) {$\Leftrightarrow$};
\node[scale=0.5] at (0,1) {\textbullet};
\node[anchor=east,scale=0.75] at (0,1) {$\Phi(g)$};
\node[draw,rectangle, draw=vert, fill=vert!10,align=center, rounded corners=2, minimum height=1em,scale=0.75] at (1,0.5) {\textcolor{vert}{\boldmath $a$-tile}};
\end{scope}

\begin{scope}[shift={(17,-8)}]
\draw (0,0) rectangle (2,1);
\node[scale=0.75] at (0.75,2.25) {$x_g=m_5(a,.)$};
\node[scale=0.75] at (0.65,1.75) {$\lettermfive$};
\node[scale=0.75,rotate=90] at (1,1.325) {$\Leftrightarrow$};
\node[scale=0.5] at (0,1) {\textbullet};
\node[anchor=east,scale=0.75] at (0,1) {$\Phi(g)$};
\node[draw,rectangle, draw=vert, fill=vert!10,align=center, rounded corners=2, minimum height=1em,scale=0.75] at (1,0.5) {\textcolor{vert}{\boldmath $a$-tile}};
\end{scope}

\end{tikzpicture}
\end{bigcenter}
\caption{Signification of the letter $x_g$ that appears in position $g\in BS(1,2)$ inside a configuration $x\in X_\sigma$. The black rectangle with top left corner labeled by $\Phi(g)$ is the $\Phi(g)$-box.}
\label{figure.alphabet_tiling_intuition}
 \end{figure}
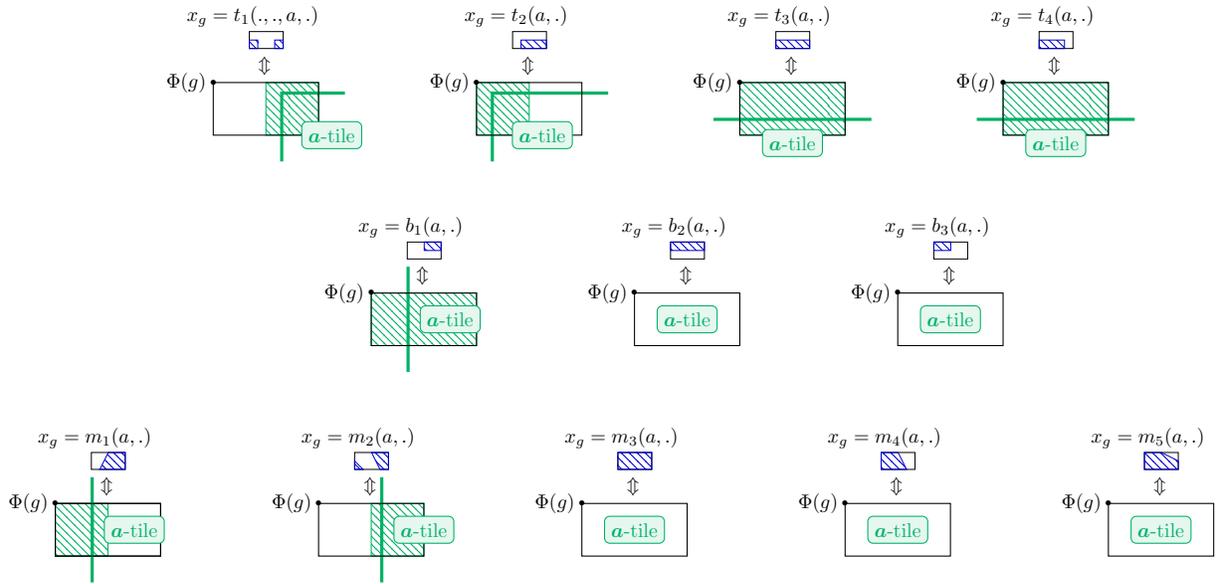

The letter $r$ may code the presence of a vertical line (letters $m_1$, $m_2$ and $b_1$), of an horizontal line (letters $t_3$ and $t_4$), of a corner (letters $t_1$ and $t_2$) or the absence of borders of $\sigma$-tiles (letters $m_3$, $m_4$, $m_5$, $b_2$ and $b_3$).  Figure~\ref{figure.alphabet_tiling_intuition} should be interpreted this way: if a letter $r\in\mathcal{R}$ codes the presence of an element (corner, vertical or horizontal line), it must be located somewhere inside the green shaded area. For instance, the fact that $x_g\in t_1$ is equivalent to the presence of a top left corner somewhere in $\Phi(g)+[2^{k-i};2^{k-i+1}[\times]-1;0]$. Similarly, the fact that $x_g\in m_2$ is equivalent to the presence of a vertical line $\{j\cdot2^{i}+x\}\times[k-i,k-i+1[$ with $x\in [2^{k-i};2^{k-i+1}[$. It follows that going down -- following generator $t$ -- inside a configuration $x\in X_\sigma$ gives a better approximation of where $\sigma$-tiles are.

\paragraph{Description of the local rules}

Denote $W:=\max_{a\in\A} w_a+1$. The set $\mathcal{A}_\sigma$ of allowed patterns defining $X_\sigma\subset \mathcal{R}^{BS(1,2)}$ that we describe below is actually a subset of $\mathcal{R}^{R_{3M,3(h+1)}}$. For the purpose of clarity, we do not present directly the set $\mathcal{A}_\sigma$, but instead list some local rules that describe which patterns are in $\mathcal{A}_\sigma$
\begin{enumerate}
 \item \label{rule.matching_rules} colors should match on edges along generators $a,a^{-1},t$ and $t ^{-1}$ -- this condition is nearest neighbor;
 \begin{bigcenter}
 \begin{tikzpicture}[scale=1]

\begin{scope}[shift={(0,0)}]
\draw (0,0) rectangle (2,1);
\draw[color=bleu,pattern color=bleu,pattern=north west lines] (0,0) -- (2,0) -- (2,0.5) -- (1,1) -- (0,1) -- cycle ;
\node[scale=0.75] at (0,1) {\textbullet};
\node[scale=0.75] at (0,1.25) {$g$};
\node[draw,rectangle, draw=bleu, fill=bleu!10, text width=2em,align=center, rounded corners=2, minimum height=1em,scale=0.6] at (1,0.5) {\textcolor{bleu}{$(b,j)$}};
\end{scope}

\begin{scope}[shift={(2,0)}]
\draw (0,0) rectangle (2,1);
\draw[color=bleu,pattern color=bleu,pattern=north west lines] (1.5,0) -- (2,0) -- (2,1) -- (1,1) -- cycle ; 
\draw[color=bleu,pattern color=bleu,pattern=north west lines] (0,0) -- (0.5,0) -- (0,0.5) -- cycle ;
\node[scale=0.75] at (0,1) {\textbullet};
\node[scale=0.75] at (0,1.25) {$g\cdot a$};
\node[draw,rectangle, draw=bleu, fill=bleu!10, text width=2em,align=center, rounded corners=2, minimum height=1em,scale=0.6] at (1.6,0.6) {\textcolor{bleu}{$(a,i)$}};
\node[draw,rectangle, draw=bleu, fill=bleu!10, text width=2em,align=center, rounded corners=2, minimum height=1em,scale=0.4] at (0.25,0.25) {\textcolor{bleu}{$(b,j)$}};
\end{scope}

\node[scale=1] at (2,-1) {$
\left\{\begin{array}{cl}
x_{g}&=m_5(b,j)\\
x_{g\cdot a}&=m_2(b,j,a,i)
\end{array}\right.$ is allowed};

\begin{scope}[shift={(-3,0)}]
\begin{scope}[shift={(10,0.5)}]
\draw (0,0) rectangle (2,1);
\draw[color=bleu,pattern color=bleu,pattern=north west lines] (0,0) rectangle (0.5,0.5);
\draw[color=bleu,pattern color=bleu,pattern=north west lines] (1.5,0) rectangle (2,0.5);
\node[scale=0.75] at (0,1) {\textbullet};
\node[scale=0.75] at (0,1.25) {$g$};
\node[draw,rectangle, draw=bleu, fill=bleu!10, text width=2em,align=center, rounded corners=2, minimum height=1em,scale=0.4] at (0.25,0.25) {\textcolor{bleu}{$(b,j)$}};
\node[draw,rectangle, draw=bleu, fill=bleu!10, text width=2em,align=center, rounded corners=2, minimum height=1em,scale=0.4] at (1.75,0.25) {\textcolor{bleu}{$(a,i)$}};
\end{scope}

\begin{scope}[shift={(10,0)},scale=0.5]
\draw (0,0) rectangle (2,1);
\draw[color=bleu,pattern color=bleu,pattern=north west lines] (0,0) -- (1.5,0) -- (1,1) -- (0,1) -- cycle ;
\node[scale=0.75] at (0,1) {\textbullet};
\node[scale=0.75] at (-0.75,1.25) {$g\cdot t$};
\node[draw,rectangle, draw=bleu, fill=bleu!10, text width=2em,align=center, rounded corners=2, minimum height=1em,scale=0.3] at (0.5,0.5) {\textcolor{bleu}{$(a,i)$}};
\end{scope}

\node[scale=1] at (12,-1) {$
\left\{\begin{array}{cl}
x_{g}&=t_1(b,j,a,i)\\
x_{g\cdot t}&=m_4(a,i)
\end{array}\right.$ is allowed};
\end{scope}

\end{tikzpicture}
\end{bigcenter}
 \item \label{rule.type_TMB} the type of letter of the first coordinate ($\mathcal{T}$, $\mathcal{M}$ or $\mathcal{B}$) is uniform on the set $\left\{ g\cdot a^k: k\in\Z\right\}$ for every $g\in BS(1,2)$ -- this condition is nearest neighbor;
 \begin{bigcenter}
\begin{tikzpicture}[scale=1]

\begin{scope}[shift={(0,0)}]
\draw (0,0) rectangle (2,1);
\draw[color=bleu,pattern color=bleu,pattern=north west lines] (0,0) rectangle (0.5,0.5);
\draw[color=bleu,pattern color=bleu,pattern=north west lines] (1.5,0) rectangle (2,0.5);
\node[scale=0.75] at (0,1) {\textbullet};
\node[scale=0.75] at (0,1.25) {$g$};
\node[draw,rectangle, draw=bleu, fill=bleu!10, text width=2em,align=center, rounded corners=2, minimum height=1em,scale=0.6] at (1.6,0.6) {\textcolor{bleu}{$(a,i)$}};
\node[draw,rectangle, draw=bleu, fill=bleu!10, text width=2em,align=center, rounded corners=2, minimum height=1em,scale=0.4] at (0.25,0.25) {\textcolor{bleu}{$(b,j)$}};
\end{scope}

\begin{scope}[shift={(2,0)}]
\draw (0,0) rectangle (2,1);
\draw[color=bleu,pattern color=bleu,pattern=north west lines] (1.5,0) -- (2,0) -- (2,1) -- (1,1) -- cycle ; 
\draw[color=bleu,pattern color=bleu,pattern=north west lines] (0,0) -- (0.5,0) -- (0,0.5) -- cycle ; 
\node[scale=0.75] at (0,1) {\textbullet};
\node[scale=0.75] at (0,1.25) {$g\cdot a$};
\node[draw,rectangle, draw=bleu, fill=bleu!10, text width=2em,align=center, rounded corners=2, minimum height=1em,scale=0.6] at (1.6,0.6) {\textcolor{bleu}{$(c,k)$}};
\node[draw,rectangle, draw=bleu, fill=bleu!10, text width=2em,align=center, rounded corners=2, minimum height=1em,scale=0.4] at (0.25,0.25) {\textcolor{bleu}{$(a,i)$}};
\end{scope}

\node[scale=1] at (2,-1) {$
\left\{\begin{array}{cl}
x_{g}&=t_1(b,j,a,i)\\
x_{g\cdot a}&=m_2(a,i,c,k)
\end{array}\right.$ is forbidden (even if colors match)};

\end{tikzpicture}
\end{bigcenter} 
  
 \item \label{rule.dimensions} the dimensions of the connected shaded zones are bounded in coherence with the substitution $\sigma$. More precisely: there are two possible heights for shaded zones: $h$ and $h-1$, independently from the letter $a\in\A$ coded. The top width of a shaded zone does depend on the letter $a\in\A$ coded, and can only take values from $(w_a-1)$ to $(2w_a-1)$. 
 
 \item \label{rule.synchro_two_sheets} the letters in $g$, $gtat^{-1}$ and $gta^{-1}t^{-1}$ are synchronized:
 \begin{enumerate}
  \item $x_{g}=t_1(b,j,a,i) \Rightarrow \left\{\begin{array}{l} x_{g\cdot tat^{-1}}=t_2(a,i)\\ x_{g\cdot ta^{-1}t^{-1}}=t_4(b,j)\end{array}\right.$
  \item $x_{g}=b_1(a,i) \Rightarrow x_{g\cdot bab^{-1}}=b_1(a,i)\text{ or } x_{g\cdot ta^{-1}t^{-1}}=b_1(a,i)$
  \item $x_{g}=b_3(a,i) \Rightarrow x_{g\cdot bab^{-1}}=b_3(a,i)\text{ or } x_{g\cdot ba^{-1}t^{-1}}=b_3(a,i)$
  \item $x_{g}=m_1(a,i) \Leftrightarrow  x_{g\cdot ta^{-1}t^{-1}}=m_2(b,j,a,i)$
 \end{enumerate}
    \begin{bigcenter}
     \begin{tikzpicture}[scale=1]

\draw (0,0) -- (2,0) -- (1.5,1) -- (-0.5,1) -- cycle;
\draw (2,0) -- (4,0) -- (3.5,1) -- (1.5,1) -- cycle;
\draw (1,0) -- (3,0) -- (3,1.25) -- (1,1.25) -- cycle;
\draw[color=vert,pattern color=vert,pattern=north west lines] (0,0) -- (1.5,0) -- (1.25,0.5) -- (-0.25,0.5) -- cycle ;
\draw[color=vert,pattern color=vert,pattern=north west lines] (2.5,0) -- (4,0) -- (3.75,0.5) -- (2.25,0.5) -- cycle ;
\draw[color=bleu,pattern color=bleu,pattern=north west lines] (1,0) rectangle (1.5,0.625);
\draw[color=bleu,pattern color=bleu,pattern=north west lines] (2.5,0) rectangle (3,0.625);
\foreach \i in {(-0.5,1),(1.5,1),(1,1.25)}{
	\node[scale=0.75] at \i {\textbullet};
	}
\node[scale=0.75] at (1,1.5) {\textcolor{bleu}{$g$}};
\node[scale=0.75] at (2,1.5) {\textcolor{vert}{$g\cdot tat^{-1}$}};
\node[scale=0.75] at (-0.5,1.5) {\textcolor{vert}{$g\cdot ta^{-1}t^{-1}$}};
\node[draw,rectangle, draw=bleu, fill=bleu!10, text width=2em,align=center, rounded corners=2, minimum height=1em,scale=0.4] at (1.25,0.25) {\textcolor{bleu}{$(b,j)$}};
\node[draw,rectangle, draw=vert, fill=vert!10, text width=2em,align=center, rounded corners=2, minimum height=1em,scale=0.4] at (0.5,0.25) {\textcolor{vert}{$(b,j)$}};
\node[draw,rectangle, draw=bleu, fill=bleu!10, text width=2em,align=center, rounded corners=2, minimum height=1em,scale=0.4] at (2.75,0.25) {\textcolor{bleu}{$(a,i)$}};
\node[draw,rectangle, draw=vert, fill=vert!10, text width=2em,align=center, rounded corners=2, minimum height=1em,scale=0.4] at (3.5,0.25) {\textcolor{vert}{$(a,i)$}};
	
\node[scale=1] at (2,-1.5) {(a) $x_{g}=t_1(b,j,a,i) \Rightarrow \left\{\begin{array}{l} x_{g\cdot tat^{-1}}=t_2(a,i)\\ x_{g\cdot ta^{-1}t^{-1}}=t_4(b,j)\end{array}\right.$
};	

\begin{scope}[shift={(7,0)}]
\draw (0,0) -- (2,0) -- (1.5,1) -- (-0.5,1) -- cycle;
\draw (2,0) -- (4,0) -- (3.5,1) -- (1.5,1) -- cycle;
\draw[color=vert,pattern color=vert,pattern=north west lines] (0.75,0.5) -- (1.75,0.5) -- (1.5,1) -- (0.5,1) -- cycle ;
\node[draw,rectangle, draw=vert, fill=vert!10, text width=2em,align=center, rounded corners=2, minimum height=1em,scale=0.6] at (1.15,0.75) {\textcolor{vert}{$(a,i)$}};
\draw (1,0) -- (3,0) -- (3,1.25) -- (1,1.25) -- cycle;
\draw[color=bleu,pattern color=bleu,pattern=north west lines] (2,0.625) rectangle (3,1.25); 
\foreach \i in {(-0.5,1),(1.5,1),(1,1.25)}{
	\node[scale=0.75] at \i {\textbullet};
	}
\node[scale=0.75] at (1,1.5) {\textcolor{bleu}{$g$}};
\node[scale=0.75] at (2,1.5) {\textcolor{vert}{$g\cdot tat^{-1}$}};
\node[scale=0.75] at (-0.5,1.5) {\textcolor{vert}{$g\cdot ta^{-1}t^{-1}$}};	
\node[draw,rectangle, draw=bleu, fill=bleu!10, text width=2em,align=center, rounded corners=2, minimum height=1em,scale=0.6] at (2.5,1) {\textcolor{bleu}{$(a,i)$}};	
\end{scope}	

\begin{scope}[shift={(12,0)}]
\draw (0,0) -- (2,0) -- (1.5,1) -- (-0.5,1) -- cycle;
\draw (2,0) -- (4,0) -- (3.5,1) -- (1.5,1) -- cycle;
\draw (1,0) -- (3,0) -- (3,1.25) -- (1,1.25) -- cycle;
\draw[color=vert,pattern color=vert,pattern=north west lines] (2.75,0.5) -- (3.75,0.5) -- (3.5,1) -- (2.5,1) -- cycle ;
\node[draw,rectangle, draw=vert, fill=vert!10, text width=2em,align=center, rounded corners=2, minimum height=1em,scale=0.6] at (3.15,0.75) {\textcolor{vert}{$(a,i)$}};
\draw[color=bleu,pattern color=bleu,pattern=north west lines] (2,0.625) rectangle (3,1.25); 
\foreach \i in {(-0.5,1),(1.5,1),(1,1.25)}{
	\node[scale=0.75] at \i {\textbullet};
	}
\node[scale=0.75] at (1,1.5) {\textcolor{bleu}{$g$}};
\node[scale=0.75] at (2,1.5) {\textcolor{vert}{$g\cdot tat^{-1}$}};
\node[scale=0.75] at (-0.5,1.5) {\textcolor{vert}{$g\cdot ta^{-1}t^{-1}$}};	
\node[draw,rectangle, draw=bleu, fill=bleu!10, text width=2em,align=center, rounded corners=2, minimum height=1em,scale=0.6] at (2.5,1) {\textcolor{bleu}{$(a,i)$}};	
\end{scope}	

\begin{scope}[shift={(10,-0.5)}]
\node[scale=1] at (1,-1) {(b) $x_{g}=b_1(a,i) \Rightarrow 
\left\{\begin{array}{l}
x_{g\cdot bab^{-1}}=b_1(a,i)\\
or \\
x_{g\cdot ta^{-1}t^{-1}}=b_1(a,i)
\end{array}\right.$}; 
\end{scope}

\end{tikzpicture}
    \end{bigcenter}
    \begin{bigcenter}
    \begin{tikzpicture}[scale=0.7]

\begin{scope}[shift={(24,0)}]
\draw (0,0) -- (2,0) -- (1.5,1) -- (-0.5,1) -- cycle;
\draw (2,0) -- (4,0) -- (3.5,1) -- (1.5,1) -- cycle;
\draw (1,0) -- (3,0) -- (3,1.25) -- (1,1.25) -- cycle;
\draw[color=vert,pattern color=vert,pattern=north west lines] (1.5,0) -- (2,0) -- (1.5,1) -- (0.5,1) -- cycle ;
\draw[color=vert,pattern color=vert,pattern=north west lines] (0,0) -- (0.5,0) -- (-0.25,0.5) -- cycle ; 

\node[draw,rectangle, draw=vert, fill=vert!10, text width=2em,align=center, rounded corners=2, minimum height=1em,scale=0.6] at (1.6,0.6) {\textcolor{vert}{$(a,i)$}};
\node[draw,rectangle, draw=vert, fill=vert!10, text width=2em,align=center, rounded corners=2, minimum height=1em,scale=0.4] at (0.25,0.25) {\textcolor{vert}{$(b,j)$}};
\draw[color=bleu,pattern color=bleu,pattern=north west lines] (1.5,0) -- (3,0) -- (3,1.25) -- (2,1.25) -- cycle ;
\foreach \i in {(-0.5,1),(1.5,1),(1,1.25)}{
	\node[scale=0.75] at \i {\textbullet};
	}
\node[scale=0.75] at (1,1.5) {\textcolor{bleu}{$g$}};
\node[scale=0.75] at (-0.5,1.25) {\textcolor{vert}{$g\cdot ta^{-1}t^{-1}$}};	
\node[draw,rectangle, draw=bleu, fill=bleu!10, text width=2em,align=center, rounded corners=2, minimum height=1em,scale=0.6] at (2.5,0.5) {\textcolor{bleu}{$(a,i)$}};
	
\node[scale=1] at (1,-1) {(d) $x_{g}=m_1(a,i) \Leftrightarrow  x_{g\cdot ta^{-1}t^{-1}}=m_2(b,j,a,i)$};	
\end{scope}

\begin{scope}[shift={(8,0)}]
\draw (0,0) -- (2,0) -- (1.5,1) -- (-0.5,1) -- cycle;
\draw (2,0) -- (4,0) -- (3.5,1) -- (1.5,1) -- cycle;
\draw (1,0) -- (3,0) -- (3,1.25) -- (1,1.25) -- cycle;
\draw[color=vert,pattern color=vert,pattern=north west lines] (-0.25,0.5) -- (0.75,0.5) -- (0.5,1) -- (-0.5,1) -- cycle ;
\node[draw,rectangle, draw=vert, fill=vert!10, text width=2em,align=center, rounded corners=2, minimum height=1em,scale=0.5] at (0.2,0.75) {\textcolor{vert}{$(a,i)$}};	
\draw[color=bleu,pattern color=bleu,pattern=north west lines] (1,0.625) rectangle (2,1.25); 
\foreach \i in {(-0.5,1),(1.5,1),(1,1.25)}{
	\node[scale=0.75] at \i {\textbullet};
	}
\node[scale=0.75] at (1,1.5) {\textcolor{bleu}{$g$}};
\node[scale=0.75] at (2,1.5) {\textcolor{vert}{$g\cdot tat^{-1}$}};
\node[scale=0.75] at (-0.5,1.5) {\textcolor{vert}{$g\cdot ta^{-1}t^{-1}$}};	
\node[draw,rectangle, draw=bleu, fill=bleu!10, text width=2em,align=center, rounded corners=2, minimum height=1em,scale=0.6] at (1.5,0.95) {\textcolor{bleu}{$(a,i)$}};	
\end{scope}	

\begin{scope}[shift={(14,0)}]
\draw (0,0) -- (2,0) -- (1.5,1) -- (-0.5,1) -- cycle;
\draw (2,0) -- (4,0) -- (3.5,1) -- (1.5,1) -- cycle;
\draw (1,0) -- (3,0) -- (3,1.25) -- (1,1.25) -- cycle;
\draw[color=vert,pattern color=vert,pattern=north west lines] (1.75,0.5) -- (2.75,0.5) -- (2.5,1) -- (1.5,1) -- cycle ;
\node[draw,rectangle, draw=vert, fill=vert!10, text width=2em,align=center, rounded corners=2, minimum height=1em,scale=0.5] at (2.2,0.75) {\textcolor{vert}{$(a,i)$}};
\draw[color=bleu,pattern color=bleu,pattern=north west lines] (1,0.625) rectangle (2,1.25); 
\foreach \i in {(-0.5,1),(1.5,1),(1,1.25)}{
	\node[scale=0.75] at \i {\textbullet};
	}
\node[scale=0.75] at (1,1.5) {\textcolor{bleu}{$g$}};
\node[scale=0.75] at (2,1.5) {\textcolor{vert}{$g\cdot tat^{-1}$}};
\node[scale=0.75] at (-0.5,1.5) {\textcolor{vert}{$g\cdot ta^{-1}t^{-1}$}};	
\node[draw,rectangle, draw=bleu, fill=bleu!10, text width=2em,align=center, rounded corners=2, minimum height=1em,scale=0.6] at (1.5,0.95) {\textcolor{bleu}{$(a,i)$}};	
\end{scope}	

\begin{scope}[shift={(12,-0.5)}]
\node[scale=1] at (1,-1) {(c) $x_{g}=b_3(a,i) \Rightarrow 
\left\{\begin{array}{l}
x_{g\cdot bab^{-1}}=b_3(a,i)\\
or \\
x_{g\cdot ba^{-1}t^{-1}}=b_3(a,i)
\end{array}\right.$}; 
\end{scope}

\end{tikzpicture}
    \end{bigcenter}

 \item \label{rule.substitution} 
 if a shaded area carries information $(a,i)$, then the shaded areas immediately below it must carry information $(a_1,1),\dots,(a_{|\sigma(a)|},|\sigma(a)|)$ in this order, with $\sigma(a)=a_1\dots a_{|\sigma(a)|}$. 
  \begin{bigcenter}
  \begin{tikzpicture}[scale=0.5]

\draw (-10,2) rectangle (10,3);
\foreach \i in {-8,-6,...,6,8}{
	\draw (\i,2) -- (\i,3);
	} 
\draw (-10,1) rectangle (11,2);
\foreach \i in {-9,...,10}{
	\draw (\i,1) -- (\i,2);
	} 
\draw[color=bleu,pattern color=bleu,pattern=north west lines] (-9,2.5) rectangle (9,3); 
\node[draw,rectangle, draw=bleu, fill=bleu!10,align=center, rounded corners=2, minimum height=1em,scale=0.75] at (0,3.25) {\textcolor{bleu}{\boldmath $(a,i)$}};
\draw[color=vert,pattern color=vert,pattern=north west lines] (-8.75,1) rectangle (-5.75,1.5);  
\node[draw,rectangle, draw=vert, fill=vert!10,align=center, rounded corners=2, minimum height=1em,scale=0.75] at (-7.25,0.75) {\textcolor{vert}{\boldmath $(a_1,1)$}};
\draw[color=rouge,pattern color=rouge,pattern=north west lines] (-5.25,1) rectangle (-0.25,1.5); 
\node[draw,rectangle, draw=rouge, fill=rouge!10,align=center, rounded corners=2, minimum height=1em,scale=0.75] at (-2.75,0.75) {\textcolor{rouge}{\boldmath $(a_2,2)$}}; 
\draw[color=violet,pattern color=violet,pattern=north west lines] (0.25,1) rectangle (4.25,1.5); 
\node[draw,rectangle, draw=violet, fill=violet!10,align=center, rounded corners=2, minimum height=1em,scale=0.75] at (2.25,0.75) {\textcolor{violet}{\boldmath $(a_3,3)$}}; 
\draw[color=orange,pattern color=orange,pattern=north west lines] (4.75,1) rectangle (10.25,1.5); 
\node[draw,rectangle, draw=orange, fill=orange!10,align=center, rounded corners=2, minimum height=1em,scale=0.75] at (7.5,0.75) {\textcolor{orange}{\boldmath $(a_4,4)$}}; 
 
\end{tikzpicture}
  \end{bigcenter}
  
\end{enumerate}

\begin{figure}[!ht]
  \begin{bigcenter}
  \includegraphics[width=1.2\linewidth]{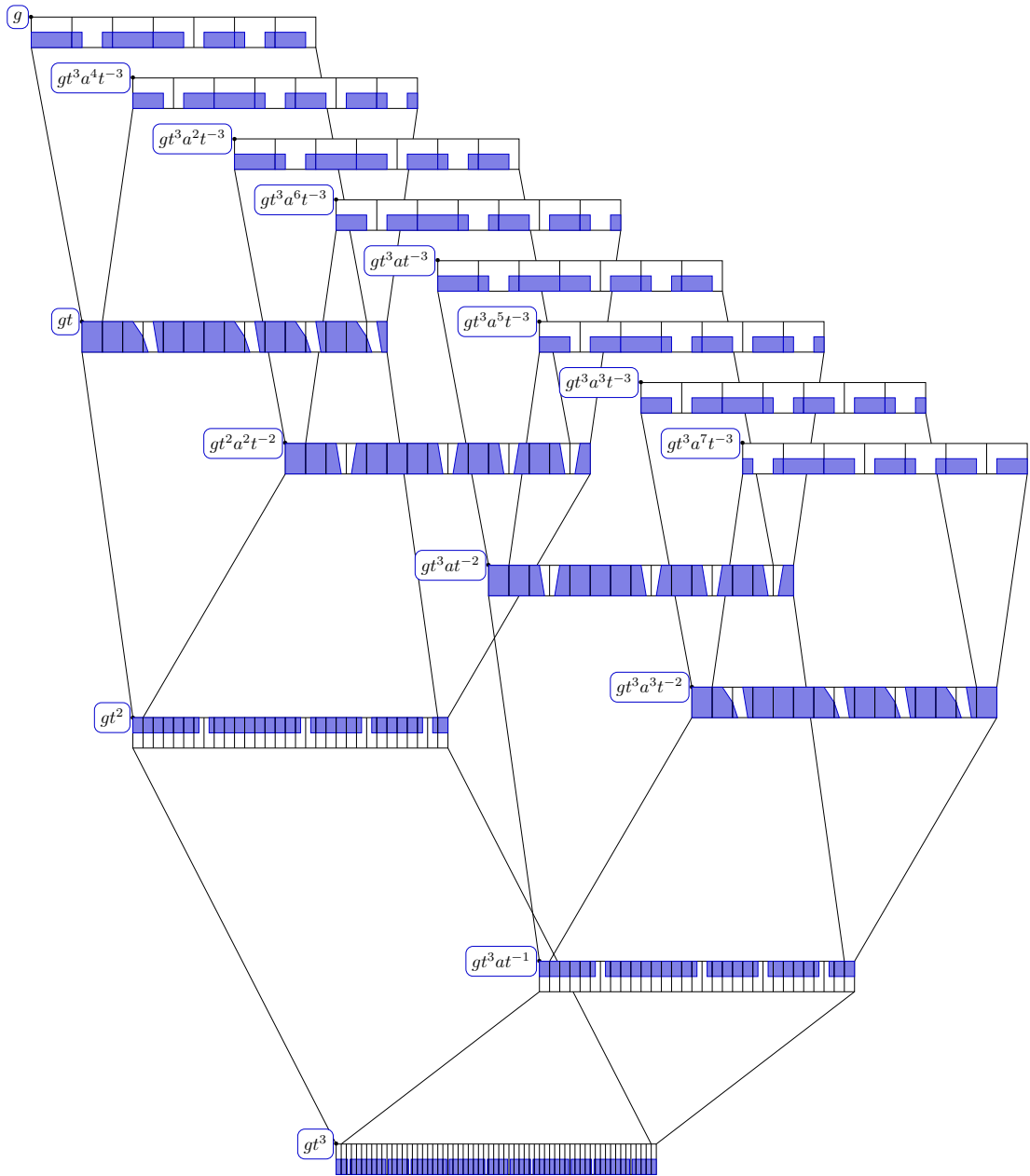}
  \end{bigcenter}
  \caption{Example of configuration $x$ in $X_\sigma$.}
  \label{figure.example_configuration_X_sigma_1}
  \end{figure}

\begin{figure}[!ht]
  \begin{bigcenter}
  \includegraphics[width=1.2\linewidth]{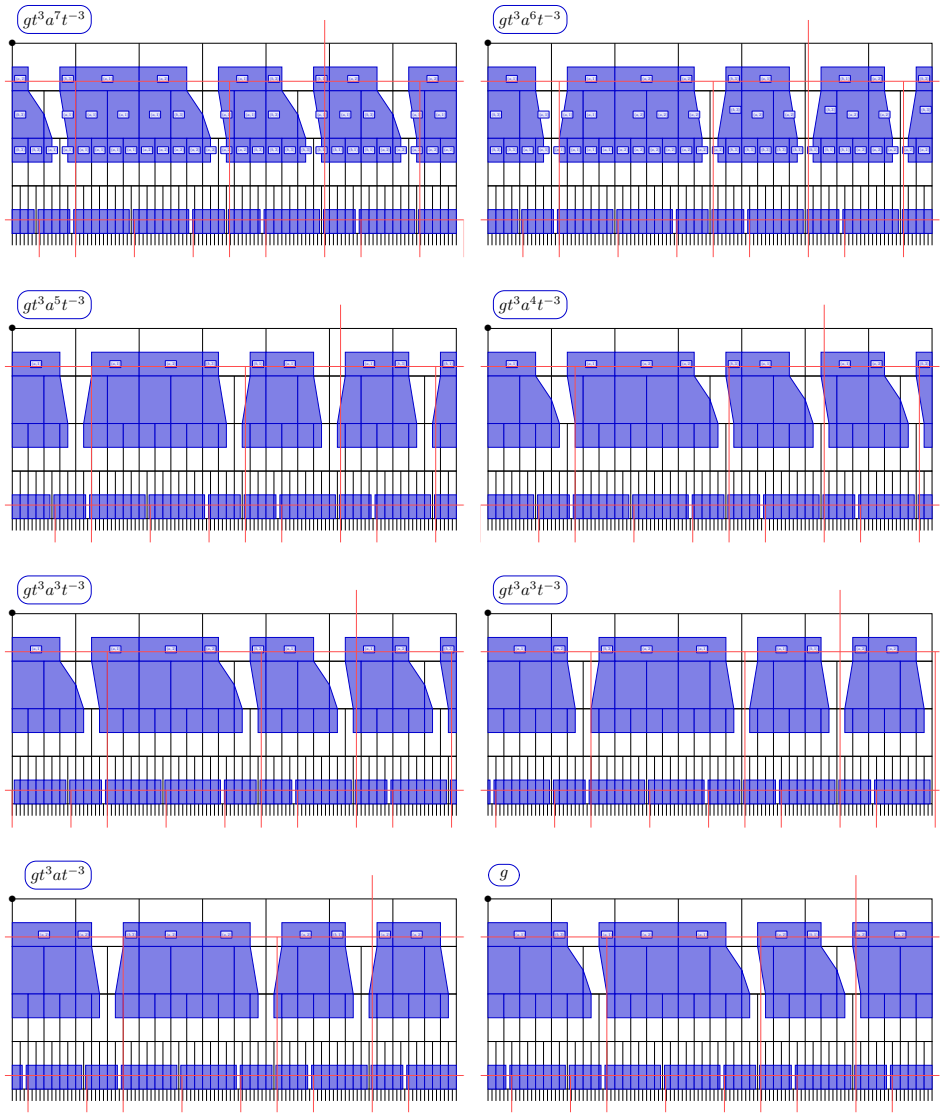}
  \end{bigcenter}
  \caption{The configuration $x$ of Figure~\ref{figure.example_configuration_X_sigma_1} pictured sheet by sheet this time. The $\sigma$-tiling corresponding to $x$ appears in red.}
  \label{figure.example_configuration_X_sigma_2}
  \end{figure}

We have now completely defined an SFT $X_\sigma\subset \mathcal{R}^{BS(1,2)}$. In the next section, we prove that a $\sigma$-tiling corresponds to a configuration in $X_\sigma$ and vice-versa, which is the statement of~Theorem~\ref{theorem.sigma_tilings_and_X_sigma}. The example depicted on Figures~\ref{figure.example_configuration_X_sigma_1} and~\ref{figure.example_configuration_X_sigma_2} should help to understand configurations in $X_\sigma$.

\newpage

\newpage

\newpage





\section{The SFT $X_\sigma$ encodes $\sigma$-tilings}
\label{section.X_sigma_SFT_cover_Y_sigma}

Remind that $\pi:\mathcal{R}\to\A$ is the mapping defined by
 \begin{align*}
 \pi\left((r,(a,i))\right) &= a\\
 \pi\left((r,(a,i),(b,j))\right) &= b
 \end{align*}
that only keeps from a letter $r\in\mathcal{R}$ the letter $a\in\A$ of the $a$-tile it encodes. Denote by $p_{\ell,a,i}$ the pattern with support $\{1,a,\dots,a^{\ell-1}\}$ such that 
\begin{align*}
 (p_{\ell,a,i})_1&=(b_1,(a,i));\\
 (p_{\ell,a,i})_{a^j}&=(b_2,(a,i)) \text{ for every }1\leq j\leq \ell-2;\\
 (p_{\ell,a,i})_{a^{\ell-1}}&=(b_3,(a,i)).
\end{align*}

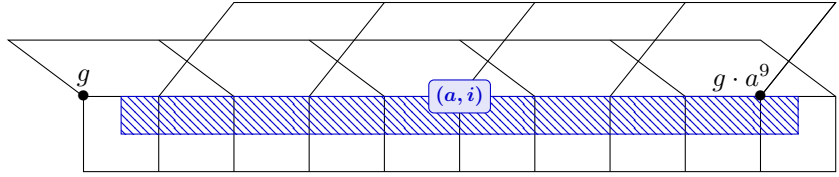
\begin{figure}[!ht]
\begin{bigcenter}
\begin{tikzpicture}[scale=1]
 
\draw (-5,1) rectangle (5,2);
\foreach \i in {-4,...,4}{
	\draw (\i,1) -- (\i,2);
	} 

\draw (-5,2) -- (-6,2.75) -- (4,2.75) -- (5,2);
\foreach \i in {-3,-1,1,3}{
	\draw (\i,2) -- (\i-1,2.75);
	}
	
\draw (-3,2) -- (-4,2) -- (-3,3.25) -- (5,3.25) -- (4,2) -- (3,2);
\foreach \i in {-2,0,2,4}{
	\draw (\i,2) -- (\i+1,3.25);
	}
	
\draw[color=bleu,pattern color=bleu,pattern=north west lines] (-4.5,1.5) rectangle (4.5,2); 
\node[draw,rectangle, draw=bleu, fill=bleu!10,align=center, rounded corners=2, minimum height=1em,scale=0.75] at (0,2) {\textcolor{bleu}{\boldmath $(a,i)$}};

\node[scale=1] at (-5,2) {$\bullet$};
\node[scale=1] at (-5,2.25) {$g$};
\node[scale=1] at (4,2) {$\bullet$};
\node[scale=1] at (3.75,2.25) {$g\cdot a^9$};

\end{tikzpicture}
\end{bigcenter}
\caption{The pattern $p_{10,a,i}$ on the bottom of the $g$-cone $\mathcal{C}_{g,10,1}$.}
\label{figure.pattern_pnai_and_g_cone}
\end{figure}

The pattern $p_{\ell,a,i}$ is locally admissible if and only if $2^{h-1}(w_a-3)+2\leq \ell\leq 2^h(2w_a-1)$. 
For $\ell,m\in\N$, denote 
$$\mathcal{C}_{g,\ell,m}:=\left\{ g\cdot a^{i}t^{-j}\mid 0\leq i \leq \ell-1,0\leq j \leq m\right\}$$
the \define{$g$-cone} of base $\ell$ and height $m$.
 
\begin{lemma}\label{lemma.bottom_forces_rectangle}
 Let $x$ and $y$ be two configurations in $X_\sigma$ such that $p_{\ell,a,i}$ appears in $x$ and $y$ in position $g\in BS(1,2)$. Then there exists $m\in\{h,h-1\}$ such that $x_{g'}=y_{g'}$ for every $g'\in\mathcal{C}_{g,\ell,m}$.
\end{lemma}
 
\begin{proof}
We only need to determine that the two configurations $x$ and $y$ coincide on some $\mathcal{C}_{g,\ell,m}$ with $m$ big enough, and the fact that $m\in\{h,h-1\}$ will immediately follow from the local rule~\ref{rule.dimensions} given on page~\pageref{rule.dimensions}. We distinguish between several cases that are summarized on Figure~\ref{figure.cases_bottom_fixes_all}. With these rules we reconstruct from bottom to top the shaded areas with bottom $p_{\ell,a,i}$, by updating the value of $\ell$ at each step as stated on the picture. Since these rules are totally deterministic, and since the case $w_a-1 \leq \ell\leq 2w_a-1$ is eventually reached, we necessarily obtain a pattern with support $\mathcal{C}_{g,\ell,m}$ with only finite shaded areas on all its sheets. These shaded areas all have the same height $m$ (and $m\geq 2$ by item~2 of Proposition~\ref{proposition.unique_size_property}), but their top widths may differ. 

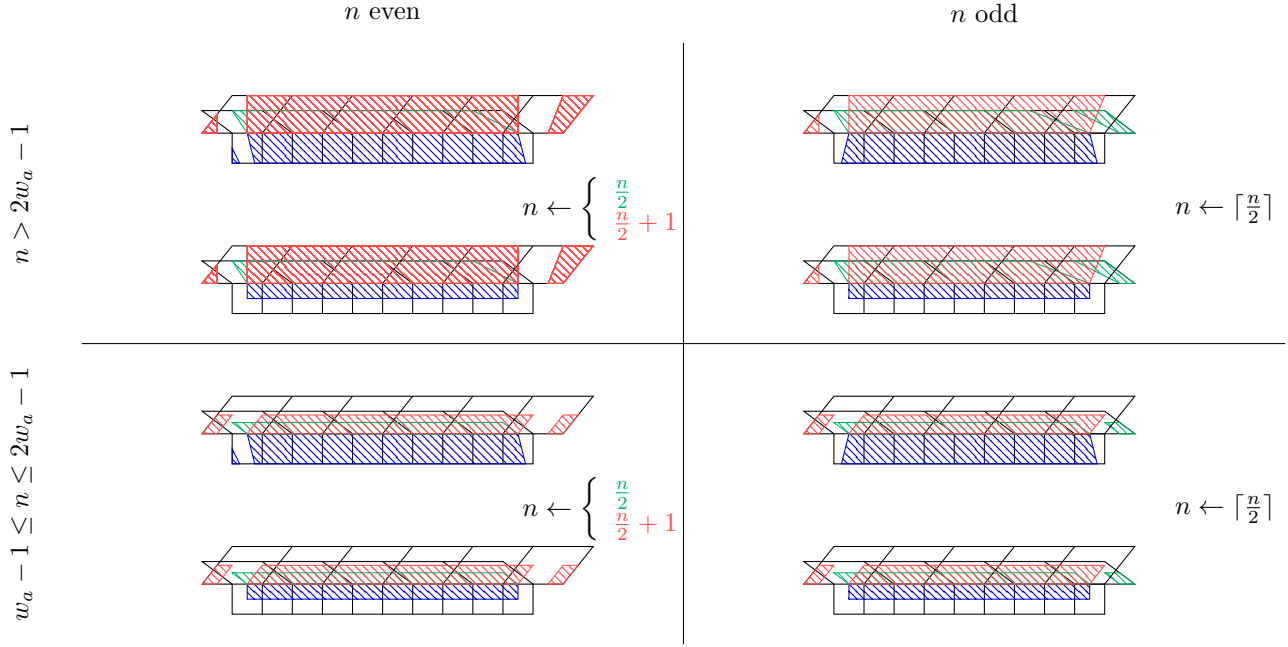
\begin{figure}[!ht]
\begin{bigcenter}
\begin{tikzpicture}[scale=0.4]

\draw (-20,0) -- (20,0);
\draw (0,-10) -- (0,10);
\draw node at (-10,11) {\scalebox{1}{$n$ even}};
\draw node at (10,11) {\scalebox{1}{$n$ odd}};

\draw node at (-22,-5) {\rotatebox{90}{\scalebox{1}{$w_a-1\leq n\leq 2w_a-1$}}};
\draw node at (-22,5) {\rotatebox{90}{\scalebox{1}{$n> 2w_a-1$}}};
 
\begin{scope}[shift={(10,0)}]

\draw (-5,1) rectangle (4,2);
\foreach \i in {-4,...,3}{
	\draw (\i,1) -- (\i,2);
	} 
\draw[color=bleu,pattern color=bleu,pattern=north west lines] (-4.5,1.5) rectangle (3.5,2); 

\draw (-5,2) -- (-6,2.75) -- (4,2.75) -- (5,2) -- (4,2);
\foreach \i in {-3,-1,1,3}{
	\draw (\i,2) -- (\i-1,2.75);
	} 	 
\draw[color=vert,pattern color=vert,pattern=north west lines] (-4.5,2) -- (-5,2.75) -- (1.5,2.75) --  (3.5,2) -- cycle; 
\draw[color=vert,pattern color=vert,pattern=north west lines] (4.5,2) -- (3,2.75) -- (4,2.75) --  (5,2) -- cycle; 

\draw (-5,2) -- (-6,2) -- (-5,3.25) -- (5,3.25) -- (4,2) -- (3,2);
\foreach \i in {-4,-2,0,2}{
	\draw (\i,2) -- (\i+1,3.25);
	} 
\draw[color=rouge,pattern color=rouge,pattern=north west lines] (-4.5,2) -- (-4.5,3.25) -- (4,3.25) --  (3.5,2) -- cycle;
\draw[color=rouge,pattern color=rouge,pattern=north west lines] (-5.5,2) -- (-5.5,2.625) --  (-6,2) -- cycle; 

\draw node at (8,4.5) {\scalebox{1}{$n\leftarrow \lceil\frac{n}{2}\rceil$}};
\end{scope}

\begin{scope}[shift={(10,5)}]

\draw (-5,1) rectangle (4,2);
\foreach \i in {-4,...,3}{
	\draw (\i,1) -- (\i,2);
	} 
\draw[color=bleu,pattern color=bleu,pattern=north west lines] (-4.75,1) -- (3.75,1) -- (3.5,2) -- (-4.5,2) -- cycle; 

\draw (-5,2) -- (-6,2.75) -- (4,2.75) -- (5,2) -- (4,2);
\foreach \i in {-3,-1,1,3}{
	\draw (\i,2) -- (\i-1,2.75);
	} 	
\draw[color=vert,pattern color=vert,pattern=north west lines] (-4.5,2) -- (-5,2.75) -- (1.5,2.75) --  (3.5,2) -- cycle; 
\draw[color=vert,pattern color=vert,pattern=north west lines] (4.5,2) -- (3,2.75) -- (4,2.75) --  (5,2) -- cycle; 

\draw (-5,2) -- (-6,2) -- (-5,3.25) -- (5,3.25) -- (4,2) -- (3,2);
\foreach \i in {-4,-2,0,2}{
	\draw (\i,2) -- (\i+1,3.25);
	} 
\draw[color=rouge,pattern color=rouge,pattern=north west lines] (-4.5,2) -- (-4.5,3.25) -- (4,3.25) --  (3.5,2) -- cycle;
\draw[color=rouge,pattern color=rouge,pattern=north west lines] (-5.5,2) -- (-5.5,2.625) --  (-6,2) -- cycle; 
\end{scope}

\begin{scope}[shift={(-10,0)}]

\draw (-5,1) rectangle (5,2);
\foreach \i in {-4,...,4}{
	\draw (\i,1) -- (\i,2);
	} 
\draw[color=bleu,pattern color=bleu,pattern=north west lines] (-4.5,1.5) rectangle (4.5,2);

\draw (-5,2) -- (-6,2.75) -- (4,2.75) -- (5,2);
\foreach \i in {-3,-1,1,3}{
	\draw (\i,2) -- (\i-1,2.75);
	}
\draw[color=vert,pattern color=vert,pattern=north west lines] (-4.5,2) -- (-5,2.75) -- (3,2.75) --  (4.5,2) -- cycle;

\draw (-5,2) -- (-6,2) -- (-5,3.25) -- (7,3.25) -- (6,2) -- (5,2);
\foreach \i in {-4,-2,0,2,4}{
	\draw (\i,2) -- (\i+1,3.25);
\draw[color=rouge,pattern color=rouge,pattern=north west lines] (-4.5,2) -- (-4.5,3.25) -- (4.5,3.25) --  (4.5,2) -- cycle;
\draw[color=rouge,pattern color=rouge,pattern=north west lines] (-5.5,2) -- (-5.5,2.625) --  (-6,2) -- cycle; 
\draw[color=rouge,pattern color=rouge,pattern=north west lines] (6,2) -- (7,3.25) -- (6,3.25) --  (5.5,2) -- cycle; 
	} 

\draw node at (7.5,4.5) {\scalebox{1}{$n\leftarrow
\left\{\begin{array}{l}                                                         \textcolor{vert}{\frac{n}{2}}\\
\textcolor{rouge}{\frac{n}{2}+1}
\end{array}\right.
$}};	
\end{scope}

\begin{scope}[shift={(-10,5)}]

\draw (-5,1) rectangle (5,2);
\foreach \i in {-4,...,4}{
	\draw (\i,1) -- (\i,2);
	} 
\draw[color=bleu,pattern color=bleu,pattern=north west lines] (-4.25,1) -- (-4.5,2) -- (4.5,2) -- (4.75,1) -- cycle;
\draw[color=bleu,pattern color=bleu,pattern=north west lines] (-5,1) -- (-5,1.5) -- (-4.75,1) -- cycle;

\draw (-5,2) -- (-6,2.75) -- (4,2.75) -- (5,2);
\foreach \i in {-3,-1,1,3}{
	\draw (\i,2) -- (\i-1,2.75);
	} 	 
\draw[color=vert,pattern color=vert,pattern=north west lines] (-4.5,2) -- (-5,2.75) -- (3,2.75) --  (4.5,2) -- cycle;

\draw (-5,2) -- (-6,2) -- (-5,3.25) -- (7,3.25) -- (6,2) -- (5,2);
\foreach \i in {-4,-2,0,2,4}{
	\draw (\i,2) -- (\i+1,3.25);
\draw[color=rouge,pattern color=rouge,pattern=north west lines] (-4.5,2) -- (-4.5,3.25) -- (4.5,3.25) --  (4.5,2) -- cycle;
\draw[color=rouge,pattern color=rouge,pattern=north west lines] (-5.5,2) -- (-5.5,2.625) --  (-6,2) -- cycle; 
\draw[color=rouge,pattern color=rouge,pattern=north west lines] (6,2) -- (7,3.25) -- (6,3.25) --  (5.5,2) -- cycle; 
	} 
\end{scope}

\begin{scope}[shift={(10,-10)}]

\draw (-5,1) rectangle (4,2);
\foreach \i in {-4,...,3}{
	\draw (\i,1) -- (\i,2);
	} 
\draw[color=bleu,pattern color=bleu,pattern=north west lines] (-4.5,1.5) rectangle (3.5,2);

\draw (-5,2) -- (-6,2.75) -- (4,2.75) -- (5,2) -- (4,2);
\foreach \i in {-3,-1,1,3}{
	\draw (\i,2) -- (\i-1,2.75);
	} 	 
\draw[color=vert,pattern color=vert,pattern=north west lines] (-4.5,2) -- (-5,2.375) -- (3,2.375) --  (3.5,2) -- cycle; 
\draw[color=vert,pattern color=vert,pattern=north west lines] (4.5,2) -- (4,2.375) -- (4.5,2.375) --  (5,2) -- cycle; 
	
\draw (-5,2) -- (-6,2) -- (-5,3.25) -- (5,3.25) -- (4,2) -- (3,2);
\foreach \i in {-4,-2,0,2}{
	\draw (\i,2) -- (\i+1,3.25);
	} 
\draw[color=rouge,pattern color=rouge,pattern=north west lines] (-4.5,2) -- (-4,2.625) -- (4,2.625) --  (3.5,2) -- cycle; 
\draw[color=rouge,pattern color=rouge,pattern=north west lines] (-5.5,2) -- (-5,2.625) -- (-5.5,2.625) --  (-6,2) -- cycle;

\draw node at (8,4.5) {\scalebox{1}{$n\leftarrow \lceil\frac{n}{2}\rceil$}};
\end{scope}

\begin{scope}[shift={(10,-5)}]

\draw (-5,1) rectangle (4,2);
\foreach \i in {-4,...,3}{
	\draw (\i,1) -- (\i,2);
	} 
\draw[color=bleu,pattern color=bleu,pattern=north west lines] (-4.75,1) -- (3.75,1) -- (3.5,2) -- (-4.5,2) -- cycle; 

\draw (-5,2) -- (-6,2.75) -- (4,2.75) -- (5,2) -- (4,2);
\foreach \i in {-3,-1,1,3}{
	\draw (\i,2) -- (\i-1,2.75);
	} 	
\draw[color=vert,pattern color=vert,pattern=north west lines] (-4.5,2) -- (-5,2.375) -- (3,2.375) --  (3.5,2) -- cycle; 
\draw[color=vert,pattern color=vert,pattern=north west lines] (4.5,2) -- (4,2.375) -- (4.5,2.375) --  (5,2) -- cycle; 
	
\draw (-5,2) -- (-6,2) -- (-5,3.25) -- (5,3.25) -- (4,2) -- (3,2);
\foreach \i in {-4,-2,0,2}{
	\draw (\i,2) -- (\i+1,3.25);
	} 
\draw[color=rouge,pattern color=rouge,pattern=north west lines] (-4.5,2) -- (-4,2.625) -- (4,2.625) --  (3.5,2) -- cycle; 
\draw[color=rouge,pattern color=rouge,pattern=north west lines] (-5.5,2) -- (-5,2.625) -- (-5.5,2.625) --  (-6,2) -- cycle; 
\end{scope}

\begin{scope}[shift={(-10,-10)}]

\draw (-5,1) rectangle (5,2);
\foreach \i in {-4,...,4}{
	\draw (\i,1) -- (\i,2);
	} 
\draw[color=bleu,pattern color=bleu,pattern=north west lines] (-4.5,1.5) rectangle (4.5,2);

\draw (-5,2) -- (-6,2.75) -- (4,2.75) -- (5,2);
\foreach \i in {-3,-1,1,3}{
	\draw (\i,2) -- (\i-1,2.75);
	} 
\draw[color=vert,pattern color=vert,pattern=north west lines] (-4.5,2) -- (-5,2.375) -- (4,2.375) --  (4.5,2) -- cycle;
	
\draw (-5,2) -- (-6,2) -- (-5,3.25) -- (7,3.25) -- (6,2) -- (5,2);
\foreach \i in {-4,-2,0,2,4}{
	\draw (\i,2) -- (\i+1,3.25);
	} 
\draw[color=rouge,pattern color=rouge,pattern=north west lines] (-4.5,2) -- (-4,2.625) -- (5,2.625) --  (4.5,2) -- cycle;
\draw[color=rouge,pattern color=rouge,pattern=north west lines] (-5.5,2) -- (-5,2.625) -- (-5.5,2.625) --  (-6,2) -- cycle; 
\draw[color=rouge,pattern color=rouge,pattern=north west lines] (6,2) -- (6.5,2.625) -- (6,2.625) --  (5.5,2) -- cycle;  

\draw node at (7.5,4.5) {\scalebox{1}{$n\leftarrow
\left\{\begin{array}{l}                                                         \textcolor{vert}{\frac{n}{2}}\\
\textcolor{rouge}{\frac{n}{2}+1}
\end{array}\right.
$}};
\end{scope}

\begin{scope}[shift={(-10,-5)}]

\draw (-5,1) rectangle (5,2);
\foreach \i in {-4,...,4}{
	\draw (\i,1) -- (\i,2);
	} 
\draw[color=bleu,pattern color=bleu,pattern=north west lines] (-4.25,1) -- (-4.5,2) -- (4.5,2) -- (4.75,1) -- cycle;
\draw[color=bleu,pattern color=bleu,pattern=north west lines] (-5,1) -- (-5,1.5) -- (-4.75,1) -- cycle;

\draw (-5,2) -- (-6,2.75) -- (4,2.75) -- (5,2);
\foreach \i in {-3,-1,1,3}{
	\draw (\i,2) -- (\i-1,2.75);
	} 
\draw[color=vert,pattern color=vert,pattern=north west lines] (-4.5,2) -- (-5,2.375) -- (4,2.375) --  (4.5,2) -- cycle;
	
\draw (-5,2) -- (-6,2) -- (-5,3.25) -- (7,3.25) -- (6,2) -- (5,2);
\foreach \i in {-4,-2,0,2,4}{
	\draw (\i,2) -- (\i+1,3.25);
	} 
\draw[color=rouge,pattern color=rouge,pattern=north west lines] (-4.5,2) -- (-4,2.625) -- (5,2.625) --  (4.5,2) -- cycle;
\draw[color=rouge,pattern color=rouge,pattern=north west lines] (-5.5,2) -- (-5,2.625) -- (-5.5,2.625) --  (-6,2) -- cycle; 
\draw[color=rouge,pattern color=rouge,pattern=north west lines] (6,2) -- (6.5,2.625) -- (6,2.625) --  (5.5,2) -- cycle; 
\end{scope} 
 
\end{tikzpicture}
\end{bigcenter}
\caption{The different cases to construct a shaded area starting from its bottom.}
\label{figure.cases_bottom_fixes_all}
\end{figure}


Finally since local rules~\ref{rule.matching_rules} and~\ref{rule.synchro_two_sheets} force a connected shaded component to wear the same information $(a,i)$, all letters $x_{g'}$ and $y_{g'}$ for $g'\in\mathcal{C}_{g,\ell,m}$ are totally determined by the fact that pattern $p_{\ell,a,i}$ appears in $x$ and $y$ in position $g$, and the lemma is proven.
\end{proof}
 
Lemma~\ref{lemma.bottom_forces_rectangle} ensures that a bottom pattern $p_{n,a,i}$ that appears in some position $g\in BS(1,2)$ forces the entire $g$-cone $\mathcal{C}_{g,n,h}$ or $\mathcal{C}_{g,n,h-1}$ above it.

\medskip

Define $\pi:\mathcal{R}\to\A$ the mapping given by
 \begin{align*}
 \pi\left((r,(a,i))\right) &= a\\
 \pi\left((r,(a,i),(b,j))\right) &= b
 \end{align*}
that only keeps from a letter $r\in\mathcal{R}$ the letter $a\in\A$ of the $a$-tile it encodes.

\begin{theorem}\label{theorem.sigma_tilings_and_X_sigma}
  If $\tau$ is a $\sigma$-tiling, then there exists a configuration $x\in X_\sigma$ such that $\tau(\Phi(gat)) = \pi(x_g)$ for all $g\in BS(1, 2)$.
  Reciprocally, for every configuration $x\in X_\sigma$, one can associate a $\sigma$-tiling $\tau$ such that $\tau(\Phi(gat)) = \pi(x_g)$ for all $g\in BS(1, 2)$.
 \end{theorem}
 
Before going into the proof of Theorem~\ref{theorem.sigma_tilings_and_X_sigma}, remind that for $g\in BS(1,2)$, the $\Phi(g)$-box is the rectangle $\Phi(g)+[0,2^{-i}[\times ]-1,0]$ in $\mathbb{R}^2$.

\begin{remark}\label{remark.at_most_one_vertical_object}
The assumption that $v(a)\geq3$ for every letter $a$ implies that at most one of the $\Phi(g)$ and $\Phi(ga)$-boxes may contain one object among vertical line, cross or $T$.
\end{remark} 

\begin{remark}\label{remark.at_most_one_horizontal_object}
The assumption that $\log_2(\lambda)\geq3$ implies that at most one of the $\Phi(g)$ and $\Phi(gt)$-boxes may contain one object among horizontal line, cross or $T$.
\end{remark} 

\begin{remark}\label{remark.intersection_Phi_boxes}
Let $g\in BS(1,2)$ and $k\in\N$. Then the rightmost $\frac{2^k-1}{2^k}$ part of the $\Phi(g)$-box is the leftmost $\frac{2^k-1}{2^k}$ part of the $\Phi(gt^kat^{-k})$-box.
\end{remark}

With all these remarks in mind, we are ready to prove Theorem~\ref{theorem.sigma_tilings_and_X_sigma}. For clarity purpose, we split the theorem into two lemmas.
 
\begin{lemma}\label{lemma.tiling_to_X_sigma}
  If $\tau$ is a $\sigma$-tiling, then there exists a configuration $x\in X_\sigma$ such that $\tau(\Phi(gat)) = \pi(x_g)$ for all $g\in BS(1, 2)$.
 \end{lemma}

\begin{proof}[Proof of Lemma~\ref{lemma.tiling_to_X_sigma}]
Let $\tau$ be a $\sigma$-tiling. We define a configuration $x\in\mathcal{R}^{BS(1,2)}$, and then prove that $x\in X_\sigma$. First, every group element $g\in BS(1,2)$ is associated the unique letter $x_g$ in $\mathcal{R}$ such that $\pi(x_g)=\tau\left(\Phi(gat)\right)$. We now need to determine the letter $x_g\in\mathcal{R}$, i.e. the type of letter, the letter $a\in\A$ and the index $i$ (plus maybe another letter $b$ and another index $j$).

We first define all the $x_g$ such that $x_g\in\mathcal{T}$. To do so, for every group element $g$, we check whether the $\Phi(g)$, $\Phi(ga)$, $\Phi(gt)$ and $\Phi(gat)$-boxes are in one of the following cases (thanks to Remark~\ref{remark.at_most_one_vertical_object}, this list is exhaustive):
\begin{enumerate}
 \item if the $\Phi(g)$-box contains a cross or a $T$, and the $\Phi(gt)$-box is empty, set $x_g\in t_1$;
 \item if the $\Phi(g)$-box contains a cross or a $T$, and the $\Phi(gt)$-box contains a vertical line, set $x_g\in t_2$;
 \item if the $\Phi(g)$-box contains an horizontal line and the $\Phi(ga)$-box contains an horizontal line, set $x_g\in t_3$;
 \item if the $\Phi(g)$-box contains an horizontal line, the $\Phi(ga)$-box contains a cross or a $T$, and the $\Phi(gat)$-box is empty, set $x_g\in t_3$;
 \item if the $\Phi(g)$-box contains an horizontal line, the $\Phi(ga)$-box contains a cross or a $T$ and the $\Phi(gat)$-box contains a vertical line, set $x_g\in t_4$.
\end{enumerate}
\begin{bigcenter}
\begin{tikzpicture}[scale=0.7]

 \clip (-1.5,-4) rectangle (19,3);
 
 \begin{scope}[shift={(0,0)}]
\draw[color=vert,very thick] (0,0.5) -- (2,0.5);
\draw[color=vert,very thick] (1.7,0.5) -- (1.7,0);
\draw[color=vert,very thick,dotted] (1.7,1) -- (1.7,0.5);
\draw (0,0) rectangle (2,1);
\draw (0,-1) rectangle (1,0);
\node[scale=0.75] at (0.75,2) {$x_g\in t_1=\lettertone$};
\node[scale=0.75,rotate=90] at (1,1.325) {$\Leftrightarrow$};
\foreach \i in {(0,0),(1,0),(2,0),(0,1),(2,1)}{
	\node[scale=0.5] at \i {\textbullet};
	}
\node[anchor=east,scale=0.75] at (0,1) {$\Phi(g)$};
\node[anchor=east,scale=0.75] at (0,0) {$\Phi(g t)$};
\end{scope}
 
\begin{scope}[shift={(4,0)}]
\draw[color=vert,very thick] (0,0.5) -- (2,0.5);
\draw[color=vert,very thick] (0.3,0.5) -- (0.3,-1);
\draw[color=vert,very thick,dotted] (0.3,1) -- (0.3,0.5);
\draw (0,0) rectangle (2,1);
\draw (0,-1) rectangle (1,0);
\node[scale=0.75] at (0.75,2) {$x_g\in t_2=\letterttwo$};
\node[scale=0.75,rotate=90] at (1,1.325) {$\Leftrightarrow$};
\foreach \i in {(0,0),(1,0),(2,0),(0,1),(2,1)}{
	\node[scale=0.5] at \i {\textbullet};
	}
\node[anchor=east,scale=0.75] at (0,1) {$\Phi(g)$};
\node[anchor=east,scale=0.75] at (0,0) {$\Phi(g t)$};
\end{scope}

\begin{scope}[shift={(8,0)}]
\draw[color=vert,very thick] (0,0.8) -- (4,0.8);
\draw (0,0) rectangle (2,1);
\draw (2,0) rectangle (4,1);
\node[scale=0.75] at (0.75,2) {$x_g\in t_3=\lettertthree$};
\node[scale=0.75,rotate=90] at (1,1.325) {$\Leftrightarrow$};
\foreach \i in {(0,0),(1,0),(2,0),(0,1),(2,1)}{
	\node[scale=0.5] at \i {\textbullet};
	}
\node[anchor=east,scale=0.75] at (0,1) {$\Phi(g)$};
\node[anchor=south,scale=0.75] at (2,1) {$\Phi(g a)$};

\begin{scope}[shift={(0,-2.5)}]
\node[anchor=east,scale=0.75] at (0,1.75) {or};

\draw[color=vert,very thick] (0,0.3) -- (4,0.3);
\draw[color=vert,very thick] (3.4,0.3) -- (3.4,0);
\draw[color=vert,very thick,dotted] (3.4,0.3) -- (3.4,1);
\draw (0,0) rectangle (2,1);
\draw (2,0) rectangle (4,1);
\draw (2,-1) rectangle (3,0);
\foreach \i in {(0,0),(1,0),(2,0),(0,1),(2,1)}{
	\node[scale=0.5] at \i {\textbullet};
	}
\node[anchor=east,scale=0.75] at (0,1) {$\Phi(g)$};
\node[anchor=south,scale=0.75] at (2,1) {$\Phi(g a)$};
\end{scope}
\end{scope}

\begin{scope}[shift={(14,0)}]
\draw[color=vert,very thick] (0,0.3) -- (4,0.3);
\draw[color=vert,very thick] (2.4,0.3) -- (2.4,-1);
\draw[color=vert,very thick,dotted] (2.4,0.3) -- (2.4,1);
\draw (0,0) rectangle (2,1);
\draw (2,0) rectangle (4,1);
\draw (2,-1) rectangle (3,0);
\node[scale=0.75] at (0.75,2) {$x_g\in t_4=\lettertfour$};
\node[scale=0.75,rotate=90] at (1,1.325) {$\Leftrightarrow$};
\foreach \i in {(0,0),(1,0),(2,0),(0,1),(2,1)}{
	\node[scale=0.5] at \i {\textbullet};
	}
\node[anchor=east,scale=0.75] at (0,1) {$\Phi(g)$};
\node[anchor=south,scale=0.75] at (2,1) {$\Phi(g a)$};
\end{scope}

\end{tikzpicture}
\end{bigcenter}

From now on, no other group element $g\in BS(1,2)$ will be designated a letter from $\mathcal{T}$. We now list all group elements $g\in BS(1,2)$ that will be designated a letter from $\mathcal{M}\setminus\{m_3\}$. Again we check whether the $\Phi(g)$, $\Phi(ga)$, $\Phi(gt)$ and $\Phi(gat)$-boxes are in one of the following cases:
\begin{enumerate}
 \item if the $\Phi(g)$ and $\Phi(gt)$-boxes contain a vertical line, set $x_g\in m_1$;
 \item if  the $\Phi(g)$-box contains a vertical line and the $\Phi(gt)$-box is empty, set $x_g\in m_2$;
 \item if the $\Phi(g)$-box is empty and the $\Phi(ga)$ and $\Phi(gat)$-boxes contain a vertical line (note that this is equivalent to have $x_{ga}\in m_1$), set $x_g\in m_4$;
 \item if the $\Phi(g)$-box and $\Phi(gat)$-boxes are empty and the $\Phi(ga)$-box contains a vertical line (note that this is equivalent to have $x_{ga}\in m_2$), set $x_g\in m_5$.
\end{enumerate}
\begin{bigcenter}
 \begin{tikzpicture}[scale=0.7]

 \clip (-1.5,-2) rectangle (19,3);
 
 \begin{scope}[shift={(0,0)}]
\draw[color=vert,very thick] (0.4,1) -- (0.4,-1);
\draw (0,0) rectangle (2,1);
\draw (0,-1) rectangle (1,0);
\node[scale=0.75] at (0.75,2) {$x_g\in m_1=\lettermone$};
\node[scale=0.75,rotate=90] at (1,1.325) {$\Leftrightarrow$};
\foreach \i in {(0,0),(1,0),(2,0),(0,1),(2,1)}{
	\node[scale=0.5] at \i {\textbullet};
	}
\node[anchor=east,scale=0.75] at (0,1) {$\Phi(g)$};
\node[anchor=east,scale=0.75] at (0,0) {$\Phi(g t)$};
\end{scope}
 
\begin{scope}[shift={(4,0)}]
\draw[color=vert,very thick] (1.4,1) -- (1.4,0);
\draw (0,0) rectangle (2,1);
\draw (0,-1) rectangle (1,0);
\node[scale=0.75] at (0.75,2) {$x_g\in m_2=\lettermtwo$};
\node[scale=0.75,rotate=90] at (1,1.325) {$\Leftrightarrow$};
\foreach \i in {(0,0),(1,0),(2,0),(0,1),(2,1)}{
	\node[scale=0.5] at \i {\textbullet};
	}
\node[anchor=east,scale=0.75] at (0,1) {$\Phi(g)$};
\node[anchor=east,scale=0.75] at (0,0) {$\Phi(g t)$};
\end{scope}

\begin{scope}[shift={(8,0)}]
\draw[color=vert,very thick] (2.4,1) -- (2.4,-1);
\draw (0,0) rectangle (2,1);
\draw (2,0) rectangle (4,1);
\draw (2,-1) rectangle (3,0);
\node[scale=0.75] at (0.75,2) {$x_g\in m_4=\lettermfour$};
\node[scale=0.75,rotate=90] at (1,1.325) {$\Leftrightarrow$};
\foreach \i in {(0,0),(1,0),(2,0),(0,1),(2,1)}{
	\node[scale=0.5] at \i {\textbullet};
	}
\node[anchor=east,scale=0.75] at (0,1) {$\Phi(g)$};
\node[anchor=south,scale=0.75] at (2,1) {$\Phi(g a)$};
\end{scope}

\begin{scope}[shift={(14,0)}]
\draw[color=vert,very thick] (3.2,1) -- (3.2,0);
\draw (0,0) rectangle (2,1);
\draw (2,0) rectangle (4,1);
\draw (2,-1) rectangle (3,0);
\node[scale=0.75] at (0.75,2) {$x_g\in m_5=\lettermfive$};
\node[scale=0.75,rotate=90] at (1,1.325) {$\Leftrightarrow$};
\foreach \i in {(0,0),(1,0),(2,0),(0,1),(2,1)}{
	\node[scale=0.5] at \i {\textbullet};
	}
\node[anchor=east,scale=0.75] at (0,1) {$\Phi(g)$};
\node[anchor=south,scale=0.75] at (2,1) {$\Phi(g a)$};
\end{scope}

\end{tikzpicture}
\end{bigcenter}
Now we force every element $g\in BS(1,2)$ which has not yet been designated a letter in $\mathcal{R}$ and which is such that some $g\cdot a^k$ has been designated a letter in $\mathcal{M}\setminus\{m_3\}$ to get letter $m_3$ (we fill in the holes inside $\mathcal{M}$-rows).

\medskip

We finally list all group elements $g\in BS(1,2)$ that will be designated a letter from $\mathcal{B}$. Again we check whether the $\Phi(g)$, $\Phi(ga)$, $\Phi(gt)$ and $\Phi(gat)$-boxes are in one of the following cases:
\begin{enumerate}
 \item if the $\Phi(gt)$-box or the $\Phi(gta)$-box contains a cross, set $x_g\in b_1$;
 \item if the $\Phi(g)$ and $\Phi(ga)$-boxes are empty and the $\Phi(gt)$-box contains an horizontal line or a $T$, set $x_g\in b_2$;
 \item if the $\Phi(g)$-box is empty, the $\Phi(gt)$-box contains an horizontal line or a $T$ and the $\Phi(ga)$-box contains a vertical line, set $x_g\in b_3$.
\end{enumerate}
\begin{bigcenter}
\begin{tikzpicture}[scale=0.7]

 \clip (-2.5,-2) rectangle (18,3);
 
 \begin{scope}[shift={(0,0)}]
\node[scale=0.75] at (0.75,2) {$x_g\in b_1=\letterbone$};
\node[scale=0.75,rotate=90] at (1,1.325) {$\Leftrightarrow$};
\node[scale=0.75] at (1.5,-0.75) {or};

\begin{scope}[shift={(-1,0)}]
\draw[color=vert,very thick] (0.4,0) -- (0.4,-1);
\draw[color=vert,very thick] (0,-0.7) -- (1,-0.7);
\draw (0,0) rectangle (2,1);
\draw (0,-1) rectangle (1,0);
\foreach \i in {(0,0),(1,0),(2,0),(0,1),(2,1)}{
	\node[scale=0.5] at \i {\textbullet};
	}
\node[anchor=east,scale=0.75] at (0,1) {$\Phi(g)$};
\node[anchor=east,scale=0.75] at (0,0) {$\Phi(g t)$};
\end{scope}

\begin{scope}[shift={(2.5,0)}]
\draw[color=vert,very thick] (1.4,0) -- (1.4,-1);
\draw[color=vert,very thick] (1,-0.7) -- (2,-0.7);
\draw (0,0) rectangle (2,1);
\draw (1,-1) rectangle (2,0);
\foreach \i in {(0,0),(1,0),(2,0),(0,1),(2,1)}{
	\node[scale=0.5] at \i {\textbullet};
	}
\node[anchor=east,scale=0.75] at (0,1) {$\Phi(g)$};
\node[anchor=north east,scale=0.75] at (1,0) {$\Phi(g ta)$};
\end{scope}

\end{scope}

\begin{scope}[shift={(7,0)}]
\draw[color=vert,very thick] (0,-0.3) -- (1,-0.3);
\draw[color=vert,very thick,dotted] (0.6,-0.3) -- (0.6,-1);
\draw (0,0) rectangle (2,1);
\draw (2,0) rectangle (4,1);
\draw (0,-1) rectangle (1,0);
\node[scale=0.75] at (0.75,2) {$x_g\in b_2=\letterbtwo$};
\node[scale=0.75,rotate=90] at (1,1.325) {$\Leftrightarrow$};
\foreach \i in {(0,0),(1,0),(2,0),(0,1),(2,1)}{
	\node[scale=0.5] at \i {\textbullet};
	}
\node[anchor=east,scale=0.75] at (0,1) {$\Phi(g)$};
\node[anchor=south,scale=0.75] at (2,1) {$\Phi(g a)$};
\node[anchor=east,scale=0.75] at (0,0) {$\Phi(g t)$};
\end{scope}

\begin{scope}[shift={(13,0)}]
\draw[color=vert,very thick] (3.2,1) -- (3.2,0);
\draw[color=vert,very thick] (0,-0.3) -- (1,-0.3);
\draw[color=vert,very thick,dotted] (0.6,-0.3) -- (0.6,-1);
\draw (0,0) rectangle (2,1);
\draw (2,0) rectangle (4,1);
\draw (0,-1) rectangle (1,0);
\node[scale=0.75] at (0.75,2) {$x_g\in b_3=\letterbthree$};
\node[scale=0.75,rotate=90] at (1,1.325) {$\Leftrightarrow$};
\foreach \i in {(0,0),(1,0),(2,0),(0,1),(2,1)}{
	\node[scale=0.5] at \i {\textbullet};
	}
\node[anchor=east,scale=0.75] at (0,1) {$\Phi(g)$};
\node[anchor=south,scale=0.75] at (2,1) {$\Phi(g a)$};
\node[anchor=east,scale=0.75] at (0,0) {$\Phi(g t)$};
\end{scope}

\end{tikzpicture}
\end{bigcenter}

A carefull observation of all possible cases treated by the three previous lists show that every element $g\in BS(1,2)$ necessarily falls into one of the cases described, so that every element $g\in BS(1,2)$ is now designated a letter from $\mathcal{R}$.

\medskip

To get a configuration $x\in\mathcal{R}^{BS(1,2)}$, it only remains to assign letters from $\A$ and indices to shaded connected areas. They are chosen to be identical for all letters in a same shaded area. Consider one shaded area. The top left corner, of type $t_1$ or $t_2$, appears in some position $g\in BS(1,2)$. The corresponding $\Phi(g)$-box contains a top left corner of a $\sigma$-tile. Denote $a$ the corresponding letter : assign the letter $a$ to every letter in $\mathcal{R}$ belonging to the shaded area. To determine the index $i$, it is enough to look at the $\sigma$-tile immediately above in the $\sigma$-tiling. Assume it is a $a'$-tile, and that our shaded area corresponds to the index $i$ such that $\sigma(a')_i=a$: assign the index $i$ to every letter in $\mathcal{R}$ belonging to the shaded area. For letters in $t_1$ or $m_2$, the additional letter $b$ and index $j$ are chosen to be the letter and index associated with the shaded area immediately to the left.

\bigskip

We now check that this configuration $x\in\mathcal{R}^{BS(1,2)}$ in actually in $X_\sigma$ (local rules are listed on page~\pageref{rule.matching_rules}). By construction, it is straightforward that local rules~\ref{rule.matching_rules},~\ref{rule.type_TMB} and~\ref{rule.substitution} are respected. The definition of $\sigma$-tiles (see page~\pageref{def.sigma_tiles}) implies that local rule~\ref{rule.dimensions} is satisfied. And a careful reading of Remark~\ref{remark.intersection_Phi_boxes} ensures that local rule~\ref{rule.synchro_two_sheets} always holds. Finally, we can associate to every $\sigma$-tiling a valid configuration $x\in X_\sigma$ such that $\tau(\Phi(gat)) = \pi(x_g)$ for all $g\in BS(1, 2)$, so Lemma~\ref{lemma.tiling_to_X_sigma} is proven.
\end{proof}

\begin{lemma}\label{lemma.X_sigma_to_tiling}
For every configuration $x\in X_\sigma$, one can associate a $\sigma$-tiling $\tau$ such that $\tau(\Phi(gat)) = \pi(x_g)$ for all $g\in BS(1, 2)$.
 \end{lemma}

\begin{figure}[!ht]
\begin{bigcenter}
\includegraphics{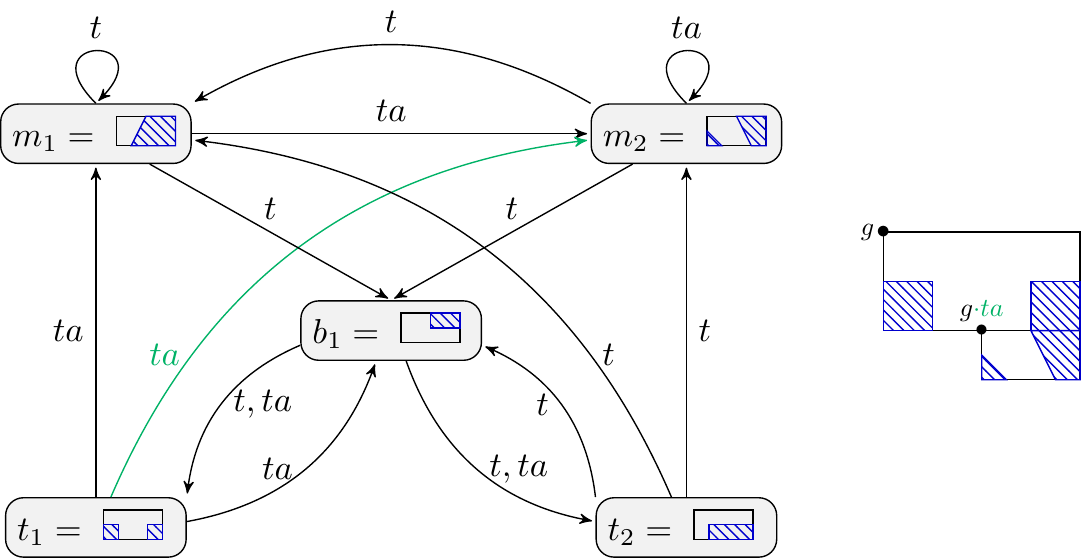}
\end{bigcenter}
\caption{The transducer $\mathcal{D}$ that summarizes how letters from alphabet $\mathcal{R}$ can be arranged to form the left border of shaded areas. The transition in green corresponds to the situation pictured to the right.}
\label{figure.transducer_down}
\end{figure}

Before going into the proof of Lemma~\ref{lemma.X_sigma_to_tiling}, we focus  on the left borders of shaded areas in configurations of $X_\sigma$. Let $\mathcal{D}$ be the transducer depicted on Figure~\ref{figure.transducer_down}. A transition $(q,w,q')$ in $\mathcal{D}$ should be interpreted as: if a letter $q$ appears in some position $g$ in a configuration $x$ from $X_\sigma$, then the letter $q'$ may appear in position $g\cdot w$ in $x$.

 We could have constructed a bigger transducer that would imply all letters from $\mathcal{R}$, but for the purpose of Lemma~\ref{lemma.X_sigma_to_tiling} we only need to focus on letters $t_1,t_2,m_1,m_2$ and $b_1$ since there are the only ones that may appear on the left border of a shaded area.  Note that this transducer $\mathcal{D}$ is non deterministic, but for our purpose the uniqueness of the infinite path described in Proposition~\ref{proposition.path_transducer} is enough.
 
\begin{proposition}\label{proposition.path_transducer}
Let $x$ be a configuration of $X_\sigma$. Let $g_0\in BS(1,2)$ be a group element such that $x_{g_0}\in t_1 \cup t_2$. Then there exists a unique infinite path $\left(q_i,w_i\right)_{i\in\N}$ in the transducer $\mathcal{D}$ such that
$$x_{g_0\cdot\prod_{i=0}^n w_i}\in q_n$$
for every $n\in\N$.
\end{proposition}

\begin{proof}
 The proof follows directly from a careful observation of the local rules, and more particularly local rule~\ref{rule.matching_rules} (see Page~\pageref{rule.matching_rules}).
\end{proof}

Let $x\in X_\sigma$ and $g_0\in BS(1,2)$ such that $x_{g_0}\in t_1\cup t_2 \cup b_1 \cup m_1 \cup m_2$. Assume that $g_0$ has normal form $g_0=t^ia^jt^{-k}$. Thanks to Proposition~\ref{proposition.path_transducer}, there exists an infinite path $\left(q_i,w_i\right)_{i\in\N}$ in the transducer $\mathcal{D}$ such that
$$x_{g_0\cdot\prod_{i=0}^n w_i}\in q_n$$
for every $n\in\N$ and this path is unique. Define 
$$v(x_{g_0}):=j \cdot 2^{-i} + 2^{k-i} \cdot \left(\sum_{\ell=0}^{\infty} \delta_{w_\ell}^{ta}\cdot 2^{-\ell-1} \right) \in \left[j\cdot 2^{-i};(j+1)\cdot 2^{k-i}\right]$$
to be the position of the \define{vertical line initialized at ${g_0}$ in $x$}, where $\delta_{q}^{r}$ is the Kronecker symbol that equals $1$ if $q=r$ and $0$ otherwise.

\begin{proposition}\label{proposition.vertical_position_well_defined}
For a given configuration $x\in X_\sigma$ and group element $g_0\in BS(1,2)$, define the finite word $h\in \left\{ ta,t\right\}^M$ as the label of the $M$ first transition of the unique path given by Proposition ~\ref{proposition.path_transducer}. Then 
$$v(x_{g_0})= v(x_{g_0\cdot h}).$$	
\end{proposition}

\begin{proof}
First note that if $\left(q_i,w_i\right)_{i\in\N}$ is the unique infinite path in the transducer $\mathcal{D}$ for $g_0$, then $\left(q_{i+M},w_{i+M}\right)_{i\in\N}$ is the unique infinite path in the transducer $\mathcal{D}$ for $g_0\cdot h$. Then, using the normal form of $g_0$ and the fact that $t^{-k}a=a^{2^k}t^{-k}$, we deduce the normal form for $g_0\cdot h$.

\begin{align*} 
g_0\cdot h &=  t^i \cdot a^j \cdot t^{-k}\cdot h_0\dots h_{M-1} \\ 
 &= t^i \cdot a^{j+2^{k-1}\cdot \delta_{h_0}^{ta}} \cdot t^{-(k-1)}\cdot h_1\dots h_{M-1} \\
 &= t^i \cdot a^{j+2^{k-1}\cdot \delta_{h_0}^{ta}+ 2^{k-2}\cdot \delta_{h_1}^{ta}} \cdot t^{-(k-2)}\cdot h_2\dots h_{M-1} \\
 &= t^i\cdot a^{j+\sum_{\ell=0}^{M-1} 2^{k-\ell-1}\cdot \delta_{h_\ell}^{ta}}\cdot t^{-(k-M)}
\end{align*}

There are now two cases. If $M-k\leq 0$, then we have the normal form for $g_0\cdot h$. By definition, 
\begin{align*}
v(g_0\cdot h) &= \left( j+\sum_{\ell=0}^{M-1} 2^{k-\ell-1}\cdot \delta_{h_\ell}^{ta}\right)\cdot 2^{-i} + 2^{k-M-i}\cdot \sum_{\ell=0}^{\infty} \delta^{ta}_{w_{M+\ell}}\cdot 2^{-\ell-1}\\
&= j\cdot 2^{-i} + 2^{k-i}\left( \sum_{\ell=0}^{M-1} 2^{-\ell-1}\cdot \delta_{h_\ell}^{ta}+ \sum_{\ell=0}^{\infty} \delta^{ta}_{w_{M+\ell}}\cdot 2^{-M-\ell-1}\right) \\
&= j\cdot 2^{-i} + 2^{k-i}\left( \sum_{\ell=0}^{M-1} 2^{-\ell-1}\cdot \delta_{h_\ell}^{ta}+ \sum_{\ell=M}^{\infty} \delta^{ta}_{w_{\ell}}\cdot 2^{-\ell-1}\right)\\
&= v(g_0)		 
\end{align*}
and we are done.

If $M-k\geq 0$, then using the fact that $at^k=t^ka^{2^k}$ we get that the normal form for $g_0\cdot h$ is this time

$$g_0\cdot h = t^{i+M-k}\cdot a^{2^{M-k}\left( j+\sum_{\ell=0}^{M-1} 2^{k-\ell-1}\cdot \delta_{h_\ell}^{ta} \right)}$$

and thus 

\begin{align*}
v(g_0\cdot h) &= 2^{M-k}\left( j+\sum_{\ell=0}^{M-1} 2^{k-\ell-1}\cdot \delta_{h_\ell}^{ta} \right)\cdot 2^{-i-M+k} + 2^{-i-M+k} \cdot \sum_{\ell=0}^{\infty} \delta^{ta}_{w_{M+\ell}}\cdot 2^{-\ell-1}\\
&= j\cdot 2^{-i}+ 2^{k-i}\left( \sum_{\ell=0}^{M-1} 2^{-\ell-1}\cdot \delta_{h_\ell}^{ta}+ \sum_{\ell=M}^{\infty} \delta^{ta}_{w_{\ell}}\cdot 2^{-\ell-1}\right)\\
&= v(g_0)		 
\end{align*}
and the proposition is proven.

\end{proof}

\begin{proof}[Proof of Lemma~\ref{lemma.X_sigma_to_tiling}]Let $x$ be a configuration in $X_\sigma$. We proceed in two steps
\begin{enumerate}
 \item every sheet of $x$ defines a $\sigma$-tiling;
 \item the $\sigma$-tiling defined by one sheet of $x$ is compatible with all the sheets.
\end{enumerate}

\paragraph{Step 1: every sheet of $x$ defines a $\sigma$-tiling}

Consider one sheet in $x$. Our strategy to define $\tau$ a $\sigma$-tiling is the following : we first partially describe a tiling through all its vertical lines, then using the fact we want a $\sigma$-tiling, we define horizontal lines and the type -- i.e. the letter $a\in\A$-- of every $\sigma$-tile.
\begin{enumerate}
 \item Thanks to Proposition~\ref{proposition.vertical_position_well_defined}, every top tile $t_1(\dots,i\geq2)$ or $t_2(\dots,i\geq2)$ in position $g$ defines a semi-infinite vertical line : its horizontal position is completely determined and equals $v(x_g)$, while its vertical position is only known up to some error.
 \item Now that all the horizontal positions of semi-infinite vertical lines are known, we deduce the exact position of all tiles composing $\tau$ from the definition of $\sigma$-tiles (see Figure~\ref{figure.example_sigma_tiles}). Note that the type of tile $a\in\A$ is given by the letter $t_1$ or $t_2$ on the top left corner. An $a$-tile of with left and right borders at horizontal positions $x_\ell$ and $x_r$ respectively should be at position $\left( x_\ell,\log_2\left( \frac{x_r-x_\ell}{v(a)}\right)\right)$.  
\end{enumerate}
With these two conditions together, the tiling $\tau$ we get is indeed a $\sigma$-tiling.

\paragraph{Step 2: all sheets define the same $\sigma$-tiling}
We need to check that two different sheets give the same $\sigma$-tiling under the construction described immediately above. This is actually ensured by local rule~\ref{rule.synchro_two_sheets} on page~\pageref{rule.synchro_two_sheets}: this condition forces left border of shaded area to be synchronized between merging sheets, hence the horizontal positions of vertical lines do not depend on the sheet chosen.

\end{proof}

\begin{proposition}\label{proposition.X_sigma_SFT_cover_Y_sigma}
  The subshift $X_\sigma$ is an SFT cover of $Y_\sigma$, and this latter is consequently a sofic subshift.
 \end{proposition}
 
\begin{proof}
 This is a direct consequence of Theorem~\ref{theorem.sigma_tilings_and_X_sigma}: the onto morphism from $X_\sigma$ to $Y_\sigma$ is the letter-to-letter morphism given by the local map $\pi$.
\end{proof}

\section{Dynamical properties of $X_\sigma$}
\label{section.dynamical_properties_X_sigma}

\begin{proposition}\label{proposition.X_sigma_minimal}
  The subshift $X_\sigma$ is minimal if and only if $\log_2(\lambda)$ is irrational.
 \end{proposition}
 
 \begin{proof}
 Assume that $\log_2(\lambda)$ is irrational. Let $x,y$ be two configurations in $X_\sigma$. We prove that the orbit of $x$ is dense in $X_\sigma$, i.e. that for every $\varepsilon>0$, there exists some $g\in BS(1,2)$ such that $\mathfrak{S}_g(x)$ is close enough to $y$, i.e. that $d(\mathfrak{S}_g(x),y)<\varepsilon$. We fix $\varepsilon=2^{-N}$ and we proceed in two steps. First we use that the set $\left\{ \lceil k\cdot\log_2(\lambda)\rceil - k\cdot\log_2(\lambda)\mid k\in\Z\right\}$ is dense in $[0;1]$, since $\log_2(\lambda)$ has be assumed irrational. Thus by a shift $t^k$, we can synchronize the tilings $\tau_{\mathfrak{S}_{t^k}(x)}$ and $\tau_{y}$ so that their sequences of height $h$ and height $h+1$ rows coincide at least on the ball of radius $N$. Second we use that since the substitution $\sigma$ is primitive, then the subshift $Z_\sigma$ is minimal~\cite{}. With a shift $t^na^\ell t^{-n}$ we can perform an horizontal shift of $\ell 2^{-n}$ on $\tau_{\mathfrak{S}_{t^k}(x)}$ so that $d(\mathfrak{S}_{t^na^\ell t^{-n} t^k}(x),y)<\varepsilon$ by minimality of $Z_\sigma$ as required. Thus $X_\sigma$ is minimal.

 \medskip
 
 Conversely, assume $\log_2(\lambda)$ is rational and can be written as $\log_2(\lambda)=\frac{p}{q}$ with $p,q\in\N^*$ and $q\nmid p$ (unless $q=1$). Then with $q$ rows of $\sigma$-tiles, we get a stripe of height $p\in\N^*$. Assume that the top of this stripe has integer vertical coordinate $y\in\Z$. Then all other stripes of $q$ rows of $\sigma$-tiles have integer coordinates. Denote $\tau$ this $\sigma$-tiling, and $x\in X_\sigma$ the corresponding configuration obtained as in the proof of Lemma~\ref{lemma.tiling_to_X_sigma}. If we consider $\widetilde{h}:=\left(h_i\right)_{i\in\Z}$ the sequence of $h$' and $h+1$' where $h_i$ is the number of $\Phi$-boxes vertically intersected by a $\sigma$-tile on the $i$th row of $\tau$, we get a periodic sequence with period $H=h_1\dots h_q$. And so $h$'s and $h+1$' have frequencies $f_h=\frac{|\left\{i\in[1;q]\mid h_i=h\right\}|}{q}$ and $f_{h+1}=1-f_h$ in $\widetilde{h}$. Also since we have chosen that all stripes of $q$ rows of $\sigma$-tiles have integer coordinates, we deduce that both $h_1=h$ and $h_q=h$.
 
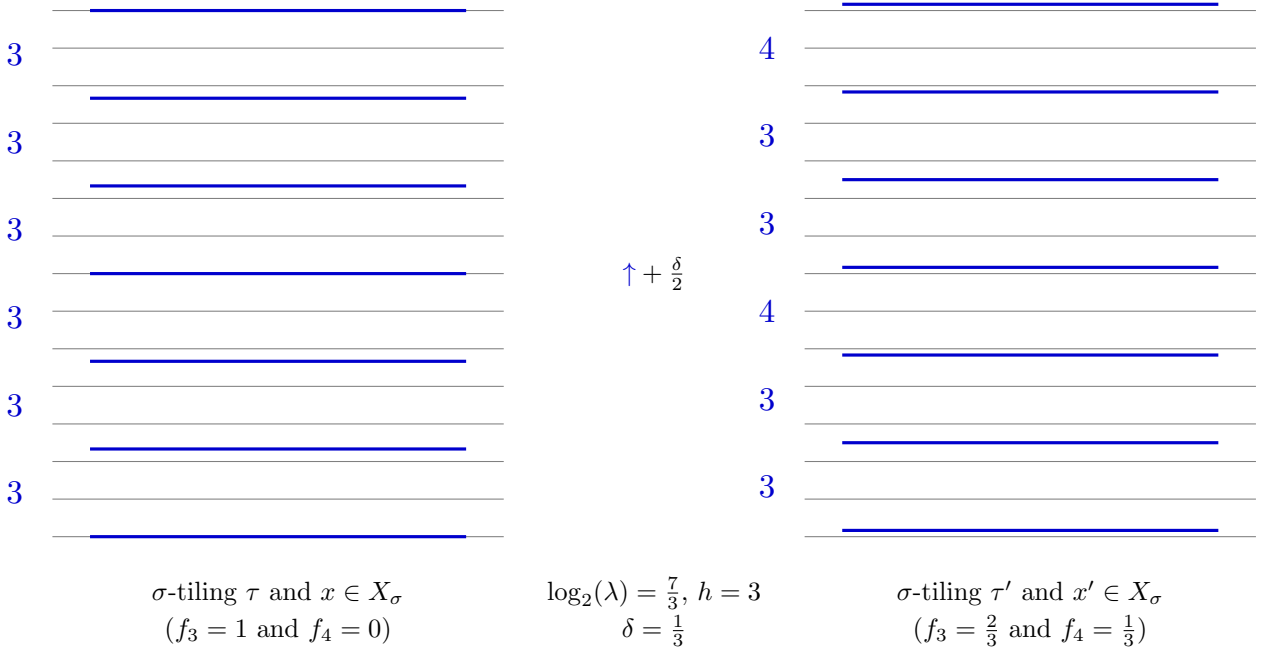
\begin{figure}[!ht]
\begin{bigcenter}
\begin{tikzpicture}[scale=0.5]

\begin{scope}[shift={(0,0)}]
\foreach \i in {-7,...,7}{
	\draw[color=black!50] (-6,\i) -- (6,\i);
	} 
	
\foreach \i in {-3,...,3}{
	\draw[very thick,color=bleu] (-5,7*\i/3) -- (5,7*\i/3);
	}
\node[color=bleu,scale=1.25] at (-7,35/6) {$3$};	
\node[color=bleu,scale=1.25] at (-7,21/6) {$3$};
\node[color=bleu,scale=1.25] at (-7,7/6) {$3$};
\node[color=bleu,scale=1.25] at (-7,-7/6) {$3$};
\node[color=bleu,scale=1.25] at (-7,-21/6) {$3$};
\node[color=bleu,scale=1.25] at (-7,-35/6) {$3$};

\node[] at (0,-8.5) {$\sigma$-tiling $\tau$ and $x\in X_\sigma$};
\node[] at (0,-9.5) {($f_3=1$ and $f_4=0$)};
\end{scope}

\begin{scope}[shift={(20,0)}]
\foreach \i in {-7,...,7}{
	\draw[color=black!50] (-6,\i) -- (6,\i);
	} 

\foreach \i in {-3,...,3}{
	\draw[very thick,color=bleu] (-5,7*\i/3+1/6) -- (5,7*\i/3+1/6);
	}	
\node[color=bleu,scale=1.25] at (-7,35/6+1/6) {$4$};	
\node[color=bleu,scale=1.25] at (-7,21/6+1/6) {$3$};
\node[color=bleu,scale=1.25] at (-7,7/6+1/6) {$3$};
\node[color=bleu,scale=1.25] at (-7,-7/6+1/6) {$4$};
\node[color=bleu,scale=1.25] at (-7,-21/6+1/6) {$3$};
\node[color=bleu,scale=1.25] at (-7,-35/6+1/6) {$3$};

\node[] at (0,-8.5) {$\sigma$-tiling $\tau'$ and $x'\in X_\sigma$};
\node[] at (0,-9.5) {($f_3=\frac{2}{3}$ and $f_4=\frac{1}{3}$)};
\end{scope}
 
\node[] at (10,-8.5) {$\log_2(\lambda)=\frac{7}{3}$, $h=3$};  
\node[] at (10,-9.5) {$\delta=\frac{1}{3}$}; 
\node[] at (10,0) {$\textcolor{bleu}{\uparrow}+\frac{\delta}{2}$}; 
 
\end{tikzpicture}
\end{bigcenter}
\caption{Illustration of the proof that if $\log_2(\lambda)$ is rational, then $X_\sigma$ is not minimal: two $\sigma$-tilings $\tau$ and $\tau'$ that give two configurations $x$ and $x'$ both in $X_\sigma$ but not in the same orbit.}
\label{figure.proof_loglambda_rational_implies_non_minimal}
\end{figure}
 
 Define $\delta:= \min \left\{d(i\cdot\frac{p}{q},\N)\mid i=1\dots q-1\right\}$. Then necessarily $\delta>0$, otherwise $q\nmid p$. Consider $\tau'$ the $\sigma$-tiling obtained by vertically shifting $\tau$ by $+\frac{\delta}{2}$, and denote $x'$ the associated configuration in $X_\sigma$ and $\widetilde{h}':=\left(h'_i\right)_{i\in\Z}$ the sequence of $h$' and $h+1$' where $h'_i$ is the number of $\Phi$-boxes vertically intersected by a $\sigma$-tile on the $i$th row of $\tau'$. Similarly to $\widetilde{h}$, the sequence $\widetilde{h}'$ is periodic with period $q$, but with $h'_1=h+1$. So the frequencies of $h$' and $h+1$' in $\widetilde{h}$ and $\widetilde{h}'$ are strictly different, thus $x'$ cannot be in the orbit of $x$. And finally the subshift $X_\sigma$ is not minimal, which terminates the proof.

 \end{proof}

\begin{proposition}\label{proposition.Y_sigma_aperiodic}
  The subshift $Y_\sigma$ is strongly aperiodic if $\log_2(\lambda)$ is irrational.
\end{proposition}
 
 \begin{proof} 
 Let $y\in Y_\sigma$, and assume that $h$ is a period for $y$ : for every $g\in BS(1,2)$, one has $y_g=y_{h^{-1}g}$. 
 
 If $h^{-1}$ has normal form $h^{-1}=t^\ell a^m t^{-p}$ and $gat$ has normal form $gat=t^i a^j t^{-k}$, the fact that $y_g=y_{h^{-1}g}$ is equivalent to
 \[
 \tau(j\cdot 2^{-i},k-i)=\tau(\left(j+m\cdot2^{i-p}\right)2^{p-\ell-i},p-\ell+k-i)
 \]
 for a certain $\sigma$-tiling $\tau$. Fix $i=j=0$ in this equality. Then for every $k\in\N$, one has
\[
 \tau(0,k)=\tau(m\cdot2^{p-\ell},p-\ell+k).
 \]
Since $\log_2(\lambda)$ is irrational, there exists some $k_0\in\N$ such that $\tau(0,k_0)$ and $\tau(m\cdot2^{p-\ell},p-\ell+k_0)$ do not belong to the same horizontal row of $\sigma$-tiles. Thus two different rows of $\sigma$-tiles in $\tau$ are the same, up to an horizontal translation by $m\cdot2^{p-\ell}$. This is possible only if the translation by $\left(m\cdot2^{p-\ell},p-\ell\right)$ is trivial, i.e. $h=1$. Finally every configuration $y\in Y_\sigma$ has a trivial stabilizer, so the subshift $Y_\sigma$ is strongly aperiodic.
 \end{proof}
 
As an immediate consequence of Proposition~\ref{proposition.Y_sigma_aperiodic} and Proposition~\ref{proposition.X_sigma_SFT_cover_Y_sigma} we get that the SFT $X_\sigma$ is also strongly aperiodic, provided $\log_2(\lambda)$ is irrational. 

\begin{proposition}\label{proposition.X_sigma_aperiodic}
The SFT $X_\sigma$ is strongly aperiodic if $\log_2(\lambda)$ is irrational.
\end{proposition}

Thanks to Lemma~\ref{lemma.bottom_forces_rectangle}, we can also prove that the SFT $X_\sigma$ has zero entropy.

\begin{proposition}\label{proposition.X_zero_entropy}
  The SFT $X_\sigma$ has zero entropy.
\end{proposition}
 
 \begin{proof}
 Choose $k$ and $\ell$ such that $\ell>>h$. Recall that 
 $$R_{k,\ell}:=\left\{ t^{\ell}a^it^{-j} \mid i\in [0;(k+1)\cdot2^{\ell-1}-1]\text{ and }j\in[0;\ell]\right\}.$$
 We define three subsets of $R_{k,\ell}$ :
 \begin{align*}
  \texttt{Bottom}&:= \left\{ t^{\ell}a^it^{-j} \mid i\in [0;(k+1)\cdot2^{\ell-1}-1]\text{ and }j\in[0;k]\right\};\\
  \texttt{Left}&:= \left\{ t^{\ell}a^it^{-j} \mid i\in [0;\cdot2^{\ell-1}-1]\text{ and }j\in[0;\ell]\right\};\\
  \texttt{Right}&:= \left\{ t^{\ell}a^it^{-j} \mid i\in [k\cdot2^{\ell-1};(k+1)\cdot2^{\ell-1}-1]\text{ and }j\in[0;\ell]\right\};\\
  \texttt{Top}&:= \left\{ t^{\ell}a^it^{-j} \mid i\in [0;(k+1)\cdot2^{\ell-1}-1]\text{ and }j\in[\ell-k;\ell]\right\}.
 \end{align*}
Denote $N:=\sharp\texttt{Bottom}+\sharp\texttt{Left}+\sharp\texttt{Right}+\sharp\texttt{Top}$. Then
\[
 N=2(k+1)^2\cdot2^{\ell-1}+2\ell\cdot 2^{\ell-1}+,
\]
and since $|R_{k,\ell}|=\ell\cdot(k+1)\cdot2^{\ell-1}$, we get that $\lim_{k,\ell\to\infty}\frac{N}{|R_{k,\ell}|}=0$. 

\medskip

Consider now a pattern $p$ with support $R_{k,\ell}$ that is in the language of $X_\sigma$. Then by Lemma~\ref{lemma.bottom_forces_rectangle} and Local rule~\ref{rule.substitution}, the pattern $p$ is entirely determined by $p_{\texttt{Bottom}\cup\texttt{Left}\cup\texttt{Right}\cup\texttt{Top}}$. Hence $\sharp \mathcal{L}_{k,\ell}(X_\sigma)\leq |\mathcal{R}|^{N}$ and consequently
\begin{align*}
h(X_\sigma)&=\lim_{k,\ell\to\infty} \frac{\log\left(\sharp \mathcal{L}_{k,\ell}(X_\sigma)\right)}{|R_{k,\ell}|}\\
&\leq \lim_{k,\ell\to\infty} |\mathcal{R}| \frac{N}{|R_{k,\ell}|}\\
&=0.
\end{align*}
So the SFT $X_\sigma$ has zero entropy.
 \end{proof}

%

%
%
%
  
\section{An application: a hierarchical strongly aperiodic SFT on $BS(1,2)$}
\label{section.robinson_BS}
  
    \subsection{Robinson hyperbolic tileset}
    \label{subsection.tuiles_Robinson_hyperboliques}

\begin{figure}[!ht]
    \begin{center}
    \includegraphics[width=0.99\linewidth]{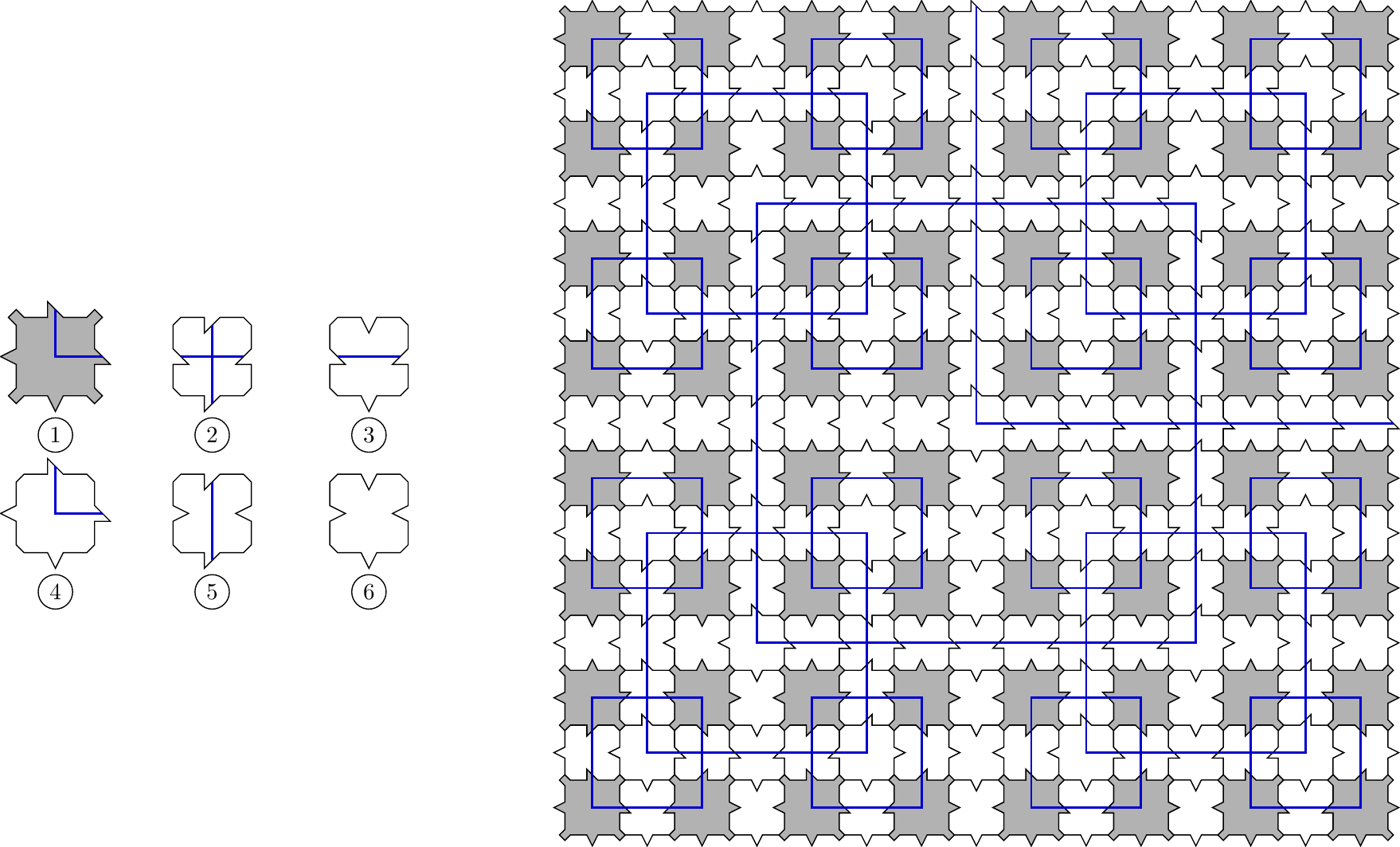}
    \end{center}
    \caption{On the left the six tiles from which one gets the whole Robinson's tileset, thanks to rotations and reflections. On the right a pattern made with Robinson's tileset.}
    \label{figure.Robinson_Z2}
    \end{figure}

    In the classical construction of a strongly aperiodic SFT on $\Z^2$~\cite{Robinson1971}, all valid tilings present a hierarchical structure of squares of increasing size, where the same process is repeated to obtain bigger and bigger squares: four squares of the same size are gathered to form a bigger square. The smallest squares are enforced by bumpy tiles (tiles~1 on Figure~\ref{figure.Robinson_Z2}).
    
    \medskip
    
\begin{figure}[!ht]
    \begin{bigcenter}
    \includegraphics[width=0.8\linewidth]{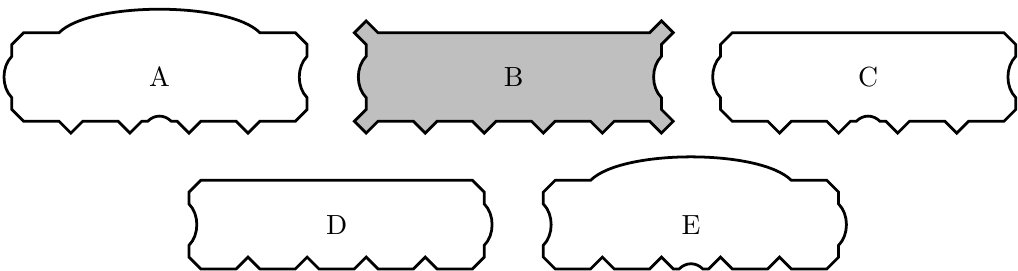}
    \end{bigcenter}
    \caption{The five shapes of tiles A, B, C, D and E on $\H_2$. Tile B is the only one with bumpy corners, tiles A and C have both bumpy and dented corners while tiles D and E have only dented corners.}
    \label{figure.bumpy_dented_H5_tiles}
    \end{figure}

\begin{figure}[!ht]
    \begin{bigcenter}
    \begin{tikzpicture}[xscale=0.05,yscale=0.25]
 

\tile{-10}{6}{50}{white}
\tile{-10}{3}{10}{white}
\tile{40}{3}{10}{black!30}
\tile{90}{3}{10}{white}
\tile{140}{3}{10}{black!30}
\tile{190}{3}{10}{white}
\foreach \x in {0,...,24} {
\tile{-10+\x*10}{0}{2}{white}
}
\foreach \x in {0,2,...,124} {
\tile{-10+\x*2}{-3}{0.4}{white}
}
\foreach \x in {1,3,...,123} {
\tile{-10+\x*2}{-3}{0.4}{black!30}
}

\draw[very thick, bleu] (65,4.5) rectangle (165,-1.5); 

\end{tikzpicture}
    \end{bigcenter}
    \caption{Location of the B  tiles (in grey) in a tiling. In blue, the expected shape of a square of level~1.}
    \label{figure.where_are_B}
    \end{figure}
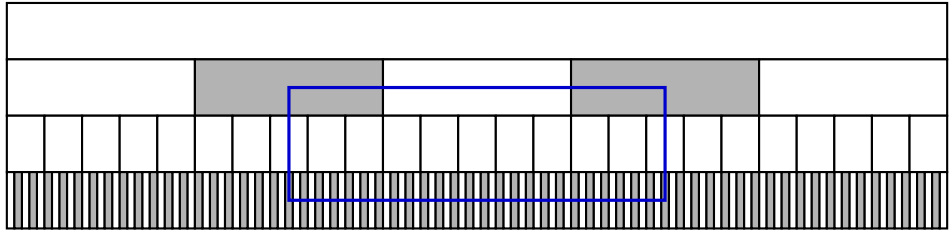

    The situation is more complex on $\H_2$, and instead of bumpy tiles and others, we use the five shapes pictured on Figure~\ref{figure.bumpy_dented_H5_tiles}. 
    These preliminary local rules already impose strong constraints of valid tilings. The global structure of A/B/C/D/E tilings appears on Figure~\ref{figure.where_are_B}. Indeed, first remark that decorations on left and right sides of the tiles impose that horizontal rows are either composed with A/B/C or with D/E. Then a careful study of the corners show an alternation of A/B/C rows and D/E rows, and that inside an A/B/C row, B tiles necessarily appear one time in two.

    \medskip
    
    On Figure~\ref{figure.where_are_B} is pictured what would be an hyperbolic analogue of a square of level~1 in the Robinson tiling on $\Z^2$: on the same horizontal row, two nearest neighbor tiles B carry the top left and top right corners of the square; the bottom corners are located two rows below, but the geometry of $\H_2$ imposes that, between these two corners carried by tiles B, we find $5^2-1=24$ other tiles B. Figure~\ref{figure.where_are_B} show that, contrary to the $\Z^2$ case, not all the B tiles carry a corner of a square. Thus we will need many more tiles than on $\Z^2$.
    
    \medskip  
    
    Before we go any further in the description of tiles of the hyperbolic Robinson tiling, one can already notice that constraints imposed by the five shapes of tiles are strong, and almost deterministic in some sense. Denote $\sigma$ the non-deterministic substitution on the alphabet $:\{$A,B,C,D,E$\}$ given by the rules
    \begin{align*}
    \sigma(\text{A}) &= \text{DDEDD}\\
    \sigma(\text{B}) &= \text{DDDDD}\\
    \sigma(\text{C}) &= \text{DDEDD}\\
    \sigma(\text{D}) &= \text{BCBCB or CBCBC}\\
    \sigma(\text{E}) &= \text{CBABC}.
    \end{align*}
    
    A careful comparison between the orbits of this substitution $\sigma$ and the A/B/C/D/E tilings shows that an orbit of $\sigma$ can be seen as an A/B/C/D/E tiling, and vice-versa.

    \medskip
    
    \begin{proposition}\label{proposition.ABCDE_sigma}
    To every tiling of $\H_2$ with tiles A/B/C/D/E, one can associate the orbit of a bi-infinite word $w\in\left\{A,B,C,D,E\right\}^\Z$ under the action of $\sigma$.
    \end{proposition}

    If we add extra decorations to the tiles A/B/C/D/E, we get new tilings of~$\H_2$. Thanks to Proposition~\ref{proposition.ABCDE_sigma}, one can think about these tilings as SFTs on orbit graphs of $\sigma$, so that the construction of Section~\ref{section.X_sigma_SFT_cover_Y_sigma} can be used.

    \medskip
    
\begin{figure}[!ht]
    \begin{bigcenter}
    \includegraphics[width=0.99\linewidth]{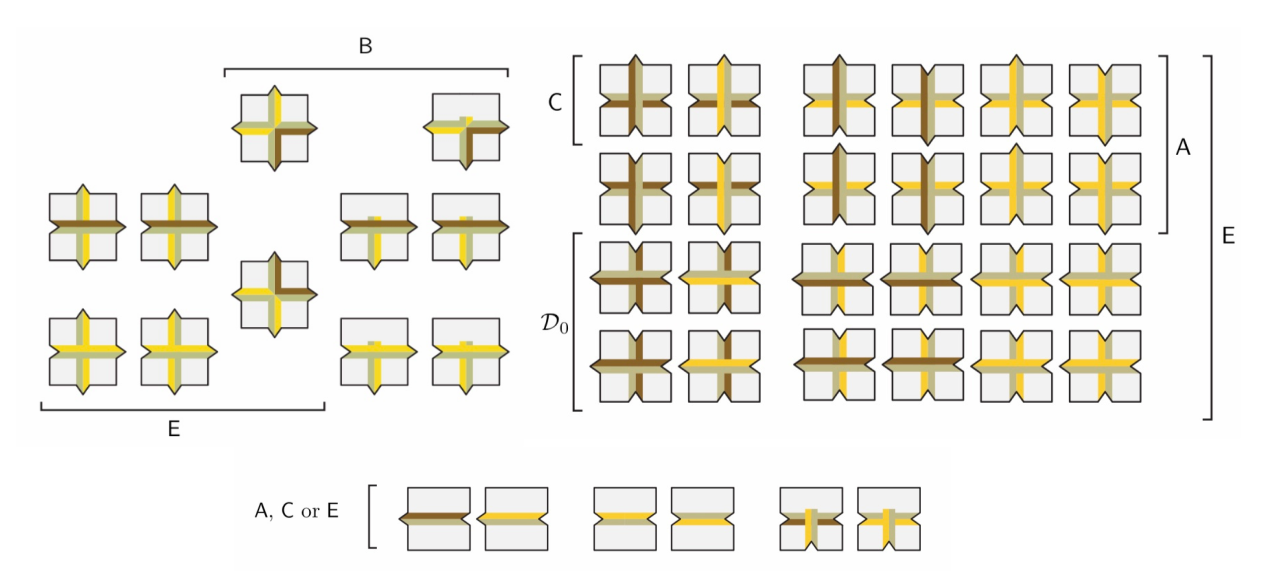}
    \end{bigcenter}
    \caption{The additional decorations on tiles A/B/C/D/E to get the hyperbolic Robinson tileset, as they are presented in~\cite{GS2010} and reproduced by kind permission of the author.}
    \label{figure.regles_GS_hyperbolique}
    \end{figure}

    We now enrich the basic tileset $\left\{ A,B,C,D,E\right\}$ by adding decorations pictured on Figure~\ref{figure.regles_GS_hyperbolique} so that the description of the hyperbolic Robinson tileset is finalized.

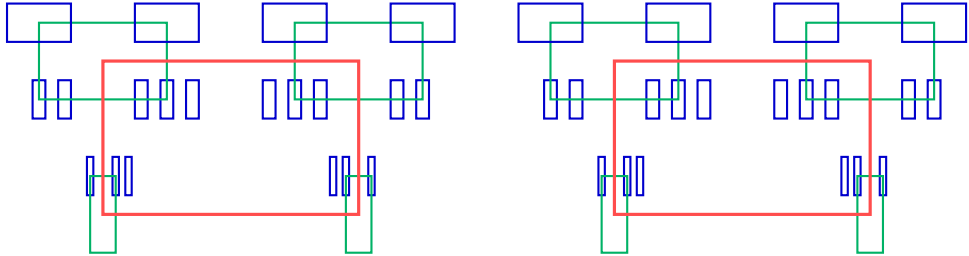
\begin{figure}[!ht]
    \begin{center}
    \begin{tikzpicture}[scale=0.17]

\foreach \x in {-8,-6,0,2,4,10,12,14,20,22}{
\draw[thick, bleu] (\x,0) rectangle (\x+1,-3); 
}

\foreach \x in {-3.75,-1.75,-0.75,15.25,16.25,18.25} {
\draw[thick, bleu] (\x,-6) rectangle (\x+0.5,-9);
} 

\draw[thick, vert] (2.5,-1.5) rectangle (-7.5,4.5); 
\draw[thick, vert] (22.5,-1.5) rectangle (12.5,4.5); 
\draw[thick, vert] (-3.5,-13.5) rectangle (-1.5,-7.5); 
\draw[thick, vert] (16.5,-13.5) rectangle (18.5,-7.5); 

\draw[thick, bleu] (-10,3) rectangle (-5,6);
\draw[thick, bleu] (0,3) rectangle (5,6);
\draw[thick, bleu] (10,3) rectangle (15,6);
\draw[thick, bleu] (20,3) rectangle (25,6);

\draw[very thick, rouge] (-2.5,1.5) rectangle (17.5,-10.5);

\begin{scope}[shift={(40,0)}]
\foreach \x in {-8,-6,0,2,4,10,12,14,20,22}{
\draw[thick, bleu] (\x,0) rectangle (\x+1,-3); 
}

\foreach \x in {-3.75,-1.75,-0.75,15.25,16.25,18.25} {
\draw[thick, bleu] (\x,-6) rectangle (\x+0.5,-9);
} 

\draw[thick, vert] (2.5,-1.5) rectangle (-7.5,4.5); 
\draw[thick, vert] (22.5,-1.5) rectangle (12.5,4.5); 
\draw[thick, vert] (-3.5,-13.5) rectangle (-1.5,-7.5); 
\draw[thick, vert] (16.5,-13.5) rectangle (18.5,-7.5); 

\draw[thick, bleu] (-10,3) rectangle (-5,6);
\draw[thick, bleu] (0,3) rectangle (5,6);
\draw[thick, bleu] (10,3) rectangle (15,6);
\draw[thick, bleu] (20,3) rectangle (25,6);

\draw[very thick, rouge] (-2.5,1.5) rectangle (17.5,-10.5); 
\end{scope}

\end{tikzpicture}
    \end{center}
    \caption{Schematic view of the first levels of the hierarchy of squares enforced by the hyperbolic Robinson's tileset. hyperbolic squares of level~1,~2 and~3 are pictured in blue, green and red. The structure that appears is an hyperbolic analogue of the one visible on Figure~\ref{figure.Robinson_Z2}.}
    \label{figure.first_level_hierarchy_H5}
    \end{figure}

    With the enriched tileset obtained by adding decorations from Figure~\ref{figure.regles_GS_hyperbolique} to the A/B/C/D/E tiles, we can tile the discrete hyperbolic plane $\H_2$, and the structure of tilings is schematized on Figure~\ref{figure.first_level_hierarchy_H5}.
    
     To obtain a tileset analogue with Robinson's and Goodman-Strauss's on every $BS(1,n)$ with $n\geq2$, it suffices to apply Theorem~\ref{theorem.sigma_tilings_and_X_sigma} to the $\sigma$-tilings with $\sigma$ the substitution on alphabet $\left\{ A,B,C,D,E\right\}$ defined above, and to enrich the SFT $X_\sigma$ with local rules so that the hyperbolic Robinson's tiling is copied out on every sheet of $BS(1,n)$.

    \begin{theorem}\label{SFT_aper_min_hierarchique}
     For every $n\geq2$ there exists a strongly aperiodic SFT on $BS(1,n)$ such that every configuration $x$ of the SFT carries an hyperbolic Robinson tiling on every single sheet of $x$.
    \end{theorem}

\section*{Acknowledgments}
The authors are thankfull to 
Chaim Goodman-Strauss for letting us use his pictures. This work was partially supported by the ECOS-SUD project C17E08 and the ANR project CoCoGro (ANR-16-CE40-0005).


\end{document}